\documentclass{report}

\usepackage{amsfonts,amsmath,eufrak}
\usepackage{amsmath}
\usepackage{amstext}
\usepackage{amsthm}
\usepackage{amssymb}
\usepackage{mathrsfs}
\usepackage{amsopn}
 \usepackage[dvips]{color}
 \usepackage{graphicx}
 \usepackage{latexsym}
 \usepackage{amsfonts}
 \usepackage{url}
 \usepackage{rotating}
 \usepackage[dvips]{color}
\usepackage{natbib}
\usepackage{subfigure}

\newtheorem{lemma}{Lemma}[chapter]
\newtheorem{theorem}[lemma]{Theorem}
\newtheorem{definition}[lemma]{Definition}
\newtheorem{corollary}[lemma]{Corollary}
\newtheorem{proposition}[lemma]{Proposition}

\newtheorem{example}[lemma]{Example}
\newtheorem{algorithm}[lemma]{Algorithm}
\newtheorem{Algorithm}[lemma]{Algorithm}
\newtheorem{question}[lemma]{Question}
\newtheorem{problem}[lemma]{Problem}

\newcommand{\maximize}{\mbox{maximize}}
\newcommand{\keyw}[1]{{\bf #1}}

\newcommand{\co}{\mbox{cone}}
\newcommand{\conv}{\mbox{conv}}
\newcommand{\ind}{\mbox{ind}}

\newcommand{\eopf}{\framebox[6.5pt]{} \vspace{0.2in}}

\newcommand{\real}{\mathbb R}
\newcommand{\jdlZ}{\mathbb Z}
\def\jdlqed{\vbox{\hrule \hbox{\vrule\hbox to
5pt{\vbox to 6pt{\vfil}\hfil}\vrule}\hrule}}
\newcommand{\Q}{{\mathbb Q}}
\newcommand{\Z}{{\mathbb Z}}
\newcommand{\N}{{\mathbb N}}

\newcommand{\C}{{\mathbb C}}
\newcommand{\R}{\mathbb R}

\newcommand{\be}{\begin{equation}\normalsize}
\newcommand{\ee}{\end{equation}\normalsizedb}
\title{Barvinok's Rational Functions: Algorithms and Applications to Optimization, Statistics, and Algebra}

\author{Ruriko Yoshida}

\date{\today}

\begin{document}
\pagenumbering{roman}
\begin{center}
{\Large \bf Barvinok's Rational Functions: Algorithms and Applications to Optimization, Statistics, and Algebra}
\vskip 0.2in
\centerline{By}
\vskip 0.1in
\centerline{\large  Ruriko Yoshida}
\centerline{\large B.A. (University of California at Berkeley) 2000}
\vskip 0.3in
\centerline{\large  DISSERTATION}
\vskip 0.1in
\centerline{\large Submitted in partial satisfaction of the requirements for the
degree of}
\vskip 0.1in
\centerline{\large  DOCTOR OF PHILOSOPHY}
\vskip 0.1in
\centerline{\large in}
\vskip 0.1in
\centerline{\large MATHEMATICS}
\vskip 0.2in
\centerline{\large in the}
\vskip 0.1in
\centerline{\large  OFFICE OF GRADUATE STUDIES}
\vskip 0.1in
\centerline{\large of the}
\vskip 0.1in
\centerline{\large  UNIVERSITY OF CALIFORNIA}
\vskip 0.1in
\centerline{\large  DAVIS}
\end{center}
\vskip 0.1in
\begin{center}
{\large Approved:}
\end{center}
\begin{center}
\centerline{\underbar{\hskip 2.5in}}
\vskip 0.15in
\centerline{\underbar{\hskip 2.5in}}
\vskip 0.15in
\centerline{\underbar{\hskip 2.5in}}
\vskip 0.2in
\centerline{\large Committee in Charge}
\vskip 0.2in
\centerline{\large 2004}
\end{center}

\pagebreak

\begin{center}

  \vspace*{50pt} 

To Pete Camagna,\\
for his love and support.
\end{center}

\newpage
\large
\tableofcontents
\listoftables
\listoffigures


\newpage
\large
{\Large \bf ACKNOWLEDGEMENTS} \\

First, I would like to thank my parents, Tasturou Yoshida and Chizuko Yoshida, in Japan for being patient with me.

Also I would like to thank my husband Pete Camagna for all support and for encouragement to write this thesis.  He supported me emotionally and gave me great inspiration to create some algorithms.  He has been always here for me no matter what happens.  

I would like to thank Prof Jes\'us A. De Loera for all his help as my thesis advisor.  He introduced me to Computational Algebra and Combinatorial generating functions in integer programming and Statistics.  Without him, I would not know these interesting topics.  He also motivated me to learn computer language {\tt C} and {\tt C++}, which I did not like when I was an undergraduate student.  I really love to implement algorithms in {\tt C++} now and create interesting algorithms to apply Statistics and Combinatorial optimization.  I really thank him for this.  I also wish to thank my thesis committee, Roger Wets and Naoki Saito for taking their time to read this thesis.  

I would also thank Raymond Hemmecke for useful conversations and encouragement to improve my thesis.  He was like my big brother, looking after me and giving me great advice. 

I would also like to thank Bernd Sturmfels.  Bernd Sturmfels gave me much great advice for improving the {\tt LattE} software, on which I was the head programmer.  {\tt LattE} is the only implementation to count the number of any rational convex polytope.  
Also when I was taking a phylogeny seminar from him in Berkeley, he gave me many great inspirations applying Algebraic algorithms to Biostatistics and Algebraic Statistics. Even though he is really busy, he looks after me and gives me great advises.  He is my greatest inspiration for Mathematics.  

Toward the end, Lior Pachter helped me on many aspects.  I am so glad that I have an opportunity to work with him.

Finally I would like to thank all my best friends, Alice Stevens and Peter Huggins, for their encouragement and support to finish this thesis.  Pete Camagna, Alice Stevens, and Peter Huggins are always here for me, every time I have problems and have been stressed out.  Without them I would not have finished this thesis.  

I am really lucky to have so many people who have been supportive. 

This thesis was partially supported by NSF Grants DMS-0309694,and DMS-0073815.


\newpage 
\begin{center}
\underline{\bf \Large Abstract} 
\end{center}

\medskip

\large
The main theme of this dissertation is the study of the lattice points in a rational convex polyhedron and their encoding in terms of Barvinok's short rational functions.  The first part of this thesis looks into theoretical applications of these rational functions to Optimization, Statistics, and Computational Algebra.
The main theorem on Chapter \ref{GBT} concerns the computation of the
\emph{toric ideal} $I_A$ of an integral $n \times d$ matrix $A$.
We encode the binomials belonging to the toric ideal 
$I_A$ associated with $A$ using Barvinok's rational functions.
If we fix $d$ and $n$, this
representation allows us to compute a universal Gr\"obner basis
and the reduced Gr\"obner basis of the ideal $I_A$, with respect
to any term order, in polynomial time. We
 also derive a polynomial time algorithm for normal form
computations which replaces in this new encoding the usual
reductions of the division algorithm.
Chapter \ref{IPT} presents three ways to use
Barvinok's rational functions to solve Integer Programs: The
\emph{$(A,b,c)$-test set algorithm}, the {\em Barvinok's binary search algorithm}, and the {\em digging algorithm}.
\vskip 0.1in

 The second part of the thesis is experimental and consists mainly of the software package {\tt LattE}, the first implementation of Barvinok's algorithm to compute short rational functions which encode the lattice points in a rational convex polytope.  
In Chapter \ref{chapsoftware} we report on
experiments with families of well-known rational polytopes: multiway
contingency tables, knapsack type problems, and rational polygons, and we present formulas for the Ehrhart quasi-polynomials of several
hypersimplices and truncations of cubes. 
We also developed a new algorithm, {\em the homogenized Barvinok's algorithm} to compute the generating function for a rational polytope. 
We showed that it runs in polynomial time in fixed dimension.
With the homogenized Barvinok's algorithm, we obtained new combinatorial formulas: the generating function for the number of $5\times
5$ magic squares and the generating function for the number of $3\times 3
\times 3 \times 3$ magic cubes as 
rational functions.


\newpage
\pagestyle{myheadings} 
\pagenumbering{arabic}
\markright{  \rm \normalsize CHAPTER 1. \hspace{0.5cm}
 Introduction}
\large 
\chapter{Introduction}\label{IntroThesis}
\thispagestyle{myheadings}



\section{Basic notations and definitions}\label{definitions}

In this section, we will recall basic notations and definitions.  
Let $\{x_1, x_2, \dots, x_m\}$ be a finite set of points in $\R^d$.  A point
$$x = \sum_{i=1}^m \alpha_i x_i, \mbox{ where } \sum_{i=1}^m \alpha_i = 1 \mbox{ and } \alpha_i \geq 0 \mbox{ for } i = 1, 2, \dots, m$$ is called a {\bf convex combination } of $x_1, x_2, \dots, x_m$.  Given two distinct points $x, \, y \in \R^d$, the set $[x, y] = \{\alpha x + (1-\alpha) y: 0 \leq \alpha \leq 1 \}$ of all convex combinations of $x$ and $y$ is called the {\bf interval} with endpoints $x$ and $y$.  A set $C \subset \R^d$ is called {\bf convex}, provided $[x, y] \subset C$ for any two $x, \, y \in C$.  For $C \subset \R^d$, the set of all convex combinations of points from $C$ is called the {\bf convex hull} of $C$ and denoted $\conv(C)$.
Let $A_1, A_2, \dots, A_n \in \R^d$ and let $b_1, b_2, \dots, b_n \in \R$.  Then the set $$P:= \{x \in \R^d: A_i \cdot x \leq b_i \mbox{ for } i = 1, 2, \dots, n\}$$ is called a {\bf polyhedron}.  The convex hull of a finite set of points in $\R^d$ is called a {\bf polytope} and the Weyl-Minkowski Theorem says that a polytope is a bounded polyhedron \citep{schrijver}.

A finite set of points $\{ x_1, x_2, \ldots, x_k\} \subseteq \R^d$ is 
{\bf affinely independent} if $\forall \lambda_j \in \R$, $\sum_{j=1}^k
\lambda_j x_j =0, \sum_{j=1}^k \lambda_j =0$ $\Rightarrow$ $\lambda_j =0$
$\forall j=1,2, \ldots , k$.

A $(d-1)$ dimensional affine set in $\R^d$ is called a {\bf hyperplane} and every hyperplane can be represented as
 $\{ x \in \R^d: ax=b$,
 $a \in \R^d$, $a \not = 0$, $b \in \R \}$.  $a$ is called a {\em
normal vector} of this hyperplane.

Let $H:= \{ x \in \R^d : h \cdot x \leq \beta\}$, where $h \in \R^d$, $h \not = 0$, and $\beta \in \R$, be an {\bf affine half space}.  Then if $P \subset H$ and $P \cap  \{ x \in \R^d : h \cdot x = \beta\} \not = \emptyset$, then $H$ is called a {\bf supporting hyperplane} of $P$.  A subset $F$ of $P$ is called a {\bf face} if $F = P$ or $F = P \cap H$, where $H$ is a supporting hyperplane.  If a face $F$ is minimal with respect to inclusion and $F$ contains only a point, then $F$ is called a {\bf vertex}.

Let $V(P)$ be the set of all vertices of $P$.  If any points in $V(P)$ are in $\Z^d$ then $P$ is called an {\bf integral} polyhedron.  If any points in $V(P)$ are in $\Q^d$ then $P$ is called a {\bf rational} polyhedron.

Now we will define our main tool, a {\bf generating function} of a polyhedron.
Let $P \subset \R^d$ be a polyhedron and let $\Z^d \subset \R^d$ be the integer lattice.  For an integral point $m = (m_1, m_2, \dots, m_d) \in \Z^d$, we can write the monomial $$z^m:= z^{m_1}_1z_2^{m_2}\dots z_d^{m_d}$$ in $d$ complex variables, $z_1, z_2, \dots, z_d$.  The generating function $f(P, z)$ of a polyhedron $P$ is the sum of monomials such that: 

\begin{equation} \label{generating}
f(P, z) = \sum_{m \in P \cap \Z^d} z^m.
\end{equation}

For example, consider the integral quadrilateral shown in Figure \ref{egbrion} with the
vertices $V_1=(0,0)$, $V_2=(5,0)$, $V_3=(4,2)$, and $V_4=(0,2)$.
Then we have the generating function $f(P, z)$ such that:

\begin{figure}[ht]
\begin{center}
\includegraphics[width=5 cm]{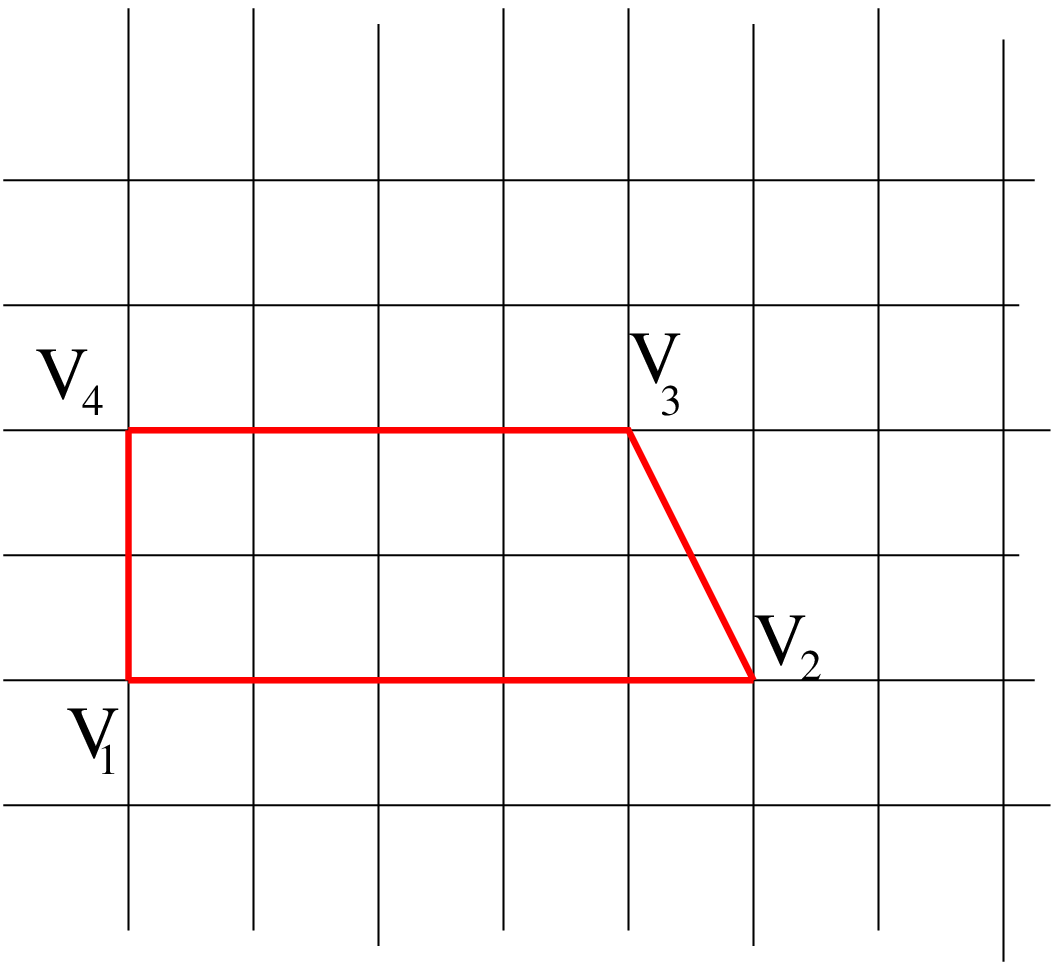} 
\caption{An example for the generating function}\label{egbrion}
\end{center} 
\end{figure}

$f(P, z) = {z_1}^{5}+{z_1}^{4}z_2+{z_1}^{4}+{z_1}^{4}{z_2}^{2}+z_2{z_1}^{3}+{z_1}^{3}+{z_1}^{3}{z_2}^{2
}+z_2{z_1}^{2}+{z_1}^{2}+{z_1}^{2}{z_2}^{2}+z_1z_2+z_1+z_1{z_2}^{2}+{z_2}^{2}+z_2+1$. 

One notices that the multivariate generating function $f(P, z)$ has exponentially many monomials even though we fixed the dimension.  So one might ask if it is possible to encode $f(P, z)$ in a ``short'' way. 
In 1994, A. Barvinok showed an algorithm that counts
the lattice points inside $P$ in  polynomial time when
$d$ is a constant \citep{bar}. The input for this algorithm
is the binary encoding of the integers $A_{ij}$ and $b_i$,
and the output is a short formula for the multivariate
generating function $f(P,\,  z)=\sum_{a \in P \cap \Z ^d} z^a$. This
long polynomial with exponentially many monomials  is
encoded as a short sum of rational functions in the form
\begin{equation}
\label{barvinokseries}
f(P, \, z) \quad
= \quad
\sum_{i \in I} \pm \frac{z^{u_i}}{(1-z^{c_{1,i}})(1-z^{c_{2,i}})\dots
(1-z^{c_{d,i}})},
\end{equation}

\noindent
where $u_i, \, c_{1, i}, \, c_{2, i}, \, \ldots ,\, c_{d, i} \in \Z^{d}$ and where $I$ is a polynomial sized index set.

\noindent
We call this short sum of rational functions of the form (\ref{barvinokseries}) {\bf Barvinok's rational function} for the generating function $f(P, z)$.  For brevity, we also call it a {\bf short rational function} for the generating function $f(P, z)$. 
For example, suppose we have the polytope in Figure \ref{egbrion}.  Then we can write:

\noindent
$f(P, z)={z_1}^{5}+{z_1}^{4}z_2+{z_1}^{4}+{z_1}^{4}{z_2}^{2}+z_2{z_1}^{3}+{z_1}^{3}+{z_1
}^{3}{z_2}^{2}+z_2{z_1}^{2}+{z_1}^{2}+{z_1}^{2}{z_2}^{2}+z_1z_2+z_1+z_1{z_2}^{2}+{z_2}^{2}+z_2+1$

\noindent
$= \frac {1}{\left (1-z_1\right )\left (1-z_2\right )} + \frac{{z_1}^{5}}{ (1-{z_1}^{-1}) (1-z_2)} + \frac{z_1^2}{(1-z_1)(1-z_2^{-1})} + \frac{z_1^5}{(1-z_1^{-1}z_2)(1-z_2^{-1})} +  \frac{{z_1}^{4}{z_2}^2}{ (1-{z_2}^{-1})
(1-z_1)} -\newline \frac{z_1^4z_2^{2}}{(1-z_1^{-1}{z_2}^2 )(1-z_1^{-1})}.$

\begin{figure}[ht]
\begin{center}
\includegraphics[width=7 cm]{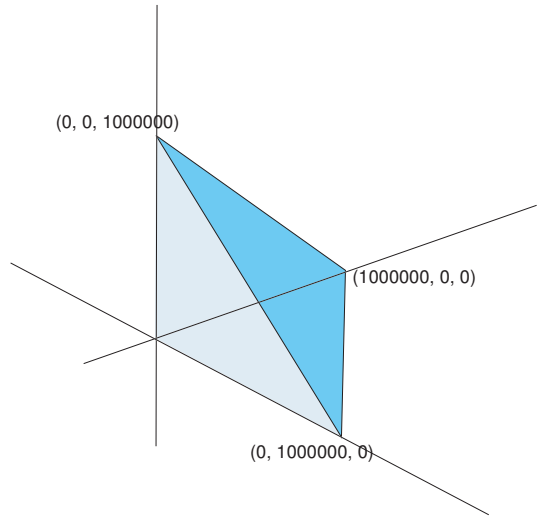} 
\caption{A tetrahedron example for the generating function}\label{tet}
\end{center} 
\end{figure}

Here is another example to clarify Barvinok's rational functions.  Suppose we have a tetrahedron $P$ with vertices $v_1 = (0, 0, 0)$, $v_2 = (1000000, 0, 0)$, $v_3 = (0, 1000000, 0)$, and $v_4 = (0, 0, 1000000)$ in Figure \ref{tet}.  Then we have the multivariate generating function $f(P, z)$ which has $166,667,666,668,500,001$ monomials.    However, if we use Barvinok's rational functions, we can represent all these monomials using a ``small'' encoding:

$$f(P, z)\quad =\quad \frac{1}{(1-z_1)(1-z_2)(1-z_3)}+ \frac{z_1^{1000000}}{(1-z_1^{-1}z_2)(1-z_1^{-1}z_3)(1-z_1^{-1})}$$

$$+\frac{z_2^{1000000}}{(1-z_1z_2^{-1})(1-z_2^{-1}z_3)(1-z_2^{-1})}+\frac{z_3^{1000000}}{(1-z_1z_3^{-1})(1-z_2z_3^{-1})(1-z_3^{-1})}.$$

One might ask how one can compute Barvinok's rational functions for the input polyhedron.  The following theorem tells us that there is an algorithm created by \cite{bar} to compute Barvinok's rational functions from the input polyhedron in polynomial time in fixed dimension.  We will describe a variation of the algorithm in Chapter \ref{chapsoftware}.

\begin{theorem}\label{decombar}\cite[Theorem 5.4]{bar}
Fix the dimension $d$.  Then there exists a polynomial time algorithm which for a given rational polyhedron $P \subset \R^d$, computes $f(P, z)$ in the form of  (\ref{barvinokseries}) in polynomial time.
\end{theorem}

We do not want to expand Barvinok's rational functions because expanding them causes exponential complexity.  So, if we want to perform operations on sets via generating functions, such as taking unions, intersections, projections, and complements, we want to do it directly with Barvinok's rational functions without expanding them.  The
{\bf Hadamard product} of Laurent power series is a very useful tool for Boolean operations on sets via Barvinok's rational functions.

\begin{definition}
 Let $g_1$ and $g_2$ be Laurent power series in $z \in \C^d$ such that $g_1(z) = \sum_{m \in \Z^d} a_{m} z^{m}$ and $g_2(z) = \sum_{m \in \Z^d} b_{m} z^{m}$.  Then the Hadamard product $g = g_1 * g_2$ is the power series such that:
$$ g(z) = \sum_{m \in \Z^d} a_{m} b_{m} z^m.$$
\end{definition}

Hadamard products of Laurent power series are one of the most important tools to prove theorems in this thesis. They are used for taking unions of sets, intersections of sets, and set difference via short rational functions without expanding them.  We will show how to compute the Hadamard product of Laurent power series $g_1$ and $g_2$ given in the form of rational functions.  Let $p_1, p_2, a_{11}, \dots, a_{1k} \in \Z^d$ and $a_{21}, \dots, a_{2k} \in \Z^d$. Suppose we are given the  Laurent power series $g_1$ and $g_2$ in the form:

\begin{equation}\label{LPS}
g_1 = \frac{z^{p_1}}{(1-z^{a_{11}}) \dots (1-z^{a_{1k}})} \mbox{ and } g_2 = \frac{z^{p_2}}{(1-z^{a_{21}}) \dots (1-z^{a_{2k}})}.
\end{equation}
 
Here is an outline of the algorithm to take the Hadamard product of two  Laurent power series via Barvinok's rational functions.

\begin{Algorithm}\label{HadamardAlg}\cite[Lemma 3.4]{newbar}

{\bf Input}: Laurent power series $g_1$ and $g_2$ in the form of (\ref{LPS}).

{\bf Output}: The Hadamard product $g_1 * g_2$ of $g_1$ and $g_2$ in the form of a rational function.

{\bf Step 1}: If $a_{1j} > 0$ or $a_{2j} > 0$, then apply the identity:

$$\frac{z^p_j}{1-z^a_j} = -\frac{z^{p-a}_j}{1-z^{-a}_j},$$
to reverse the direction of $a_{1j}$ or $a_{2j}$.

{\bf Step 2}: Let $P \subset \R^{2k} = \{(\xi_1, \xi_2, \dots, \xi_{2k}): \xi_i \in \R \mbox{ for } i = 1, 2, \dots, 2k\}$ be a rational polyhedron defined by the equation:

$$p_1 + \xi_1 a_{11} + \dots + \xi_k a_{1k} = p_2 + \xi_{k+1} a_{21} + \dots + \xi_{2k} a_{2k},$$

and inequalities $$ \xi_i \geq 0 \mbox{ for } i = 1, 2, \dots, 2k.$$

{\bf Step 3}: Using Theorem \ref{decombar}, we compute:
$$f(P, x) = \sum_{i \in I} \pm \frac{x^{u_i}}{(1-x^{v_{1,i}})(1-x^{v_{2,i}})\dots
(1-x^{v_{2k,i}})},$$ for $u_i, v_{i, j} \in \Z^{2k}$. 

{\bf Step 5}: Apply the monomial substitution $\Phi: \C^{2k} \to \C^{k}$ such that:
$$x_1 = z^{a_{11}}, \dots, x_k = z^{a_{1k}}, x_{k+1} = 1, \dots, x_{2k} = 1$$ to the function $f(P, x)$.

{\bf Step 4}: Return $g_1 * g_2 = z^{p_1} \Phi(f(P, x))$.

\end{Algorithm}

\section{Computer Algebra and applications to Statistics}\label{appAlg}

One focus of this thesis is applying Barvinok's rational functions to Statistics and Mixed Integer Programming. 
First we consider the connection between Computational Algebra and contingency tables in Statistics.

\begin{definition}\label{table}
A {\bf $s$-table of size $(n_1,\dots,n_s)$} is an array of non-negative
integers $v=(v_{i_1,\dots,i_s})$, $1\leq i_j\leq n_j$. For $0\leq L
<s$, an {\bf $L$-marginal} of $v$ is any of the $s\choose L$ possible
$L$-tables obtained by summing the entries over all but $L$
indices. 
\end{definition}
 
\begin{example}\label{tableEx}

Consider a $3$-table $X=(x_{i j k})$ of size $(m \mbox{, }n \mbox{, } p)$, where $m$, $n$, and $p$ are
natural numbers.  Let the integral
matrices $M_1=(a_{j k})$, $M_2=(b_{i k})$, and
$M_3=(c_{i j})$ be $2$-marginals of $X$, where $M_1$, $M_2$, and $M_3$ are integral matrices of type $n \times p$, $m \times p$, and
$m
\times n$
respectively.  
Then, a $3$-table $X=(x_{i j k})$ of size $(m \mbox{, }n \mbox{, } p)$ with given marginals satisfies the system of equations and inequalities:
\begin{equation}\label{tableEqu}
 \begin{array}{ll}
\sum _{i=1}^m x_{i j k}=a_{j k} \mbox{, } (j=1,2,...,n\mbox{, } 
k=1,2,...,p),\\
\sum _{j=1}^n x_{i j k}=b_{i k}\mbox{, } (i=1,2,...,m\mbox{, } 
k=1,2,...,p),\\
\sum _{k=1}^p x_{i j k}=c_{i j}\mbox{, } (i=1,2,...,m\mbox{, } 
j=1,2,...,n),\\
x_{i j k} \geq 0 \mbox{, } (i=1,2,...,m\mbox{, } 
j=1,2,...,n \mbox{, } k = 1,2,...,p).\\
\end{array}
\end{equation}



\end{example}
 
Such tables appear naturally in Statistics and Operations Research
under various names such as {\em multi-way contingency tables}, or
{\em tabular data}.  We consider the {\em table counting problem} and {\em table sampling problem}:

\begin{problem} (Table counting problem)

Given a prescribed collection of marginals, how many $d$-tables
are there that share these marginals?

\end{problem}

\begin{problem} (Table sampling problem)

Given a prescribed collection of marginals, generate typical tables that share these marginals.

\end{problem}

The table counting problem and table sampling problem have several
applications in statistical analysis, in particular for independence
testing, and have been the focus of much research \citep{Fienberg, De Loera, mcmc2, seth, mcmc}.  Given a specified
collection of marginals for $d$-tables of size $(n_1,\dots,n_d)$
(possibly together with specified lower and upper bounds on some of
the table entries) the associated {\em multi-index transportation
polytope} is the set of all non-negative {\em real valued} arrays
satisfying the given marginals and entry bounds specified in the system of equations and inequalities, such as formulas (\ref{tableEqu}) for a $3$-table given in Example \ref{tableEx}. The counting problem is the same as counting the number of integer points in the
associated multi-index transportation polytope.  

In this thesis, one of the main tools to solve the table counting problems and table sampling problems is Computational Algebra.  We consider a special ideal in the multivariate polynomial ring, namely a {\em toric ideal} $I_A$ associate to the given integral matrix $A$.  We compute the {\em Gr\"obner basis} associate to the toric ideal $I_A$.  Then we apply Gr\"obner bases of the toric ideal $I_A$ to the table counting problem and the table sampling problem. 
Here, we would like to remind the reader of some definitions.
\cite{coxlittleoshea} and \cite{sturmfels} are very good references for details. Let us denote $\Z^d_+ := \{ x \in \Z^d: x \geq 0 \}$ and $\Z_+ := \{x \in \Z : x \geq 0 \}$.
\begin{definition}
Let $K$ be any field and let $K[x]= K[x_1, x_2, \dots , x_d]$ be the polynomial ring in $d$ indeterminates.  A {\bf monomial} is a product of powers of variables in $K[x]$, i.e. $x^{\alpha_1}_1x_2^{\alpha_2} \dots x^{\alpha_d}_d$, where $\alpha_1, \, \alpha_2, \, \dots, \, \alpha_d \in \Z_+$.  
\end{definition}

\begin{definition}
Let $K$ be any field and let $K[x]= K[x_1, x_2, \dots , x_d]$ be the polynomial ring in $d$ indeterminates. 
Let $I \subset K[x]$.  Then we call $I$ an {\bf ideal} if it satisfies the following:

\begin{itemize}
\item $f + g \in I$ for all $f, \, g \in I$.
\item $af \in I$ for all $f \in I$ and all $a \in K[x]$.
\end{itemize}

\end{definition}

Note that by Hilbert basis theorem \cite[Chapter 2, section 5, Theorem 4]{coxlittleoshea} 
every ideal in $K[x]$ is generated by finitely many elements in $K[x]$.

\begin{definition}\label{termorder}
Let $\prec$ be a total order on $\Z_+^d$.  We call $\prec$ a {\bf term order} if it satisfies the following:

\begin{itemize}
\item For any $\alpha, \, \beta, \, \delta \in \Z_+^d$, $\alpha \prec \beta \rightarrow \alpha + \delta \prec \beta + \delta$.
\item For any $\alpha \in \Z^d_+ \backslash \{0\}$, $0 \prec \alpha$.
\end{itemize}
\end{definition}

A term order on $\Z^d_+$ gives a term order on the monomials of $K[x]$ by setting a bijection map from $\Z^d_+ $ to $ S \subset K[x]$, where $S$ is the set of all monomials in $K[x]$, such that 
$(a_1, \, a_2, \, \dots , \, a_d) \to x_1^{a_1}x_2^{a_2}\dots x_d^{a_d}$.  

\begin{definition}(The lexicographic term ordering)
Let $\alpha, \, \, \beta \in \Z^d_+$.  We say $\alpha \prec_{lex} \beta$ if the left most non-zero entry of $\alpha - \beta \in \Z^d$ is negative.  We write $x^{\alpha} \prec_{lex} x^{\beta}$ if $\alpha \prec_{lex} \beta$.
\end{definition}

For example, if we have $(3, 2, 7)$ and $(3, 5, 2)$ in $\R^3$,  then we have $(3, 2, 7) - (3, 5, 2) = (0, -3, 5)$.  So, $(3, 2, 7) \prec_{lex} (3, 5, 2)$ and $x_1^3x_2^2x_3^7 \prec_{lex} x_1^3x_2^5x_3^2$.
One notices that the lexicographic term ordering is a term order. 

We can also define a term order from a vector $c$, as we described in Definition \ref{termorder}, by the following method:
we make this vector $c$ into a term order $\prec_c$ such that
for all $\alpha, \, \beta \in \Z^d_+$, $\alpha \prec_c \beta$ if 
\begin{itemize}
\item  $c \cdot \alpha < c \cdot  \beta$ or
\item  $c  \alpha = c  \beta$ and $\alpha \prec_{lex} \beta$.
\end{itemize} 

For example, suppose $c = (1, 0, 2)$ and if we have $(3, 2, 7)$ and $(3, 5, 2)$ in $\R^3$,  then we have $(1, 0, 2)\cdot (3, 2, 7)= 17$ and $(1, 0, 2)\cdot (3, 5, 2) = 13$.  So, since $(1, 0, 2)\cdot (3, 5, 2) < (1, 0, 2)\cdot (3, 2, 7)$, we have $(3, 5, 2) \prec_c (3, 2, 7)$ and $x_1^3x_2^5x_3^2 \prec_c x_1^3x_2^2x_3^7$.

In general, any term order is defined by a $d \times d$ integral matrix $W$. 
We represent a term order $\prec$ on monomials in $x_1,\ldots,x_d$
by an integral $d \times d $-matrix $W$ as in
\citep{Mora+Robbiano}.  Two monomials satisfy $\,x^\alpha\prec
x^\beta \,$ if and only if $W\alpha$ is lexicographically smaller
than $W\beta$.  In other words, if $w_1,\ldots,w_d$ denote the
rows of $W$, there is some $j\in\{1,\ldots,d\}$ such that
$w_i\alpha=w_i\beta$ for $i<j$, and $w_j\alpha<w_j\beta$. For
example, $W=I_d$ describes the lexicographic term ordering. 
We will denote by $\prec_W$ the term order
defined by $W$.

\begin{definition}\label{defGB}
Let $K$ be any field and let $K[x]= K[x_1, x_2, \dots , x_d]$ be the polynomial ring in $d$ indeterminates.  Given a term order $\prec$, every non-zero polynomial $f \in K[x]$ has a unique initial monomial, denoted $in_{\prec}(f)$.  
If $ I$ is an ideal in $K[x]$, then its {\bf initial ideal} is the monomial ideal
$$ in_{\prec}(I):= < in_{\prec}(f): f \in I >. $$

The monomials which do not lie in $in_{\prec}(I)$ are called {\bf standard monomials}.  A finite subset $G \subset I$ is called a {\bf Gr\"obner basis} for $I$ with respect to $\prec$ 
if $in_{\prec}(I)$ is generated by $\{in_{\prec}(g): g \in G\}.$
A Gr\"obner basis is called {\bf reduced} if for any two distinct elements $g, \bar g \in G$, no terms of $\bar g$ is divisible by $in_{\prec}(g)$.
 
\end{definition}

\begin{proposition}\cite[Proposition 1, Chapter 6]{coxlittleoshea}

Let $G:=\{g_1, g_2, \dots , g_k\}$ be a Gr\"obner basis for an ideal $I \subset K[x]$ and let $f \in K[x]$.
Then there exists a unique $r \in K[x]$ such that:

\begin{itemize}
\item No term of $r$ is divisible by any of leading term of $g_i$, for all $i = 1, 2, \dots , k$.
\item There is $g \in I$ such that $f = g + r$.
\end{itemize}

In particular $r$ is the remainder on division of $f$ by $G$, and $r$ is unique no matter how the elements of $G$ are listed when using the division algorithm.
\end{proposition}

The remainder $r$ for $f \in K[x]$ is called the {\em normal form} of $f$.  Note that the reduced Gr\"obner basis is unique.  
This thesis concentrates in a special kind of ideals $I$ in $K[x]= K[x_1, x_2, \dots , x_d]$, which are called {\em toric ideals}.  {\em Toric ideals} find applications in Integer Programming, Computational Algebra, and Computational Statistics \citep{sturmfels}.

\begin{definition}
Fix a subset $A = \{a_1, a_2, \dots , a_d\}$ of $\Z^n$.  \\Each vector $a_i$ is identified with a monomial in the Laurent polynomial ring $K[\pm t] := K[t, t^2, \dots , t^d, t^{-1}, t^{-2}, \dots , t^{-d}]$.  
Consider the homomorphism induced by the monomial map $$\hat \pi : K[x] \to K[\pm t], \, x_i \to t^{a_i}.$$
Then the kernel of the homomorphism $\hat \pi$ is called the {\bf toric ideal} of $A$. 
\end{definition}

The following lemma describes the set of generators of a toric ideal $I_A$ associated to the integral matrix $A$.
 
\begin{lemma} \cite[Lemma 4.1]{sturmfels}\label{binomiallemma}
The toric ideal $I_A$ is spanned as a $K$-vector space by the set of binomials
$$\{x^u-x^v\, :\, \,  u, \, v \in \Z^d_+, \, \, Au = Av \}.$$
\end{lemma}

The main theorem on Chapter \ref{GBT} concerns a new way to compute the
\emph{toric ideal} $I_A$ of the integral matrix $A$.


Now we are ready to discuss applications of Computational Algebra to Computational Statistics. As we mentioned earlier, a toric ideal $I_A$ and the Gr\"obner basis associated to $I_A$ find applications to Computational Statistics \citep{sturmfels}. Here we would like to discuss how we can apply Gr\"obner bases to solve the table counting and table sampling problems.  First of all, we will remind the reader of the definition of {\em Markov bases} associate to the given integral matrix $A$ \citep{mcmc2}.

\begin{definition}
Let $P = \{x \in \R^d: Ax = b, x \geq 0 \} \not = \emptyset$, where $A \in \jdlZ^{n \times d}$ and $b \in \Z^n$, and let $M$ be a finite set such that $M \subset \{x \in \Z^d: Ax = 0 \}$.  Then we define the graph $G_b$ such that:
\begin{itemize}
\item Nodes of $G_b$ are lattice points inside $P$.
\item Draw a undirected edge between a node $u$ and a node $v$ if and only if $u - v \in M$.
\end{itemize}
Then we call $M$ a  Markov basis of the toric ideal associate to a matrix $A$ if $G_b$ is connected for all $b$ with $P \not = \emptyset$.  If $M$ is minimal with respect to inclusion, then we call $M$ a {\bf minimal} Markov basis.
\end{definition}

Note that, in general, a minimal Markov basis is not necessarily unique.
A Markov basis can be used for randomly sampling data and random walks on 
contingency tables \citep{DG, mcmc2}.  We will describe the {\em Monte Carlo Markov Chain} algorithm which uses Markov bases to create random walks on contingency tables. 
We can also define a Gr\"obner basis using a graph $G_b$.  

\begin{lemma}\label{GBdef} \cite[Theorem 5.5]{sturmfels}
Let $P = \{x \in \R^d: Ax = b, x \geq 0 \} \not = \emptyset$, where $A \in \jdlZ^{n \times d}$ and $b \in \Z^n$.   Let $M$ be a finite set such that $M \subset \{x \in \Z^d: Ax = 0 \}$ and let $\prec$ be any term order on $\N^d$.  Then we define the graph $G_b$ such that:
\begin{itemize}
\item Nodes of $G_b$ are lattice points inside $P$.
\item Draw a directed edge between a node $u$ and a node $v$ if and only if $u \prec v$ for $u - v \in M$.

\end{itemize}
If $G_b$ is acyclic and has a unique sink for all $b$ with $P \not = \emptyset$,
then $M$ is a Gr\"obner basis for a toric ideal associate to a matrix $A$ with respect to $\prec$.
\end{lemma}

Notice that if $M$ is a Gr\"obner basis then this implies $M$ is a Markov basis because if we have an acyclic directed graph with a unique sink, then it has to be connected.   However, note that not all Markov bases are Gr\"obner bases.

{\bf Remark:} Gr\"obner bases provide a way to generate Markov bases for a wide variety of problems where no natural set of moves were known.

\begin{example}\label{MBExample}

Suppose we have $2 \times 3$ tables with given marginals.

\begin{table}
\begin{center}
\begin{tabular}{|c|c|c|c||c|} \hline
     & & & & Total \\ \hline
     & ? ? ?  &    ? ? ?  &    ? ? ?  & 6\\ \hline
     & ? ? ?  &    ? ? ?  &   ? ? ?  & 6 \\ \hline
 Total    & 4 & 4 & 4 & \\ \hline
\end{tabular}
\caption{$2 \times 3$ tables with $1$-marginals.}\label{T23}
\end{center}
\end{table}

There are $19$ tables with these marginals for $2 \times 3 $ tables in Table \ref{T23}.  

 \begin{figure}[h]
 \begin{center}
     \includegraphics[width=7 cm]{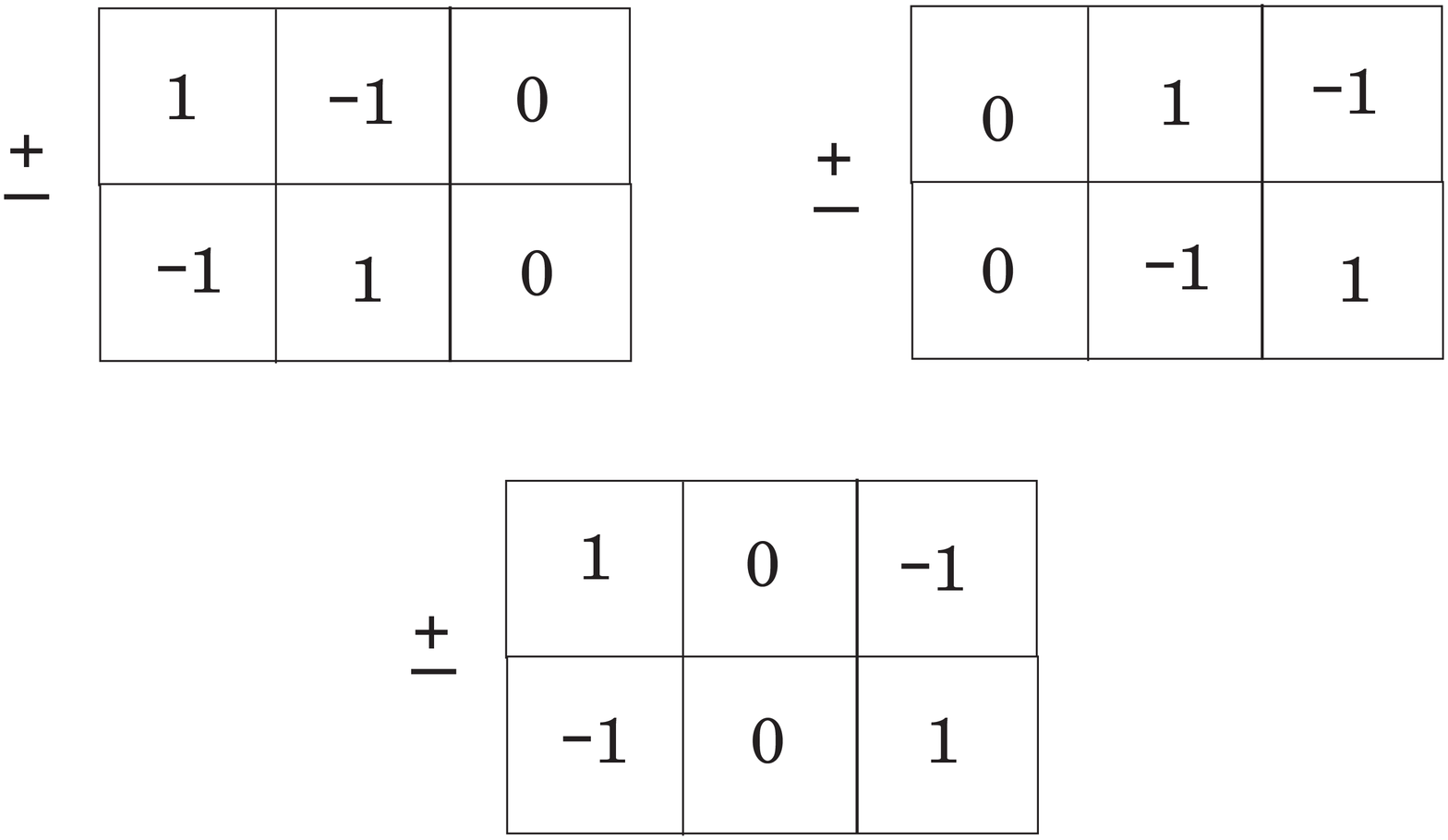}
\caption{Markov basis elements for $2 \times 3$ tables}
 \end{center}
 \end{figure} 

Up to signs, there are $3$ elements in the Markov basis.

 \begin{figure}[h]
 \begin{center}
     \includegraphics[width=11 cm]{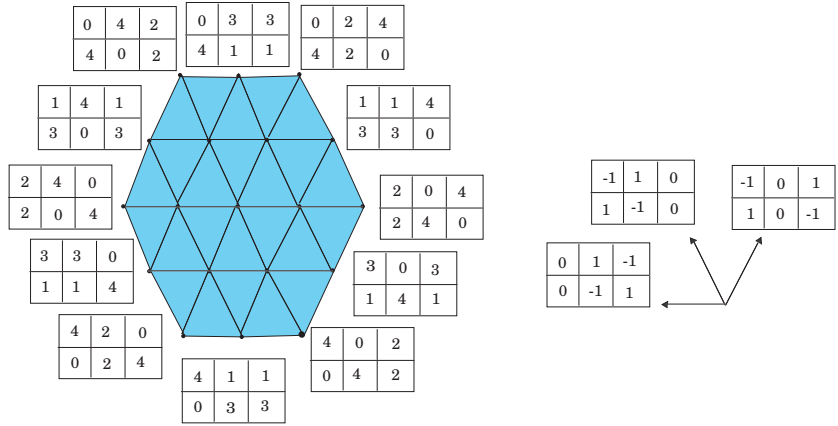}
 \end{center}
\caption{A connected graph $G_b$ for given $b$.}\label{MB}
 \end{figure} 

Figure \ref{MB} gives a connected graph for a Markov basis for $2 \times 3$ tables.  An element of the Markov basis is a undirected edge between integral points in the polytope.

\end{example}

From here on, we focus on Markov bases for contingency tables.
Why do we care about Markov bases for contingency tables? We care about them because using Markov bases, we can estimate 
the number of tables by Monte Carlo Markov Chain (MCMC) algorithm (\cite{mcmc2}).  \cite{DS} also showed that the rate of convergence for MCMC is $\delta^2$, where $\delta$ is the diameter of the graph $G_b$.  The outline of MCMC algorithm is the following:

\begin{Algorithm} \label{mcmcalg} \cite[Lemma 2.1]{mcmc2} 

(Random walk on a graph)

{\bf Input:} A Markov basis $M$ of a graph $G_b$ defined by the set of contingency tables with given marginals and an initial node $f_0$ in $G_b$.

{\bf Output:} A sample from the hypergeometric distribution $\sigma(\cdot)$ on $G_b$.

\begin{enumerate}
\item Set $f:= f_0$ and set $t = 0$.

\item {\em While} ($t < \delta^2$)

\begin{itemize}

\item Choose $u \in M$ uniformly and a sign $\epsilon = \pm 1$ with probability $1/2$ each independently from $u$.

\item If $f + \epsilon u \geq 0$ then move the chain from $f$ to $f + \epsilon u$ with probability $\min \{\sigma(f+\epsilon u)/\sigma(f), 1\}$.  If not, stay at $f$. Set $t = t+1$.

\end{itemize}

\end{enumerate}

\end{Algorithm}

\cite{mcmc2} showed that this random walk on the graph $G_b$ is a connected, reversible, aperiodic Markov chain on $G_b$ which converges to the hypergeometric distribution \cite[Lemma 2.1]{mcmc2}.  With the random walk on the graph $G_b$, we apply it to approximating the number of lattice points in a convex rational polytope $P$.  The idea for approximating the number of lattice points in a convex polytope $P$ is the following:

Suppose we take
a sequence of draws $H:= \{n_1, n_2 \dots, n_m\}$ randomly from the uniform distribution over $P$.  Let $p(n_i) = 1/|P|$ be the uniform distribution over $P$.  If we can simulate a lattice point $n_i \in P$ from a distribution $q( \cdot )$, where $q(t) > 0$ for all $t \in P$, then we have $$E[\frac{1}{q(t)}] = \sum_{t \in P}\frac{1}{q(t)}q(t) = |P|.$$  Hence by the Strong Law of Large Number \cite[(7.1) on page 56]{probability}, the estimation of $|P|$ is: $$\frac{1}{m} \sum_{i=1}^m\frac{1}{q(n_i)} .$$

\section{Algorithms for counting}

Enumerating the lattice points in a given polytope and counting the number of lattice points in a given polytope is very useful to Computational Statistics \citep{De Loera, DG, mcmc2,seth, mcmc}. So, in the next section, we will briefly discuss methods to count exactly the number of lattice points in a given convex polyhedron. Barvinok's method will be fully discussed in Chapter \ref{chapsoftware}

Using the multivariate generating function $f(P, z)$ for a polytope $P$, we can count the number of lattice points inside $P$.  In fact, the number of lattice points inside $P$ is $f(P, (1, 1, \dots , 1))$.

\noindent {\bf Example:} Let $P$ be the quadrangle with vertices 
$V_1=(0,0)$, $V_2=(5,0)$, $V_3=(4,2)$, and $V_4=(0,2)$. 

\begin{figure}[ht]
\begin{center}
\includegraphics[width=5 cm]{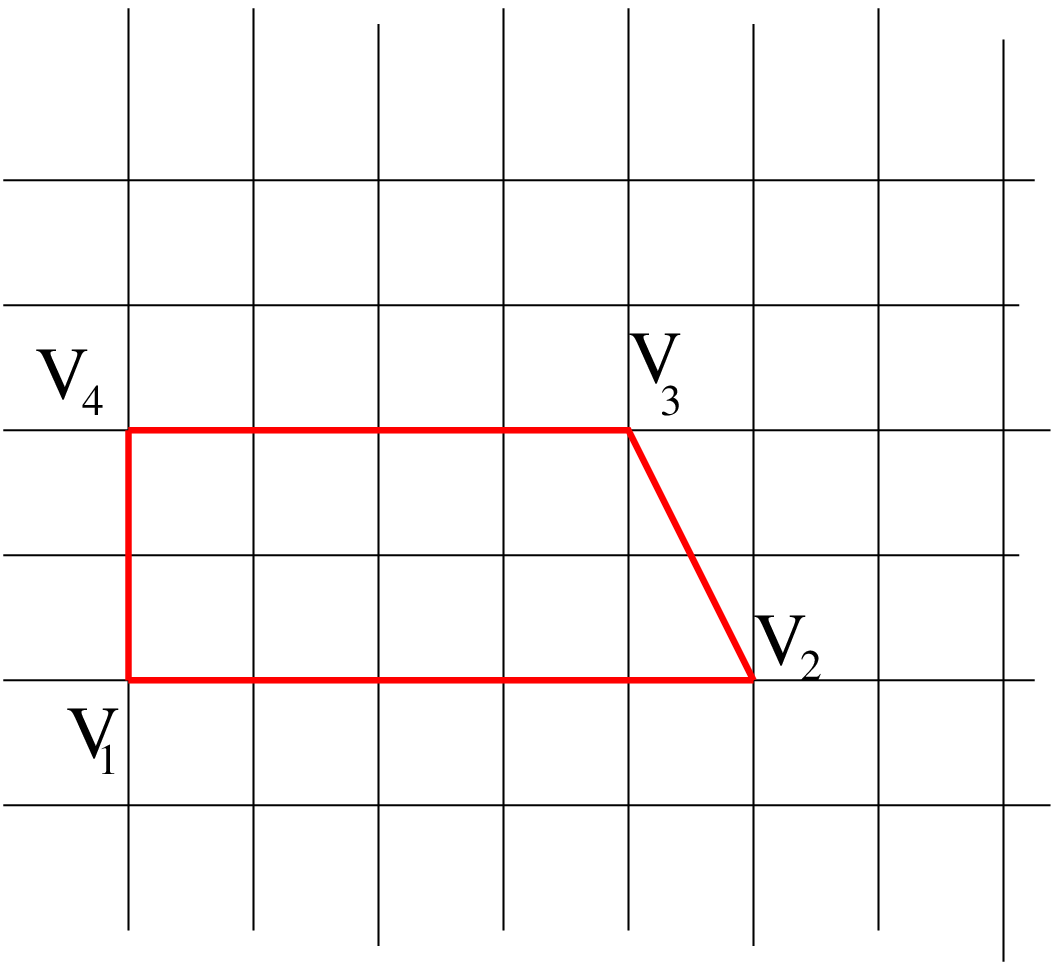} 
\end{center} 
\end{figure}

Then we have:

$f(P, z) = {z_1}^{5}+{z_1}^{4}z_2+{z_1}^{4}+{z_1}^{4}{z_2}^{2}+z_2{z_1}^{3}+{z_1}^{3}+{z_1}^{3}{z_2}^{2
}+z_2{z_1}^{2}+{z_1}^{2}+{z_1}^{2}{z_2}^{2}+z_1z_2+z_1+z_1{z_2}^{2}+{z_2}^{2}+z_2+1$. 

If we substitute $z_1 = 1$ and $z_2 = 1$ in $f(P, z)$, then we have $f(P, (1, 1)) = 16$, which is also the number of lattice points inside $P$.

We can use Barvinok's rational functions to count the number of lattice points inside a polytope.  Notice that, unfortunately, the point $(z_1, z_2, \dots , z_d) = (1, 1, \dots , 1)$ is a pole of Barvinok's rational functions.  Thus, we cannot directly substitute $(z_1, z_2, \dots , z_d) = (1, 1, \dots , 1)$ into $f(P, z)$.
Instead,  we compute $\lim_{z \to (1, \dots , 1)} f(P, z)$, which is the number of lattice points inside $P$.
In this thesis, we apply the residue calculus to compute $\lim_{z \to (1, \dots , 1)} f(P, z)$. We will discuss how to compute this limit via the residue calculus in Section \ref{limit}.

Several ``analytic'' algorithms have been proposed by many authors
\cite{vergne,Beck2,lasserre2, lasserre3, macmahon,pemantle2}. A couple of these
methods have been implemented and appear as the fastest for unimodular
polyhedra. However, only Barvinok's method has been implemented for arbitrary rational
polytopes. Consider, for example, Beck's method: let $M_i$ be the
columns of the matrix $M$.  We can interpret $P(M,b) \cap \jdlZ^d$ as
the Taylor coefficient of $z^b$ for the function $\Pi_{j=1}^d
\frac{1}{(1-z^{M_j})}$. One approach to obtain the particular
coefficient is to use the residue theorem. For example, it was seen
in \cite{Beck1} that

\[  P(M,b) \cap \jdlZ^d= \frac{ 1 }{ (2 \pi i)^m } \int_{
\left| z_1 \right| = \epsilon_1 } \cdots \int_{ \left| z_m \right| =
\epsilon_m } \frac{ z_1^{ - b_1 - 1 } \cdots z_m^{ - b_m - 1 } }{
\left( 1 - z^{ M_1 } \right) \cdots \left( 1 - z^{ M_d } \right) }
\ d z \ . 
\] 

Here $ 0 < \epsilon_1, \dots , \epsilon_m < 1 $ are different numbers
such that we can expand all the $ \displaystyle\frac{ 1 }{ 1 - z^{ M_k
} } $ into the power series about $0$. It is possible to do a partial
fraction decomposition of the integrand into a sum of simple
fractions. This was done very successfully to carry out very hard
computations regarding the Birkhoff polytopes \citep{Beck2}. Vergne and
collaborators have recently developed a powerful general theory about
the multivariate rational functions $\Pi_{j=1}^d
\frac{1}{(1-z^{M_j})}$ \citep{vergne,szenes}. Experimental results show
that it is a very fast method for unimodular polytopes
\citep{baldonideloeravergne}. Pemantle and Wilson \citep{pemantle2} have
pursued an even more general computational theory of rational
generating functions where the denominators are not necessarily
products of linear forms.  
\vskip 0.1in

Recently, \cite{lasserre3} introduced another method to enumerate the lattice points in a rational convex polyhedron.  Suppose we have a rational convex polyhedron $P_y=\{x \in \R^d: Ax = y, \, \, x \geq 0 \}$, where $A =(A_{ij}) \in \Z^{n \times d}$ and $b \in \Z^n$.  Then we define the function 
\begin{equation}\label{la1}
f(z) := \sum_{x \in P_y \cap \Z^d} e^{c\cdot x},
\end{equation}
 where $c \in \Z^d$ is small enough so that $f(z)$ is well defined.  A tool which \cite{lasserre3} use is a generating function $F: \C^{n} \to \C$:
\begin{equation}\label{la2}
z \to F(z):= \sum_{y \in \Z^n} f(y) z^y.
\end{equation}
Note that the generating functions in (\ref{la1}) and (\ref{la2}) are different from Barvinok's.  

Let us define $P^*_0 := \{b \in \R^d: b\cdot x \geq 0, \, x \in P_0\}$ and $\Gamma:=\{c \in \R^d: -c > b, \, \mbox{ for some } b \in P^*_0\}$.  Suppose $A_i$ is the $i$th column of a matrix $A$.  Then we have the following lemma.
\begin{lemma}\cite[Proposition 2.4]{lasserre3}

Let $f$ and $F$ be functions defined in (\ref{la1}) and (\ref{la2}), and let $c \in \Gamma$.
Then:

$$F(z) = \prod_{k=1}^d \frac{1}{(1-e^cz_1^{A_{1k}}\dots z_n^{A_{nk}})},$$ on the domain $(|z_1|, \dots , |z_n|) \in \{y \in \R^n: y > 0, \, e^{c_k}x^{A_k} < 1, k = 1, \dots, d\}$. 
\end{lemma}

\begin{definition}\cite[Definition 2.1]{lasserre3}
Let $p \in \N$ satisfy $n \leq p \leq d$ and let $\nu = \{\nu_1, \dots, \nu_p\} \subset \N$ be an ordered set with cardinality $|\nu| = p$ and $1 \leq \nu_1 \leq \dots \leq \nu_p \leq d$.  Then 
\begin{enumerate}
\item $\nu$ is said to be a basis of order $p$ if the $n \times p$ submatrix $A_{\nu} = [A_{\nu_1}| \dots | A_{\nu_p}]$ has the maximum rank.

\item For $n \leq p \leq d$, let $J_p :=\{\nu \subset \{1, \dots , d\}| \nu \mbox{ is a basis of order } p\}$.
\end{enumerate}

\end{definition}

Then \cite{lasserre3} show how to invert the generating function $F(z)$ in order to obtain the exact value of $f(y)$.  First, they determine an appropriate expansion of the generating function in the form: 

\begin{equation}\label{la3}
F(z) = \sum_{\sigma \in J_n} \frac{Q_{\sigma}}{\prod_{k \in \sigma}(1-e^{c_k}z^{A_k})},
\end{equation}

 where the coefficient $Q_{\sigma}: \C^n \to \C$ are rational functions with a finite Laurent series 

\begin{equation}\label{la4}
z \to Q_{\sigma} = \sum_{\beta \in \Z^n, ||\beta||\leq M} Q_{\sigma, \beta} z^{\beta},
\end{equation}

 for a strictly positive integer $M$.  
Then they apply the following theorem to obtain $f(y)$. 

\begin{theorem}\cite[Theorem 2.6]{lasserre3}

Let $A \in \Z^{n \times d}$ be of maximal rank, let $f$ be as in (\ref{la1}) with $c \in \Gamma$.  Assume that the generating function $F$ in (\ref{la2}) satisfies (\ref{la3}) and (\ref{la4}).  Then, $$f(y) = \sum_{\sigma \in J_n} \sum_{\beta \in \Z^n, ||\beta||\leq M} Q_{\sigma, \beta} E_{\sigma}(y - \beta)$$ with $$E_{\sigma}(y - \beta) =  \left \{ \begin{array}{ll} e^{c\cdot_{\sigma}x}, \mbox{ if } x := A^{-1}_{\sigma}(y-\beta) \in \N^n,\\
 0 \mbox{ otherwise;}\\ 
\end{array}\right .$$ where $c_{\sigma} = (c_{\sigma_1}, \dots, c_{\sigma_n})$.
\end{theorem}   

The main point is that as soon as we have $f(y)$ with sufficiently small $c$ vector, if we send $c_i \to 0$ for $i = 1, 2, \dots d$, then we can obtain the number of lattice points inside a convex rational polytope $P_y$ (if $P_y$ is a unbounded polyhedron, then this limit does not converge).

\section{Applications to Mixed Integer Programming}

Now we explain the connections between Barvinok's rational functions and integer programming problems.
We consider an integer linear programming problem:

\begin{question}\label{IP} (Integer Programming)

Suppose $A \in Z^{n \times d}$, $c \in Z^d$, and $b \in Z^n$.  We assume that the rank of $A$ is $n$.  Given a polyhedron $P=\{x \in \R^d: Ax= b, \, x \geq 0 \}$, we want to solve the following problem:

$$\mbox{(IP)  maximize} \ c \cdot x \ \mbox{subject to} \quad x \in P , \  \, x \in \Z^d.$$

\end{question}


These problems are called {\em integer programming problems} and we know that this problem is NP-hard \citep{karp}.  However, \cite{lenstra} showed that if we fixed the dimension, we can solve (IP) in polynomial time. 
Originally, Barvinok's counting algorithm relied
on H. Lenstra's polynomial time algorithm for Integer Programming in a
fixed number of variables \citep{lenstra}, but shortly after Barvinok's
breakthrough, \cite{Dyer} showed that this step can be
replaced by a short-vector computation using the $LLL$
algorithm. Therefore, using binary search, one can turn Barvinok's
counting oracle into an algorithm that solves integer programming
problems with a fixed number of variables in polynomial time (i.e. by
counting the number of lattice points in $P$ that satisfy $c\cdot x
\geq \alpha$, we can narrow the range for the maximum value of $c\cdot
x$, then we iteratively look for the largest $\alpha$ where the count is
non-zero). This idea was proposed by Barvinok in \citep{BarviPom}. We
call this IP algorithm the {\em BBS algorithm}.

To solve more general problems one might want to consider some variables as real numbers, instead of all integers.  Thus one can set the following problems:

\begin{question}\label{MIP} (Mixed Integer Programming)

Suppose $A \in Z^{n \times d}$, $c \in Z^d$, and $b \in Z^n$.  We assume that the rank of $A$ is $n$.  Given a polyhedron $P=\{x \in \R^d: Ax= b, \, x \geq 0 \}$, we want to solve the following problem:

$$\mbox{(MIP)   maximize} \ c \cdot x \ \mbox{subject to} \quad x \in P , \ \,  x_i \in {\mathbb Z} \ \mbox{if the index} \ i \in \ J \subset [d].$$

\end{question}

One notices that linear programming problems form a subset of (MIP) problems, i.e. $J = \emptyset$ and also integer programming problems form a subset of (MIP) problems, i.e. $J = [d]$.

As we mentioned earlier, 
we  can define the term order from a cost vector $c$, as we described in Definition \ref{termorder}.
From this term order $\prec_c$, as soon as we have the reduced Gr\"obner basis with the term order
$\prec_{c}$, we can solve the integer programming problem $\min\{c \cdot x : x \in P \cap \Z^d\}$.  A sketch of an algorithm is the following:

\begin{algorithm}\label{SturmfelsIP}\cite[Algorithm 5.6]{sturmfels}

{\bf Input:} A cost vector $c \in \Z^d$, a matrix $A \in \Z^{n \times d}$, a vector $b \in \Z^n$ and a feasible solution $v_0 \in P \cap \Z^d$, where $P:= \{x \in \R^d: Ax = b, \, \, x \geq 0\}$.

{\bf Output:} An optimal solution and the optimal value of minimize $c \cdot x$ subject to $x \in P \cap \Z^d$. 

{\bf Step 1:} Compute the Gr\"obner basis with the term order $\prec_c$.

{\bf Step 2:} Compute the normal form $x^u$ of $x^{v_0}$ and return $u$ and $cu$, which are an optimal solution and 
the optimal value, respectively.

\end{algorithm}

One notices that Algorithm \ref{SturmfelsIP} outputs the optimal value and an optimal solution for minimization.  Trivially if one wants to have the optimal value and an optimal solution for maximization such as (IP), one can set $c = -c$ and apply Algorithm \ref{SturmfelsIP}.

\begin{definition}
Suppose we have a rational convex polyhedron $P \subset \R^d$ and suppose we have an integer programming problem such that maximize $c \cdot x$ subject to $x \in P \cap \Z^d$.
Also suppose a feasible solution $x_0 \in P \cap \Z^d$ is given.
Then we call an integral vector $t \in \Z^d$ an {\bf augmenting vector} if $c(x_0 + t) > c x_0$.  
A finite set which contains all augmenting vectors is called a {\bf test set}. 
\end{definition}

There has been
considerable activity in the area of test sets and augmentation
methods (e.g. Graver, integral basis method, etc. See
\cite{survey2,survey1}).  One notices from Algorithm \ref{SturmfelsIP} that the Gr\"obner basis of a toric ideal $I_A$ associated to the integral matrix $A$ with respect to a term order $\prec_c$ is a test set for an integer programming problem $\min \{c\cdot x: Ax = b, \, \, x \geq 0, , x \in \Z^d\}$.
\vskip 0.1in

In 2003, Lasserre observed a new method for solving integer programming problems using Barvinok's short rational functions, which is different from the BBS algorithm.
We consider the integer programming
problem $maximize \{c\cdot x : Ax\leq b, x\geq 0, x\in\Z^d\}$, where
$c\in\Z^d$, $A\in\Z^{m\times d}$, and $b\in\Z^m$. This problem is equivalent to Problem \ref{IP}, because using Hermite normal form, we can project the polytope $P$ down into a lower dimension until the dimension of $P$ equals to the dimension of ambience space.  
We assume that the input system of inequalities $Ax\leq b,
x\geq 0$ defines a bounded polytope $P\subset\R^d$, such that
$P\cap\Z^d$ is nonempty. As before, all integer points are encoded as
a short rational function $f(P,z)$ in Equation (\ref{barvinokseries}) for $P$,
where the rational function is given in Barvinok's form.  Remember
that if we were to expand Equation (\ref{barvinokseries}) into monomials
(generally a very bad idea!) we would get $f(P,z) = \sum_{\alpha\in P
\cap\Z^d} z^{\alpha}$. For a given
$c\in\Z^d$, we make the substitution $z_i=t^{c_i}$, Equation (\ref{barvinokseries})
yields a univariate rational function in $t$:

\begin{equation}\label{Univariate Polynomial in t}
f(P,t)=\sum_{i\in I}{\pm \frac{t^{c\cdot u_i}}{\prod_{j=1}^d(1-t^{c\cdot v_{ij}})}}.
\end{equation}

The key observation is that if we make that substitution
directly into the monomial expansion of $f(P,z)$, we have
$z^\alpha\to t^{c\cdot\alpha}$. Moreover we would obtain the
relation 
\begin{equation}\label{eq:substituted}
f(P,t)=\sum_{\alpha\in P\cap\Z^d} t^{c\cdot\alpha} =
k t^M + \mbox{(lower degree terms)},
\end{equation}
where $M$ is the optimal value of our Integer Program and where $k$
counts the number of optimal integer solutions. Unfortunately, in
practice, $M$ and the number of lattice points in $P$ may be huge and we need
to avoid the monomial expansion step altogether. All computations
have to be done by manipulating short rational functions in (\ref{Univariate Polynomial in t}). 

\cite{lasserre} suggested the following approach: for $i\in
I$, define sets $\eta_i$ by $\eta_i = \{j\in\{1,...,d\} : c\cdot
v_{ij} > 0\}$ , and define vectors $w_i$ by $w_i = u_i -
\sum_{j\in\eta_i} v_{ij}.$ Let $n_i$ denote the cardinality of
$\eta_i$. Now we define $M=\max\{c\cdot w_i : i\in I\}$, $S=\{i\in I :
c\cdot w_i=M\}$ and we set $\sigma =\sum_{i\in S} E_i (-1)^{n_i}$. Note
that $M$ simply denotes the highest exponent of $t$ appearing in the
expansions of the rational functions defined for each $i\in I$ in
(\ref{Univariate Polynomial in t}). The number $\sigma$ is in fact the
sum of the coefficients of $t^M$ in these expressions, that is,
$\sigma$ is the coefficient of $t^M$ in $f(P,t)$.  Now with these
definitions and notation we can state the following result proved by
\cite{lasserre}.

\begin{theorem}\label{thm:aa} \cite[Theorem 3.1]{lasserre} 

If $c\cdot v_{ij}\neq 0$ for all $i\in I, j\in\{1,\ldots,d\}$, and if
$\sigma\neq 0$, then M is the optimal value $\pi$ of Integer
Program $maximize \{c\cdot x : Ax\leq b, x\geq 0, x\in\Z^d\}$.
\end{theorem}

\section{Summary of results in this thesis}\label{results}

In this section, we would like to summarize all of the results in this thesis.  In Chapter \ref{GBT}, the first theorem concerns the computation of the
\emph{toric ideal} $I_A$ of the matrix $A$.  

\begin{theorem} \label{mainsec1}\label{mainsec}
Let $A \in \jdlZ^{n \times d}$ and a term order $\prec_W$ specified by
a matrix $W$. Assuming that $n$ and $d$ are fixed, then there are
algorithms, that run in polynomial time in the size of the input data,
to perform the following four tasks:

\begin{enumerate}
\item Compute a short rational function $G$ which represents the
     reduced Gr\"obner basis of the toric ideal $I_A$ with respect to
     the term order $\prec_W$.

\item Decide whether the input monomial $x^a$ is in 
       normal form with respect to $G$.
\item Perform one step of the division algorithm modulo $G$.
\item Compute the normal form of the input monomial $x^a$ modulo the
Gr\"obner basis $G$.
\end{enumerate}
\end{theorem}

The proof of Theorem \ref{mainsec1} will be given in Section
\ref{toric}. Special attention will be paid to the Projection
Theorem \cite[Theorem 1.7]{newbar} since the projection of short
rational functions is the most difficult step to implement. Its
practical efficiency has yet to be investigated.

Theorem \ref{mainsec1} can be applied to prove many interesting theorems and corollaries in Computational Statistics and Integer Programming.  The following corollary can be proved by Theorem \ref{mainsec1} and the fact that a Gr\"obner basis associated to a toric ideal $I_A$ of the integral matrix $A$ is a Markov basis associate to $A$.

\begin{corollary}
Let $A \in \Z^{n \times d}$, where $d$ and $n$ are fixed. There
is a polynomial time algorithm to compute a multivariate rational
generating function for a Markov basis $M$ associated to $A$. This
is presented as a short sum of rational functions.
\end{corollary}

From Algorithm \ref{SturmfelsIP}, one can see that Theorem \ref{mainsec1} proves the following theorem:

\begin{corollary} \label{GBIP2}\label{GBIP}
Let $A \in \Z^{n \times d}$, $b \in \Z^n$, and $c \in \Z^d$. Given a polyhedron $P=\{x \in \R^d : Ax= b, \, x \geq 0 \}$, compute the mixed integer programming problem, via the Gr\"obner basis associated a toric ideal $I_A$ with the term order $\prec_c$ obtained by Barvinok's rational functions in Theorem \ref{mainsec},
$$\mbox{maximize} \ c \cdot x \ \mbox{subject to} \quad x \in P , \  \, x \in \Z^d,$$
in polynomial time if we fix $n$ and $d$.
\end{corollary}



Chapter \ref{GBT} introduces a new algorithm, improving on Barvinok's original work, which we call the {\em
homogenized Barvinok algorithm}. 
Like the original version in \citep{bar}, it runs in
polynomial time when the dimension is fixed. 
Then we will apply the homogenized Barvinok's algorithm to Commutative Algebra. 
The \emph{Hilbert series} of $S$ is the rational generating function
$\,\sum_{ a \in S} x^a$. 
Barvinok and Woods (2003) showed that
this Hilbert series can be computed as a short rational generating
function in polynomial time for fixed dimension. 
We show that this computation can be done without the
Projection Theorem (Lemma \ref{project}) when the semigroup is known to be normal.

\begin{theorem} \label{mainsec2} Under the hypothesis that
the ambient dimension $d$ is fixed,

\begin{enumerate}
\item the Ehrhart series  of a rational
 convex polytope given by linear inequalities can
be computed in polynomial time. The  Projection Theorem is not
used in the algorithm.

\item The same applies to computing the Hilbert series of a normal
semigroup $S$.
\end{enumerate}
\end{theorem}


Chapter \ref{IPT} discusses Integer and Mixed Integer Programming.  We describe a new mixed integer programming algorithm via Barvinok's rational functions.  This is a different approach from the method via the reduced Gr\"obner basis of the toric ideal $I_A$ associated to the integral matrix $A$.  
It is based on Theorem \ref{thm:aa} in Chapter \ref{IntroThesis}.
Also in Chapter \ref{IPT} we give a new algorithm to compute the optimal value and an optimal solution for any Integer Linear Program via Barvinok's rational functions. 
In Section \ref{Binary Search and the Digging
Algorithm}, we will show the performance of the {\em BBS algorithm} in some knapsack problems.  Chapter \ref{IPT} will show a proof of the following theorem:

\begin{theorem}\label{Main Theorem}
Let $A\in\Z^{m \times d}$, $b\in\Z^m$, $c\in\Z^d$, and assume that
number of variables $d$ is fixed. Suppose $P:= \{x \in \R^d: Ax\leq b, x\geq 0\}$ is a rational convex polytope in $\R^d$.  Given the mixed-integer programming problem
\[
\maximize \ c\cdot x \ \mbox{subject to} \quad x\in \{x \in \R^d: x \in P,\ x_i\in\Z\ 
\mbox{for}\ \ i\in J \subset \{1,\ldots,d\} \},
\]

\noindent {\bf (A)} We can use rational functions to encode the set of vectors (the $(A,b,c)$-test set): 

$\{u-v : u \ \hbox{is a} \ c-\hbox{optimal solution},\  v \ \hbox{feasible solution}, u,v \in \Z^d\},$

and then solve the MIP problem in time polynomial in the size of the input.

\noindent {\bf (B)} More strongly, the $(A,b,c)$ test set can be replaced by 
smaller test sets, such as Graver bases or reduced Gr\"obner bases.
\end{theorem}

We improve Lasserre's heuristic and give a third deterministic
IP algorithm based on Barvinok's rational function algorithms, the
\emph{digging algorithm}. In this case the algorithm can have an
exponential number of steps even for fixed dimension, but performs
well in practice. See Section \ref{Binary Search and the Digging
Algorithm} for details.

Chapter \ref{chapsoftware} concentrates on computational experiments and explains details of the implementation of the software package {\tt LattE}. 
We implemented the \emph{BBS algorithm} and the \emph{digging
algorithm} in the second release of the computer software {\tt LattE}.
We solved several challenging
knapsack problems and compared the performance of {\tt LattE} with the
mixed-integer programming solver CPLEX version 6.6. In fact the
digging algorithm is often surpassed by what we call the {\em single cone
digging algorithm}. See Section \ref{Computational Experiments} for
computational tests.
In Section \ref{experiments} we present some computational experience
with our current implementation of {\tt LattE}. We report on
experiments with families of well-known rational polytopes: multiway
contingency tables, knapsack type problems, and rational polygons. We
demonstrate that {\tt LattE} competes with commercial branch-and-bound
software and solves very hard instances and enumerates some examples
that had never been done before. We also tested the performance in the
case of two-way contingency tables and Kostant's partition function
where special purpose software has been written already
\cite{vergne,Beck2, deloerasturmfels,mount}. In Section \ref{formulas}
we present formulas for the Ehrhart quasi-polynomials of several
hypersimplices and truncations of cubes (e.g. the 24 cell). We show
solid evidence that Barvinok's ideas are practical and can be used to
solve non-trivial problems, both in Integer Programming and Symbolic
Computing. 

In Section \ref{theory} we present some experimental results with the {\em
homogenized Barvinok algorithm}.  It was recently implemented in {\tt
LattE}. Like we will show in Chapter \ref{GBT}, it runs in
polynomial time when the dimension is fixed. But it performs much
better in practice (1) when computing the Ehrhart series of
polytopes with few facets but many vertices; (2) when computing
the Hilbert series of normal semigroup rings.  We show its
effectiveness by solving the classical counting problems for $5
\times 5$ \emph{magic squares} (all row, column and diagonal sums
are equal) and $3 \times 3 \times 3 \times 3$ \emph{magic
cubes} (all line sums in the 4 possible coordinate
directions and the sums along main diagonal entries are equal).
Our computational results are presented in Theorem \ref{magic}.

\newpage
\pagestyle{myheadings}
\chapter{Gr\"obner bases of toric ideals via short rational functions}\label{GBT}
\thispagestyle{myheadings} 
\markright{  \rm \normalsize CHAPTER 2. \hspace{0.5cm}
 Gr\"obner bases of toric ideals via short rational functions}

The main techniques used in this thesis came from the Algebra of polynomial ideals.  We use special sets of generators called {\em Gr\"obner bases}. We deal with ideals associated with polyhedra that are called {\em toric ideals} (see definitions \ref{definitions}).
In this chapter we present polynomial-time algorithms
for computing with toric ideals and semigroup rings in fixed dimension. For
background on these algebraic objects and their interplay with
polyhedral geometry see \citep{stanley0,sturmfels,villarreal}. 
Our results are a direct application of recent
results by Barvinok and Woods (2003) on short encodings of
rational generating functions (such as Hilbert series).

\section{Computing Toric Ideals}  \label{toric}
From now on and without loss of generality we will assume that $ker(A) \cap
\real_{\geq 0}^d=\{0\}$.  This condition is not restrictive because
toric ideal problems can be reduced to this particular case via
homogenization of the problem. Our assumption implies that for all
$b$, the convex polyhedron $\, P \,= \, \{\, u \in \real^d \,\, :\,\,
A \cdot u = b \,\,\hbox{and}\,\, u \geq 0 \,\}\,$ is a polytope
(i.e. a bounded polytope) or the empty set. We begin by recalling
some useful results of Barvinok and Woods (2003):

\begin{lemma}\label{intersect}\cite[Theorem 3.6]{newbar}
Let $S_1, S_2$ be finite subsets of $\Z^d$, for $d$ fixed. Let
$f(S_1,x)$ and $f(S_2,x)$ be their
 generating functions, given as short
rational functions with at most $k$ binomials in each denominator.
Then there exists a polynomial time algorithm, which, given
$f(S_i,x)$, computes
$$ f(S_1 \cap S_2, x) \quad =  \quad \sum_{i \in I} \gamma_i \cdot \frac { x^{u_i} } {  (1-x^{v_{i1}})  \dots (1-x^{v_{is}}) }$$
with $s \leq 2k$, where the $\gamma_i$ are rational numbers,
$u_i,v_{ij}$ are nonzero integer vectors, and $I$ is a polynomial-size index
set.
\end{lemma}

The following lemma was proved by Barvinok and Woods using Lemma
\ref{intersect}:

\begin{lemma}\label{unioncomplement}\cite[Corollary 3.7]{newbar}
Let $S_1, S_2, \dots,S_m$ be finite subsets of $\Z^d$, for $d$
fixed. Let $f(S_i,x)$ for $i=1 \dots m$ be their generating
functions, given as short rational functions with at most $k$
binomials in each denominator.  Then there exists a polynomial
time algorithm, in the input size, which computes
$$ f(S_1 \cup S_2 \cup \dots S_m, x) \quad =  \quad \sum_{i \in I} \gamma_i \cdot \frac { x^{u_i} } {  (1-x^{v_{i1}})  \dots (1-x^{v_{is}}) }$$
with $s \leq 2k$, where the $\gamma_i$ are rational numbers,
$u_i,v_{ij}$ are nonzero integer vectors, and $I$ is a polynomial-size index
set. Similarly one can compute in polynomial time $f(S_1
\backslash S_2,x)$ as a short rational function.
\end{lemma}

We will use the \emph{Intersection Lemma} and the \emph{Boolean
Operation Lemma} to extract special monomials present in the
expansion of a generating function. The essential step in the
intersection algorithm is the use of the {\em Hadamard product}
(see Algorithm \ref{HadamardAlg}) and a special monomial substitution.
The Hadamard product is a bilinear operation on rational functions
(we denote it by $*$). The computation is carried out for pairs of
summands as in (\ref{barvinokseries}).  Note that the Hadamard
product $m_1 * m_2$ of two monomials $m_1,m_2$ is zero unless
$m_1=m_2$. We present an example of computing intersections.

\begin{example} \rm
Let $S_i=\{\,x \in \R: i-2 \leq x \leq i \,\}\cap \Z$ for $i=1,2$.
We rewrite their rational generating functions as in the proof of
Theorem 3.6 in \citep{newbar}: $f(S_1, z) = \frac{z^{-1}}{(1-z)} +
\frac{z}{(1-z^{-1})} = \frac{-z^{-2}}{(1-z^{-1})} +
\frac{z}{(1-z^{-1})} = g_{11} + g_{12},$ and $f(S_2, z) =
\frac{1}{(1-z)} + \frac{z^2}{(1-z^{-1})} =
\frac{-z^{-1}}{(1-z^{-1})} + \frac{z^2}{(1-z^{-1})} = g_{21} +
g_{22}$.

We need to compute four Hadamard products between rational
functions $g_{ij}$,whose denominators are products of binomials and whose
numerators are monomials. Lemma 3.4 in \cite{newbar} says that, 
these Hadamard products are essentially the same
as computing the rational function, as in Equation (\ref{barvinokseries}), 
of the auxiliary polyhedron $\{(\epsilon_1,\epsilon_2) |
p_1+a_1\epsilon_1=p_2+a_2\epsilon_2, \, \epsilon_i \geq 0 \}$.
Here $p_1,p_2$ are the exponents of numerators of ${g_{ij}} \ 's$
involved and $a_1,a_2$ are the exponents of the binomial
denominators. For example, the Hadamard product $g_{11} * g_{22}$
corresponds to the polyhedron $\,\{(\epsilon_1,\epsilon_2) |
\epsilon_2=4+\epsilon_1, \, \epsilon_i \geq 0 \}$. The
contribution of this half line is $- \frac{z^{-2}}{(1-z^{-1})}$.
We find
\begin{eqnarray*}
f(S_1, z) * f(S_2, z) \quad  = &  g_{11}*g_{21} +
g_{12} * g_{22} + g_{12}*g_{21} + g_{11} * g_{22} \\
= & \quad \frac{z^{-2}}{(1-z^{-1})} + \frac{z}{(1-z^{-1})} -
\frac{z^{-1}}{(1-z^{-1})} - \frac{z^{-2}}{(1-z^{-1})} \\
= & \frac{z - z^{-1}}{1 - z^{-1}} \quad = \quad 1 + z \quad  =
\quad f(S_1 \cap S_2, z).
\end{eqnarray*}
\end{example}

Another key subroutine introduced by Barvinok and Woods is the
following \emph {Projection Theorem}. In Lemmas
\ref{intersect}, \ref{unioncomplement}, and \ref{project}, 
the dimension $d$ is assumed to be fixed.

\begin{lemma} \label{project}\cite[Theorem 1.7]{newbar} 
Assume the dimension $d$ is a fixed constant. Consider  a rational
polytope $P \subset \real^d$ and a linear map $T: \Z^d \rightarrow
\Z^k$. There is a polynomial time algorithm which computes a short
representation of the generating function $\, f \bigl(T(P \cap
\Z^d),x\bigr) $.
\end{lemma}

Defining a term order $\prec_W$ by a $d \times d$ integral matrix $W$ (see details in Section \ref{appAlg}), we have the following lemma. 

\begin{lemma} \label{extractmonomial}
Let $S \subset \Z^d_+$ be a finite set of lattice points in the positive
orthant. Suppose the polynomial $\,f(S,x)
= \sum_{\beta \in S} x^\beta \, $ is represented as a short
rational function and let $\prec_W$ be a term order. We can
extract the (unique) leading monomial of $f(S,x)$ with respect to
$\prec_W$ in polynomial time.
\end{lemma}

\noindent {\em Proof:} The term order $\prec_W$ is represented by
an integer matrix $W$. For each of the rows $w_j$ of $W$ we
perform a monomial substitution $x_i:=x'_it^{w_{ji}}$. Note that
$t$ is a ``dummy variable'' that we will use to keep track of
elimination. Such a monomial substitution can be computed in
polynomial time by \cite[Theorem 2.6]{newbar}. The effect is that
the polynomial $f(S,x)$ gets replaced by a polynomial in the $t$
and the $x's$. After each substitution we determine the degree in
$t$.  This is done as follows: We want to do calculations in
univariate polynomials since this is faster so we consider the
polynomial $g(t)=f(S,1,t)$, where all variables except $t$ are set
to the constant one. Clearly the degree of $g(t)$ in $t$ is the
same as the degree of $f(S,x',t)$.  We create the \emph{interval
polynomial} $i_{[p,q]}(t)=\sum^q_{i=p} t^i$ which obviously has a
short rational function representation.  Compute the Hadamard
product of $i_{[p,q]}(t)$ with $g(t)$. This yields those
monomials whose degree in the variable $t$ lies between $p$ and
$q$. We will keep shrinking the interval $[p,q]$ until we find the
degree.  We need a bound for the degree in $t$ of $g(t)$ to start
a binary search.  An upper bound $U$ can be found via linear programming
or via the estimate in Theorem 3.1 of \citep{lasserre} which is an easy
manipulation of the numerator and denominator of the fractions in
$g(t)$. It is clear that $\log(U)$ is polynomially bounded. 
In no more than $\log(U)$ steps one can determine the
degree in $t$ of $f(S,x,t)$ by using a standard binary search
algorithm.

Let $\alpha$ be a polynomial-size upper bound on the highest total
degree of a monomial appearing in the generating function $f(S,x)$. We
can again apply linear programming or the estimate of \citep{lasserre}
to compute such an $\alpha$ (just as we computed $U$ before). Once the
highest degree $r$ in $t$ is known, we compute the Hadamard product of
$f(S,x,t)$ and $t^r h(x)$, where $h(x)$ is the rational generating
function encoding the lattice points contained inside the box
$[0,\alpha]^d$. This will capture only the desired monomials. Then
compute the limit as $t$ approaches $1$.  This can be done in
polynomial time using residue techniques. The limit represents the
subseries $\,H(S,x) = \sum_{ \beta \cdot w_j = r} x^{\beta} $.  Repeat
the monomial and highest degree search for the row
$w_{j+1}$,$w_{j+2}$, etc.  Since $\prec_W$ is a term order, after
doing this $d$ times we will have only one single monomial left, the
desired leading monomial.  \jdlqed

One has to be careful when using earlier Lemmas (especially the projection
theorem) that the sets in question are finite. We need the following
well-known bound:

\begin{lemma}
\label{gbexponentbound}\cite[Lemma 4.6 and Theorem 4.7]{sturmfels} 
Let $M$ be equal to $(n+1)(d-n)D(A)$, where $A$ is an $n \times d$
integral matrix and $D(A)$ is the biggest $n \times n$ subdeterminant
of $A$ in absolute value. Any entry of an exponent vector of any
reduced Gr\"obner basis for the toric ideal $I_A$ is less than $M$.
\end{lemma}

\begin{proposition} \label{nonreducedGB}
Let $A \in \jdlZ^{n \times d}$, $W \in \Z^{d \times d}$ specifying
a term order $\prec_W$. Assume that $n$ and $d$ are fixed.

1) There is a polynomial time algorithm to compute a short
rational function $G$ which represents a universal Gr\"obner basis
of $I_A$.

2) Suppose we are given the term order $\prec_W$ and a short rational
function encoding a finite set of binomials $x^u-x^v$ now expressed as
the sum of monomials $\sum x^uy^v$. Assume $M$
is an integer positive bound on the degree of any variable for any of
the monomials. One can compute in polynomial time a short rational
function encoding only those binomials $x^u-x^v$ that satisfy $x^v
\prec_W x^u$.

3) Suppose we are given a sum of short rational functions $f(x)$ which
is identical, in its monomial expansion, to a single monomial
$x^a$. Then in polynomial time we can recover the (unique) exponent
vector $a$.
\end{proposition}

\noindent {\em Proof:} 1) Set
 $\,M=(n+1)(d-n)D(A)\,$ where $D(A) $ is again the
largest absolute value of any $n \times n$-subdeterminant of $A$.
Using Barvinok's algorithm in \citep{bar}, we compute the
 following generating function in $2d$ variables:
$$ G(x,y) \quad = \quad \sum \bigl\{ \, x^u y^v \, \, : \,\, A u = A
v\, \,\,\hbox{and} \, \, 0 \leq u_i,v_i \leq M \, \bigr\}. $$ This is
the sum over all lattice points in a rational polytope. Lemma
\ref{gbexponentbound} above implies that the toric ideal $I_A$ is
generated by the finite set of binomials $x^u-x^v$ corresponding to
the terms $x^u y^v$ in $G(x,y)$. Moreover, these binomials are a
universal Gr\"obner basis of $I_A$.

2) Denote by $w_i$ the $i$-th row of the matrix $W$ which
specifies the term order. Suppose we are given a short rational
generating function $\,G_0(x,y) \, = \, \sum x^u y^v \,$
representing a set of binomials $\,x^u - x^v$ in $I_A$, for
instance $G_0 = G$ in part (1). In the following steps, we will
alter the series so that a term $x^u y^v$ gets removed whenever
$u$ is not bigger than $v$ in the term order $\prec_W$.  Starting
with $H_0 = G_0$, we perform Hadamard products with short rational
functions $f(S;x,y)$ for $ S \subset \Z^{2d}$.

Set $H_i = H_{i-1}* f(\{(u,v):w_iu=w_iv, \ 0 \leq u_j,v_j \leq M, \
j=1 \dots d\})$, and $G_i = H_{i-1} * f(\{(u,v):w_iu\geq w_iv+1, \ 0
\leq u_j,v_j \leq M \ j=1 \dots d\}).$ All monomials $x^uy^v\in G_j$
have the property that $w_iu=w_iv$ for $i<j$, $w_ju>w_jv$, and thus
$v\prec_W u$. On the other hand, if $v\prec_W u$ then there is some
$j$ such that $w_iu=w_iv$ for $i<j$, $w_ju>w_jv$, and we can conclude
that $x^uy^v \in G_j$.  Note that $H=G_1\cup G_2\cup\ldots\cup G_d$ is
actually a disjoint union of sets. The rational function that gives
the union, can be computed in polynomial time by Lemma
\ref{unioncomplement}. In practice, the rational generating functions
representing the $G_i$'s can be simply added together. The short
rational function $H$ encodes exactly those binomials in $G_0$ that
are correctly ordered with respect to $\prec_W$. We have proved our
claim since all of the above constructions can be done in polynomial
time.

3) Given $f(x)$ we can compute in polynomial time the partial
   derivative $\partial f(x)/\partial x_i$. This puts the exponent of
   $x_i$ as a coefficient of the unique monomial. Computing the
   derivative can be done in polynomial time by the quotient and
   product derivative rules.  Each time we differentiate a short
   rational function of the form

   $$\frac{x^{b_i}}{(1-x^{c_{1,i}})(1-x^{c_{2,i}})\dots
   (1-x^{c_{d,i}})}$$

  we add polynomially many (binomial type) factors to the
   numerator. The factors in the numerators should be expanded into
   monomials to have again summands in short rational canonical form
   $\frac{x^{b_i}}{(1-x^{c_{1,i}})(1-x^{c_{2,i}})\dots
   (1-x^{c_{d,i}})}$. Note that at most $2^d$ monomials
   appear each time  ($d$ is a constant).  Finally, if we take the limit when
   all variables $x_i$ go to one we will get the desired exponent.
\jdlqed

\begin{example} \rm Using {\tt LattE} we compute the set of
all binomials of degree less than or equal $10000$ in the toric
ideal $I_A$  of the matrix $\,A \, = \, \left [\begin
{array}{cccc} 1&1&1&1\\\noalign{\smallskip}0&1&2&3
\end {array}\right ]$. This matrix represents the
\emph{Twisted Cubic Curve} in algebraic geometry. We find that
there are exactly $195281738790588958143425$ such binomials. Each
binomial is encoded as a monomial $\,x_1^{u_1} x_2^{u_2} x_3^{u_3}
x_4^{u_4}
    y_1^{v_1} y_2^{v_2} y_3^{v_3} y_4^{v_4} $. The computation
takes about $40$ seconds. The output is a
 sum of $538$ simple rational functions of the form
a monomial divided by a product such as $ \left (1-{\frac
{x_{{3}}y_{{4}}}{x_{{1}}y_{{2}}}}\right ) \left (1-{ \frac
{x_{{1}}x_{{4}}y_{{2}}}{x_{{3}}}}\right ) \left (1-x_{{1}}y_{{1
}}\right ) \left (1-x_{{1}}x_{{3}}{y_{{2}}}^{2}\right ) \left (1
-x_{{3}}y_{{3}}\right ) \left (1-x_{{2}}y_{{2}}\right ) $. \jdlqed
\end{example}

\subsection*{\bf Proof of Theorem \ref{mainsec1}}

The proof of Theorem \ref{mainsec1} will require us to  project
and intersect sets of lattice points represented by rational
functions. We cannot, in principle, do those operations for
\emph{infinite} sets of lattice points. Fortunately, in our setting it
is possible to restrict our attention to finite sets. Besides Lemma
\ref{gbexponentbound} for the size of exponents of Gr\"obner bases, we
need a bound for the exponents of normal form monomials:

\begin{lemma} \label{nfbox}
Let $x^u$ be the normal form of $x^a$ with respect to the reduced
Gr\"obner basis $G$ of a toric ideal $I_A$ for the term order
$\prec_W$ (associated to the matrix $W$). Every coordinate of $u$ is
bounded above by $L = (n+1) d D(A) {\tilde a}$, where
$D(A)$ is the biggest subdeterminant of $A$ in absolute value,
$\tilde a$ denotes the largest coordinate of the exponent vector $a$.
\end{lemma}

\noindent {\em Proof:} We note that $u$ is a point in the (bounded)
convex polytope defined by the following inequalities in $v$: $A v = A a$, and $v \geq 0$ (it is forced to be bounded
for all $a$ because we assumed $ker(A) \cap \real_{\geq
0}^d=\{0\}$). Thus each coordinate of $u$ is bounded above by the
corresponding coordinate of some vertex of this polytope. Let $v$ be
such a vertex.  The non-zero entries of $v$ are given by 
$B^{-1} A a $ where $B $ is a maximal non-singular square 
submatrix of $A$. Clearly, each entry of $B^{-1} A$ is bounded
above by $D(A)$, and hence each entry of $v$ is bounded above by $L$.
We conclude that $L$ is an upper bound for the coordinates of $u$. \jdlqed

\noindent {\em Proof of Theorem \ref{mainsec1}:} Proposition
\ref{nonreducedGB} gives a Gr\"obner basis for the toric ideal $I_A$
in polynomial time. We now show how to get the reduced Gr\"obner basis
from it in three easy polynomial time steps.  The input is the the $n
\times d$ integral matrix $A$ and the $d \times d$ term order matrix
$W$.  The algorithm for claim (1) of Theorem \ref{mainsec1} has three
steps:

\noindent {\bf Step 1.} Let $M$ be equal to $(n+1)(d-n)D(A)$, as in Lemma
\ref{gbexponentbound}, for given input matrix $A$. As in Proposition
\ref{nonreducedGB}, compute the generating function which encodes
binomials of highest degree $M$ on variables that generate $I_A$:

$$ f(x,y) \quad = \quad \sum \bigl\{ \, x^u y^v \, \, : \,\, A u = A v
\,\,\hbox{and} \,\, 0 \leq u_j,v_j \leq M  \,\, \hbox{for} \,\,
j=1 \dots d \bigr\},$$

Next we wish to remove from $f(x,y)$ all incorrectly ordered binomials
(i.e. those monomials $x^uy^v$ with $u \prec_W v$ instead of the other
way around).  We do this using part 2 of Proposition
\ref{nonreducedGB}. We obtain from it a collection $G_0,G_1,\dots,G_d$
of rational functions encoding disjoint sets of lattice points. We
call $\bar f(x,y)$ the generating function representing the union of
$G_0,\dots,G_d$. This can be computed in polynomial time by adding the
rational functions of the $G_i$ together (since they are
disjoint). The reader should notice that this updated $\bar f(x,y)$
contains only those monomials of the old $f(x,y)$ that are now correctly
ordered.

Let $g_i(x)$ be the projection of $G_i$ onto the first group of
$x$-variables and denote by $g(x)$ the rational function that
represents the union of the $g_i(x)$. The rational function $g(x)$ can
be computed in polynomial time by the projection theorem of
Barvinok-Woods, i.e. Lemma \ref{project}.  It is important to note
that $g(x)$ is the result of projecting $\bar f(x,y)$ into the first
group of variables. This is true because a linear projection of the
union of disjoint lattice point sets (i.e. those represented by $G_i$)
equals the union of the projections of the individual sets. In
conclusion, $g(x)$ is the sum over all non-standard monomials having
degree at most $M$ in any variable.

\noindent {\bf Step 2.} Write $r(x,M)= \prod\limits_{i=1}^{d}(
\frac{1}{1-x_i} + \frac{x_i^M}{1-x_i^{-1}})$ for the generating
function of all $x$-monomials having degree at most  $M$ in any
variable. Note that this is a large, but finite, set of monomials.
We compute the following  Hadamard product of $d$ rational
functions in $x$ and Boolean complements (we denote them by
$\backslash$):
$$ \biggl( r(x,M) \backslash x_1 \cdot g(x) \biggr) *
   \biggl( r(x,M) \backslash x_2 \cdot g(x) \biggr) * \cdots *
   \biggl( r(x,M) \backslash x_d \cdot g(x) \biggr). $$ 

This is the generating function over those monomials all of whose
proper factors are standard modulo the toric ideal $I_A$ and whose
degree in any variable is at most $M$.

\noindent {\bf Step 3.} Let $h(x,y)$ denote the ordinary product of the resulting
rational function from Step 2 with 

$$ r(y,M) \backslash g(y) = \sum \bigl\{ \,y^v \,\,: \,\, v \,\,
\hbox{standard monomial modulo} \, I_A \, \hbox{of highest degree
$M$}\bigr\}.$$

Thus $h(x,y)$ is the sum of all monomials $x^uy^v$ such that $x^v$
is standard and $x^u$ is a monomial all of whose proper factors are
standard monomials modulo the toric ideal $I_A$ and, finally, the
highest degree in any variable is at most $M$.

Compute the Hadamard product $G(x,y):=\, \bar f(x,y) * h(x,y)$.
This is a short rational representation of a polynomial, namely,
it is the sum over all monomials $x^u y^v$ such that the binomial
$x^u - x^v$ is in the reduced Gr\"obner basis of $I_A$ with
respect to $W$ and $x^v \prec_W x^u$. This completes the proof of
the first claim of Theorem \ref{mainsec1}.

We next give the algorithm that solves claims 2 and 3 of Theorem
\ref{mainsec1}.  This will be done in four steps (1,2,3,4). We
are given an input monomial $x^a$ for which we aim to determine
whether it is already in normal form.

\noindent {\bf Step 4} Perform Steps 1,2,3.  Let $G(x,y)$ be the
reduced Gr\"obner basis of $I_A$ with respect to the term order $W$
encoded by the rational function obtained at the end of Step 3. Let
$r(x,{\tilde a})$ be, as before, the rational function of all
monomials having degree less than ${\tilde a}$ on any variable. Thus
$G'(x,y)=r(x,{\tilde a})\cdot G(x,y)$ consists of all monomials of the
form $x^s(x^uy^v)$ where $x^u-x^v$ is a binomial of the Gr\"obner
basis and where $0 \leq s \leq \tilde a$. Thus $x^sx^u$ is a monomial
divisible by some leading term of the Gr\"obner basis. 

Given a monomial $x^a$ consider $b(x,y)$, the rational function
representing the lattice points of $\,\{(u,v): u=a, \, 0\leq v_j\leq L
\,\,\hbox{for}\,\, j=1 \dots d \}$.  The Hadamard product
$\,\bar{G}(x,y) = G'(x,y)* b(x,y)$ is computable in polynomial time
and corresponds to those binomials in $G(x,y)$ that can reduce
$x^a$. If $\bar{G} (x,y) $ is empty then $x^a$ is in normal form
already, otherwise we use Lemma \ref{extractmonomial} and part 3 of
Proposition \ref{nonreducedGB} to find an element
$x^uy^v\in\bar{G}(x,y)$ and reduce $x^a$ to $x^{a-u+v}$. We may assume
that the coefficient of the encoded monomial is one, because we can compute
the coefficient in polynomial time using residue techniques, and divide
our rational function through by it.

Finally, we present the algorithm for claim 4 in Theorem
\ref{mainsec1} in four steps (1,2,5,6). A curious byproduct of
representing Gr\"obner bases with short rational functions is that the
reduction to normal form need not be done by dividing several times
anymore.

\noindent {\bf Step 5.}  Redo all the calculations of the Steps 1,2,3
using $L = (n+1) d D(A) {\tilde a}$ from Lemma \ref{nfbox}
instead of $M$. Note that the logarithm of $L$ is still bounded by a
polynomial in the size of the input data ($A,W,a$).  Let $\bar f(x,y)$
and $g(x)$ from Step 1,2 (now recomputed with the new bound $L$) and
compute the Hadamard product

$$ H(x,y) \quad := \quad \bar f(x,y) * \biggl( \,  r(x,L) \cdot
\bigl( r(y,L) \backslash g(y) \bigr) \biggr). $$ 

This is the sum over all monomials $x^u y^v$ where $x^v$ is the normal
form of $x^u$ and highest degree of $x^u$ on any variable is
$L$. Since we took a high enough degree, by Lemma \ref{nfbox}, the
monomial $x^ay^p$, with $x^p$ the normal form of $x^a$, is sure to be 
present.

\noindent {\bf Step 6.} 
We use $H(x,y)$ as one would use a traditional Gr\"obner
basis of the ideal $I_A$.  The normal form of a monomial $x^a$ is
computed by forming the Hadamard product $H(x,y) * (x^a \cdot r(y,L)).$
Since this is strictly speaking a sum of rational functions equal
to a single monomial, applying Part 3 of Proposition
\ref{nonreducedGB} completes the proof of Theorem \ref{mainsec1}.
\jdlqed

\section{Computing Normal Semigroup Rings}
\label{homogenized}

We will show in Chapter \ref{chapsoftware} that a major practical bottleneck of
the original Barvinok algorithm in \citep{bar} is the fact that a
polytope may have too many vertices. Since originally one visits
each vertex to compute a rational function at each tangent cone,
the result can be costly. For example, the well-known polytope of
semi-magic cubes in the $ 4\times 4 \times 4$ case has over two
million vertices, but only 64 linear inequalities describe the
polytope.  In such cases we propose a homogenization of Barvinok's
algorithm working with a single cone.

There is a second motivation for looking at the homogenization.
Barvinok and Woods \citep{newbar} proved that the Hilbert series
of semigroup rings can be computed in polynomial time. We show
that for {\em normal semigroup rings} this can be done simpler and
more directly, without using the Projection Theorem.

Given a rational polytope $P$ in $\real^{d-1}$, we set
 $i(P,m)=\# \{ z \in \Z^{d-1} : z\in m P \}$.
The {\em Ehrhart series} of $P$ is the generating function
$\,\sum_{m=0}^\infty i(P,m) t^m$.

\begin{Algorithm}[Homogenized Barvinok  algorithm]  \ \
\label{homogbarvinok}

\noindent {\bf Input:} A  full-dimensional, rational convex
polytope $P$ in $\real^{d-1}$ specified by  linear inequalities
and linear equations. \ \

{\bf Output:} The Ehrhart series of $P$. \rm
\begin{enumerate}
\item Place the polytope $P$ into the hyperplane defined by $x_d = 1$ in $\real^d$.
Let $K$ be the $d$-dimensional cone over $P$, that is,
$K=cone(\{(p,1) : p \in P\})$.

\item Compute the polar cone $K^*$.  The normal vectors of the facets
of $K$ are exactly the extreme rays of $K^*$.  If the polytope $P$
has far fewer facets then vertices, then the number of rays of the
cone $K^*$ is small.

\item Apply Barvinok's decomposition of $K^*$ into unimodular cones.
Polarize back each of these cones. It is known, e.g. Corollary 2.8
in \citep{BarviPom}, that by dualizing back we get a unimodular
cone decomposition of $K$. All these cones have the same dimension
as $K$. Retrieve a signed sum of multivariate rational functions
which represents the series $ \,\sum_{a \in K \cap \Z^{d}} x^a$.

\item Replace the variables $x_i$ by $1$ for $i \leq d-1$
and output the resulting series in $t = x_d$. This can be done
using the methods in Chapter \ref{chapsoftware}.
\end{enumerate}
\end{Algorithm}

One of the key steps in Barvinok's algorithm is
that any cone can be decomposed as the signed sum of (indicator
functions of) unimodular cones.  We will talk about this in detail on Section \ref{decomposition}, Chapter \ref{chapsoftware}.

\begin{theorem}[see \citep{bar}] \label{barvil} Fix the dimension $d$.
Then there exists a polynomial time algorithm which decomposes a
rational polyhedral cone $K \subset \real^d$ into unimodular cones
$K_i$ with numbers $\epsilon_i \in \{-1, 1\}$ such that $$f(K)
\,\, = \,\, \sum_{i \in I} \epsilon_i f(K_i) \mbox{, } \quad |I| <
\infty .$$
\end{theorem}

Originally, Barvinok had pasted together such formulas, one for each vertex of a polytope, using a result of Brion. 
Using Algorithm \ref{homogbarvinok}, we can prove Theorem \ref{mainsec2}.

\noindent{\em Proof of Theorem \ref{mainsec2}:} We first prove
part (1). The algorithm solving the problems is Algorithm
\ref{homogbarvinok}. Steps 1 and 2 are polynomial when the
dimension is fixed. Step 3 follows from Theorem \ref{barvil}. We
require a special monomial substitution, possibly with some poles.
This can be done in polynomial time by \citep{newbar}.

Part (2): Recall the characterization of the integral closure of
the semigroup $S$ as the intersection of a pointed polyhedral cone
with the lattice $\Z^d$. From this it is clear that Algorithm
\ref{homogbarvinok} computes the desired Hilbert series, with the
only modification that the input cone $K$ is given by the rays of
the cone and not the facet inequalities. The rays are the
generators of the monomial algebra. But, in fixed dimension, one
can transfer from the extreme rays representation of the cone to
the facet representation of the cone in polynomial time.  \jdlqed

\newpage
\pagestyle{myheadings}
\chapter{Theoretical applications of rational functions to Mixed Integer Programming}\label{IPT}
\thispagestyle{myheadings} 
\markright{  \rm \normalsize CHAPTER 3. \hspace{0.5cm}
 Theoretical applications of rational functions to MIP}

We now discuss how all these ideas can be used in Discrete Optimization.

\section{The $(A,b,c)$ test set algorithm}
\label{complexity}

In all our discussions below, the input data are an $m\times d$
integral matrix $A$ and an integral $m$-vector $b$. For simplicity we
assume it describes a polytope $P=\{x \in \real^d| Ax\leq b, x\geq
0\}$.
We assume that there are no redundant inequalities and no hidden equations in the system.
This polytope $P=\{x \in \real^d| Ax\leq b, x\geq
0\}$ is equivalent to the expression of $\{x \in \real^{\bar d}| \bar Ax = b, x\geq
0\}$, where $\bar A \in \Z^{m \times \bar d}$ and $\bar d = m + d$.  We can transform the expression $\{x \in \real^{\bar d}| \bar Ax = b, x\geq
0\}$ to the expression $\{x \in \real^d| Ax\leq b, x\geq
0\}$ by projecting down $P$ to a full dimensional polytope with Hermite normal form and we can also transform the expression $\{x \in \real^d| Ax\leq b, x\geq
0\}$ to the expression $\{x \in \real^{\bar d}| \bar Ax = b, x\geq
0\}$ by introducing slack variables.

First we would like to remind a reader of Barvinok's rational functions (see details on Chapter \ref{IntroThesis}).  By Theorem \ref{decombar}, with a given $P$, if we fix $d$ there is a polynomial time algorithm to compute Barvinok's short rational functions in the form of
\begin{equation} \label{eq:aa}
f(P,z) = \sum_{i\in I} {E_i \frac{z^{u_i}} {\prod\limits_{j=1}^d
(1-z^{v_{ij}})}}, 
\end{equation}
\noindent where $I$ is a polynomial sized finite indexing set, and where $E_i \in \{1,-1\}$ and $u_i, v_{ij} \in {\mathbb{Z}}^d$ for all $i$ and
$j$.
In this section we will show how to apply Barvinok's short rational functions in (\ref{eq:aa}) to Mixed Integer Programming.

\noindent {\em Proof of Theorem \ref{Main Theorem}:} We only show the proof of part (A). The proof of part (B) appears in
Chapter \ref{GBT}. We first explain how to solve integer programs
(where all variables are demanded to be integral). This part of the proof is essentially
the proof of Lemma 3.1 given in \cite{Hosten+Sturmfels:IPgap} for the
case $Ax=b, \ x \geq 0$, instead of $Ax\leq b$, but we emphasize the
fact that $b$ is fixed here. We will see how the techniques can be extended
to mixed integer programs later.
For a positive integer $R$, let $$r(x,R)= \prod\limits_{i=1}^{n}(
\frac{1}{1-x_i} + \frac{x_i^R}{1-x_i^{-1}})$$
be the generating
function encoding all $x$-monomials in the positive orthant, 
having degree at most  $R$ in any variable.  
Note that this is a large, but finite, set of monomials.
Suppose $P = \{ x \in \R^d: Ax \leq b, \, x \geq 0\}$ is a nonempty polytope.
Using Barvinok's algorithm in \cite{BarviPom}, compute the
following generating function in $2d$ variables:
\[ 
f(x,y) =  \sum \bigl\{
\, x^u y^v \, \, : \,\, Au \leq b, \, Av \leq b, \, u, \, v \geq 0, \mbox{ } c\cdot u-c\cdot v\geq 1, \mbox{ and } u, v \in \Z^d \bigr\}.
\]
This is possible because we are clearly dealing with the lattice points of a 
rational polytope. The monomial expansion of $f(x,y)$ exhibits a clear
order on the variables: $x^uy^v$ where $c\cdot u>c\cdot v$. Hence $v$ is not an
optimal solution. In fact, optimal solutions will never appear as
exponents in the $y$ variables. 

Now let $g(y)$ be the projection of $f(x,y)$ onto the $y$-variables
variables. Thus $g(y)$ is encoding all non-optimal feasible integral
vectors (because the exponent vectors of the $x$'s are better feasible
solutions, by construction), and it can be computed from $f(x,y)$ in
polynomial time by Lemma \ref{project}. Let $V(P)$ be the vertex set of $P$ and choose an integer $R \geq \max\{v_i: v \in V(P), \, 1 \leq i \leq d \}$ (we can find such an integer $R$ via linear programming).  Define $f(x,y)$ and $g(x)$ as above and compute the Hadamard product
\[ 
H(x,y) \quad := \quad f(x,y) * 
\left[
\left(r(x, R)- g(x)\right) r(y, R)
\right]. 
\] 
This is the sum over all monomials $x^u y^v$ where $u, \, v \in P$ and 
where $u$ is an optimal
solution. The reader should note that the vectors $u-v$ form a test
set (an enormous one), since they can be used to improve from any
feasible non-optimal solution $v$. This set is what we called the
$(A,b,c)$-test set.  It should be noted that one may replace $H(x,y)$
by a similar encoding of other test sets, like the Graver test set or
a Gr{\"o}bner basis (see Chapter \ref{GBT} for details).

We now use $H(x,y)$ as one would use a traditional test set for
finding an optimal solution: Find a feasible solution $a$ inside the
polytope $P$ using Lemma
\ref{extractmonomial} and Barvinok's Equation (\ref{eq:aa}).
Improve or augment to an optimal solution by computing
the Hadamard product
\[
H(x,y) * \left(y^a  r(x, R)\right).
\] 

The result is the set of monomials of the form $x^uy^a$ where $u$ is
an optimal solution. One monomial of the set, say the lexicographic
largest, can be obtained by applying Lemma \ref{extractmonomial}. This
concludes the proof of the case when all variables are integral.

Now we look at the mixed integer programming case, where only $x_i$
with $i\in J\subset \{1,\ldots,d\}$ are required to be integral.
Without loss of generality, we may assume that $J=\{r,\ldots,d\}$ for
some $r$, $1\leq r\leq d$. Thus, splitting $A$ into $(B|C)$, we may
write the polytope $P$ as $\{(x,x'):Bx+Cx'\leq b, \ x,x'\geq 0\}$ where
the variables corresponding to $B$ are not demanded to be integral.
Consider a vertex optimal solution $\bar{x}$ to the mixed integer
problem. The first key observation is that its fractional part can be
written as $\bar{x}_J=\hat{B}^{-1}(b-C\bar{x}')$ where $b-C\bar{x}'$
is an integer vector. Here $\hat{B}^{-1}$ denotes the inverse of a
submatrix of $B$. This follows from the theory of linear programming,
when we solve the mixed integer program for fixed $x'=\bar{x}'$.

The denominators appearing are then contributed by
$\hat{B}^{-1}$. Then every appearing denominator is a factor of $M$,
the least common multiple of all determinants of a square submatrix of
$A$. It is clear $M$ can be computed in polynomial time in the size of the
input. This complexity bound holds, since the number of such square
submatrices is bounded by a polynomial in $m$, the number of rows of
$A$, of degree $d$, the number of columns of $A$. Moreover, each of
these determinants can be computed in time polynomial in the size of
the input, and therefore, $M$ itself can be computed in time
polynomial in the size of the input in fixed dimension $d$.
Thanks to this information, we know that if we dilate the
original polytope $P$ by $M$, the optimal solutions of the mixed
integer program become, in the dilation $MP$, optimal integral
solutions of the problem
\[
\maximize \ c\cdot x \ \mbox{subject to} \quad x \in MP, \ x\in\Z^d 
\]
with the additional condition that the coordinates with index in $J$
are multiples of $M$. Ignoring this condition involving multiples of
$M$ for a moment, we see that, as we did before, we can obtain an encoding of
\emph{all} optimal improvements as a generating function $H(x,y)$. 

Let $\bar r (x, R) = \left[\prod\limits_{i\not\in J}\left(\frac{1}{1-x_i} + \frac{x_i^R}{1-x_i^{-1}}\right)\right] \left[\prod\limits_{i\in J}\left(\frac{1}{1-x_i^M} + \frac{x_i^{RM}}{1-x_i^{-M}}\right)\right]$.
To extract those vectors whose coordinates indexed by $J$ are
multiples of $M$, we only need to intersect (Hadamard product
again) our generating function $H(x,y)$ with the generating function 
$\bar r(x, R) \bar r(y, R)$. 
Then only those vectors whose coordinates indexed by $J$ are multiples
of $M$ remain. This completes the proof of the theorem. 
\jdlqed 

\section{The Digging Algorithm}
\label{Binary Search and the Digging Algorithm}

In what follows we present a strengthening of Lasserre's heuristic and
discuss how to use Barvinok's short rational functions to solve
integer programs using digging. Suppose $A \in \Z^{m \times d}$, $b \in \Z^m$ and finite $\Xi \subset \Z^d$ are given.  
We consider the family of integer programming
problems of the form $\maximize \{c\cdot x : Ax\leq b, x\geq 0, x\in\Z^d\}$, where
$c\in \Xi$. We assume that the input system of inequalities $Ax\leq b,
x\geq 0$ defines a polytope $P\subset\R^d$, such that
$P\cap\Z^d$ is nonempty. 
When the hypotheses of Theorem \ref{thm:aa} are met,
from an easy inspection, we could recover the optimal value of an
integer program. If we assume that $c$ is chosen randomly from some large
cube in $\Z^d$, then the first condition
is easy to obtain. Unfortunately, our computational
experiments (see Section \ref{Computational Experiments}) indicate
that the condition $\sigma\neq 0$ is satisfied only occasionally. Thus
an improvement on the approach that Lasserre proposed is needed to
make the heuristic terminate in all instances. Here we explain the
details of an algorithm that \emph{digs} for the coefficient of the
next highest appearing exponent of $t$. For simplicity our explanation
assumes the easy-to-achieve condition $c\cdot v_{ij}\neq 0$, for all $ v_{ij}$.

As before, take Equation (\ref{eq:aa}) computed via Barvinok's
algorithm.  Now, for the given $c$, we make the substitutions
$z_k=y_k t^{c_k}$, for $k = 1, \ldots, d$. Then substitution into
(\ref{eq:aa}) yields a sum of multivariate rational functions in the vector
variable $y$ and scalar variable $t$:

\begin{equation} \label{eq:ab}
g(P;y,t) = \sum_{i\in I} {E_i\frac{y^{u_i} t^{c\cdot u_i}} 
{\prod_{j=1}^d (1-y^{v_{ij}} t^{c\cdot v_{ij}})}}. 
\end{equation}

On the other hand, the substitution on the left-side of Equation (\ref{eq:aa}) 
gives the following sum of monomials.

\begin{equation} \label{eq:d}
g(P;y,t)=\sum_{\alpha \in P \cap \Z^d} y^\alpha t^{c\cdot\alpha}.
\end{equation}

Both equations, (\ref{eq:d}) and (\ref{eq:ab}), represent the same
function $g(P;y,t)$. Thus, if we compute a Laurent series of
(\ref{eq:ab}) that shares a region of convergence with the series in
(\ref{eq:d}), then the corresponding coefficients of both series must
be equal. In particular, because $P$ is a polytope, the series in
(\ref{eq:d}) converges almost everywhere.  Thus if we compute a
Laurent series of (\ref{eq:ab}) that has any nonempty region of
convergence, then the corresponding coefficients of both series must
be equal. Barvinok's algorithm provides us with the right hand side of
(\ref{eq:ab}). We need to obtain the coefficient of highest
degree in $t$ from the expanded Equation (ref{eq:d}). We compute a
Laurent series for it using the following procedure: Apply the
identity
 
\begin{equation}\label{eq:i}
\frac{1}{1-y^{v_{ij}} t^{c\cdot v_{ij}}}=\frac{-y^{-v_{ij}} 
t^{-c\cdot v_{ij}}}{1 - y^{-v_{ij}} t^{-c\cdot v_{ij}}}
\end{equation}

\noindent to Equation (\ref{eq:ab}), so that any $v_{ij}$ such that $c\cdot
v_{ij} > 0$ can be changed in ``sign'' to be sure that, for all
$v_{ij}$ in (\ref{eq:ab}), $c \cdot v_{ij} < 0$ is satisfied (we may have to
change some of the $E_i$, $u_{i}$ and $v_{ij}$ using our identity, but
we abuse notation and still refer to the new signs as
$E_i$ and the new numerator vectors as $u_{i}$ and the new denominator
vectors as $v_{ij}$). Then, for each of the rational functions in the
sum of Equation (\ref{eq:ab}) compute a Laurent series of the form

\begin{equation}\label{eq:j}
E_i \, y^{u_i} t^{c\cdot u_i} \prod_{j=1}^d (1+y^{v_{ij}} t^{c\cdot
  v_{ij}}+ {(y^{v_{ij}} t^{c\cdot v_{ij}})}^2 + \ldots).
\end{equation}

Multiply out each such product of series and add the resultant
series.  This yields precisely the Laurent series in (\ref{eq:d}). 
Thus, we have an algorithm to solve
integer programs:

\noindent{\bf Algorithm:} (\emph{Digging Algorithm}):

\noindent{\bf Input:} $A \in \Z^{m \times d}, \, b \in \Z^m, \, c \in \Xi$.

\noindent{\bf Output:} optimal value and optimal solution of $\maximize \{c\cdot x : Ax\leq b, x\geq 0, x\in\Z^d\}$ for all $c \in \Xi$.

\noindent
{\bf Procedure:} for each $c \in \Xi$, do
\begin{enumerate}
 
 \item Use the identity (\ref{eq:i}) as
	necessary to enforce that all $v_{ij}$ in (\ref{eq:ab})
	satisfy $c \cdot v_{ij} < 0$.
  
\item Via the expansion
	formulas (\ref{eq:j}), find
	(\ref{eq:d}) by calculating the terms' coefficients.
	Proceed in decreasing order with respect to the degree of
	$t$.  This can be done because, for each series appearing in
	the expansion formulas (\ref{eq:j}), the terms of the series
	are given in decreasing order with respect to the degree of
	$t$.  

\item Continue calculating the terms of the expansion
	(\ref{eq:d}), in decreasing order with respect to the degree
	of $t$, until a degree $n$ of $t$ is found such that for some
	$\alpha\in\Z^d$, the coefficient of $y^{\alpha}t^n$
	is non-zero in the expansion (\ref{eq:d}).  
\item Return ``$n$'' as the optimal value of the integer program
	and return $\alpha$ as an optimal solution.
\end{enumerate}

We close this section by noticing that one nice feature of the digging
algorithm is if one needs to solve a family of integer programs where
only the cost vector $c$ is changing, then Equation (\ref{eq:ab}) can be
computed once and then apply the steps of the algorithm above for each
cost vector to obtain all the optimal values.

Given the polytope $P:= \{x \in \R^d: Ax \leq b, \, x \geq 0\}$, the
\emph{tangent cone} $K_v$ at a vertex $v$ of $P$ is the pointed
polyhedral cone defined by the inequalities of $P$ that turn into
equalities when evaluated at $v$. We will show in Chapter \ref{chapsoftware} that a
major practical bottleneck of the original Barvinok algorithm in
\cite{bar} is the fact that a polytope may have too many
vertices. Since originally one visits each vertex to compute a
rational function at each tangent cone, the result can be
costly. Therefore a natural idea for improving the digging algorithm
is to compute with a single tangent cone of the polytope and revisit
the above calculation for a smaller sum of rational functions. Let
vertex $v^*$ give the optimal value for the given linear programming
problem and we only deal with the tangent cone $K_{v^*}$. 
Suppose we have the following integer programming problem:

\[
\mbox{(IP) }\maximize \ c\cdot x \ \mbox{subject to} \quad x\in P \cap \Z^d,
\]

where $P:= \{x \in \R^d: Ax \leq b, \, x \geq 0\}$, $A \in \Z^{m \times d}$ and $b \in \Z^m$.

Then we have the following linear programming relaxation problem for the given integer programming problem:

\[
\mbox{(LP) }\maximize \ c\cdot x \ \mbox{subject to} \quad x \in P.
\]

One of the vertices of $P$ gives the optimal value for (LP) \cite{schrijver}. Let $V(P)$ be the vertex set of $P$ and $v \in V(P)$ be a vertex such that $c \cdot v$ is the optimal value for (LP).  Then, clearly, the tangent cone $K_v$ at $v$ contains $P$.  So, if we can find an integral point $x^* \in K_v$ such that $c \cdot x^* \geq c \cdot x, \, \forall x \in P$ and $x^* \in P \cap \Z^d$, then $x^*$ is an optimal solution for (IP).  The outline for the single cone digging algorithm is the following:

\noindent{\bf Algorithm:} (\emph{Single Cone Digging Algorithm}):

\noindent{\bf Input:} $A \in \Z^{m \times d}, \, b \in \Z^m, \, c \in \Z^d$.

\noindent{\bf Output:} optimal value and optimal solution of $\maximize \{c\cdot x : Ax\leq b, x\geq 0, x\in\Z^d\}$.

In the following steps, we replace $P$ by $K_v$ in the notation.

\begin{enumerate}
 
\item Compute a vertex $v$ of $P$ such that $c \cdot v =\maximize \{c\cdot x : Ax\leq b, x\geq 0 \}$.
 
\item Compute the tangent cone $K_v$ at $v$ and compute the short rational function (\ref{eq:ab}) encoding the lattice points inside $K_v$. 

 \item Use the identity (\ref{eq:i}) as
	necessary to enforce that all $v_{ij}$ in (\ref{eq:ab})
	satisfy $c \cdot v_{ij} < 0$.

\item Via the expansion
	formulas (\ref{eq:j}), find
	(\ref{eq:d}) by calculating the terms' coefficients.
	Proceed in decreasing order with respect to the degree of
	$t$.  This can be done because, for each series appearing in
	the expansion formulas (\ref{eq:j}), the terms of the series
	are given in decreasing order with respect to the degree of
	$t$.  

\item Continue calculating the terms of the expansion
	(\ref{eq:d}), in decreasing order with respect to the degree
	of $t$, until a degree $n$ of $t$ is found such that:
\begin{itemize}
\item for some $\alpha\in \Z^d$, the coefficient of $y^{\alpha}t^n$
	is non-zero in the expansion (\ref{eq:d}),
  
\item $A \alpha \leq b, \, \alpha \geq 0$.
\end{itemize}

\item Return ``$n$'' as the optimal value of the integer program
	and return $\alpha$ as an optimal solution.

\end{enumerate}

 From Table \ref{Running times} and Table \ref{IPtable5}, one can find that the single cone digging algorithm is very practical compared to the BBS algorithm and the original digging algorithm.  This algorithm is faster and more memory efficient than the original digging algorithm in practice, since the number of unimodular cones for the single cone digging algorithm is much less than the number of unimodular cones for the original digging algorithm.






\newpage
\pagestyle{myheadings} 
\markright{  \rm \normalsize CHAPTER 4. \hspace{0.5cm}
Experimental results: development of {\tt LattE}}
\chapter{Experimental results: development of {\tt LattE}}\label{chapsoftware}

\thispagestyle{myheadings} 

\section{{\tt LattE}'s implementation of Barvinok's algorithm} \label {software}

In this section, we go through the
steps of Barvinok's algorithm, showing how we implemented them in
{\tt LattE}. Barvinok's algorithm relies on two important new
ideas: the use of rational functions as efficient data
structures and the signed decompositions of cones into unimodular cones.

Let $P \subset \R^d$ be a rational convex polyhedron and let $f(P, z)$ be the multivariate generating function defined in (\ref{generating}). Let $v$ be a vertex of $P$.  
Then, the {\em supporting cone} $K(P, v)$ of $P$ at
$v$ is $K(P, v)= v+ \,\{u \in \real^d : v+\delta u \in P \ {\rm for \ all \ sufficiently \ small} \ \delta>0\}$.
Let $V(P)$ be the vertex set of $P$. One crucial component of
Barvinok's algorithm is the ability to distribute the computation on
the vertices of the polytope. This follows from the seminal theorem of
Brion \citep{brion88}:

\begin{theorem} \citep{brion88}
Let $P$ be a rational polyhedra and let $V(P)$ be the vertex set of $P$.  
Then, $$f(P, z)=\sum_{v \in V(P)}
f(K(P,v), z).$$
\end{theorem}

\begin{example}\label{example_1}
Consider the integral quadrilateral shown in Figure \ref{brionfig}. The
vertices are $V_1=(0,0)$, $V_2=(5,0)$, $V_3=(4,2)$, and $V_4=(0,2)$.

\begin{figure}[ht]
\begin{center}
\includegraphics[width=5 cm]{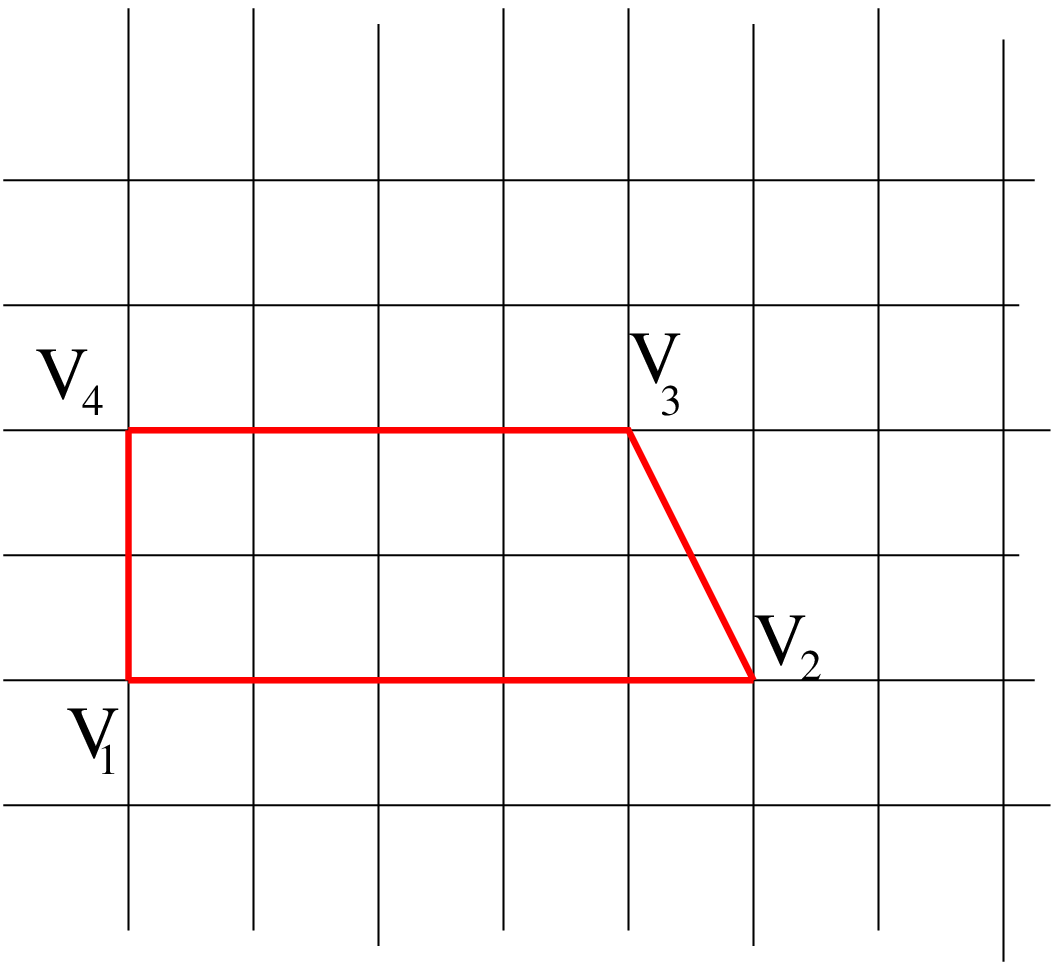} 
\caption{A quadrilateral in Example \ref{example_1}.}\label{brionfig}
\end{center} 
\end{figure}

We obtain four rational generation functions whose formulas are

$f(K_{V_1}, z)={\frac {1}{\left (1-z_1\right )\left (1-z_2\right )}}, \quad 
f(K_{V_2}, z)=\frac{({z_1}^{5}+{z_1}^{4}z_2)}{ (1-{z_1}^{-1}) (1-z_2^2z_1^{-1})}, $

$f(K_{V_3}, z)=\frac{({z_1}^{4}{z_2}+{z_1}^{4}{z_2}^{2})}{ (1-{z_1}^{-1})
(1-z_1z_2^{-2})},
 \quad  f(K_{V_4}, z)=\frac{z_2^{2}}{(1-{z_2}^{-1} )(1-z_1 )}.$

Indeed, the result of adding the rational functions is equal to the
polynomial

${z_1}^{5}+{z_1}^{4}z_2+{z_1}^{4}+{z_1}^{4}{z_2}^{2}+z_2{z_1}^{3}+{z_1}^{3}+{z_1}^{3}{z_2}^{2
}+z_2{z_1}^{2}+{z_1}^{2}+{z_1}^{2}{z_2}^{2}+z_1z_2+z_1+z_1{z_2}^{2}+{z_2}^{2}+z_2+1$. \eopf
\end{example}

In order to use Brion's theorem for counting lattice points in convex
polyhedra, we need to know how to compute the rational generating
function of convex rational pointed cones.  For polyhedral cones this
generating function is a rational function whose numerator and
denominator have a well-understood geometric meaning (see in
\citet[Chapter 4]{stanley} and in \citet[Corollary 4.6.8]{stanley2} for a clear explanation). We
already have a ``simple'' formula when the cone is a simple cone: Let
$\{u_1, u_2, \ldots, u_k\}$ be a set of linearly independent integral
vectors of $\real^d$, where $k \leq d$.  Let $K$ be a cone which is
generated by $\{u_1, u_2, \ldots, u_k\}$, in other words, $K=\{
\lambda_1 u_1 + \lambda_2 u_2 + \ldots + \lambda_k u_k,\mbox{ for some
} \lambda_i \geq 0 \mbox{ and } i=1,2, \ldots , k\}.$ Consider the
parallelepiped $S=\{ \lambda_1 u_1+ \lambda_2 u_2+ \ldots +\lambda_k
u_k, \ 0 \leq \lambda_i < 1 \mbox{, } i=1,2,\ldots, k\}.$

It is well-known \citep{stanley} that the generating function 
for the lattice points in $K$ equals

$$\sum_{\beta \in K \cap Z^d} z^{\beta} = (\sum_{\tau \in S \cap Z^d} 
z^{\tau})
\prod^k_{i=1} \frac{1}{1-z^{u_i}}. \quad (*) $$

Thus, to derive a formula for arbitrary pointed cones one could
decompose them into simple cones, via a triangulation, and then
apply the formula above and the inclusion-exclusion principle 
in \citet[Proposition 1.2]{stanley2}. Instead, Barvinok's idea is that
it is more efficient to further decompose each simple cone into
simple unimodular cones. A {\em unimodular} cone is a simple
cone with generators $\, \{u_1,\dots,u_k\} \,$ that form an integral
basis for the lattice $\, \real \{u_1,\dots,u_k\} \cap \Z^d $. Note
that in this case the numerator of the formula has a single monomial,
in other words, the parallelepiped has only one lattice point.

\subsection{Simple signed decompositions} \label{decomposition}

We now focus our attention on how the cone decomposition is done.  To
decompose a cone into simple cones the first step is to do a
triangulation ({\em triangulation} of a cone $K$ in dimension $d$ is a
collection of $d$-dimensional simple cones such that their union
is $K$, their interiors are disjoint, and any pair of them intersect
in a (possibly empty) common face).  There are efficient algorithms,
when the dimension is fixed, to carry a triangulation (see
\citet{auren,lee} for details). In {\tt LattE} we use the well-known
Delaunay triangulation which we compute via a convex hull calculation.
The idea is to ``lift'' the rays of the cone into a higher dimensional
paraboloid by adding a new coordinate which is the sum of the
squares of the other coordinates, take the lower convex hull of the
lifted points, and then ``project'' back those simple facets. We
use Fukuda's implementation in {\tt CDD} \citep{fukuda} of this
lift-and-project algorithm.  This is not the only choice of
triangulation, and definitely not the smallest one. In Section
\ref{theory} we discuss some situations when the choice of
triangulation in fact gives a better rational function.

In principle, one could at this point list the points of the
fundamental parallelepiped, for example, using a fast Hilbert bases
code such as {\tt 4ti2} \citep{Hemmecke:Hilbert} or {\tt NORMALIZ}
\citep{brunskoch}, and then use formula $(*)$ for a general simple
cone.  Theoretically this is bad because the number of lattice points
in the parallelepiped is exponentially large already for fixed
dimension. In practice, this can often be done and in some situations
 is useful.  Barvinok instead decomposes each simple cone as a
(signed) sum of simple {\em unimodular} cones. To be more formal,
for a set $A \subset \real^d$, the indicator function $[ A ]: \real^d
\rightarrow \real$ of $A$ is defined as $$[ A ] (x) = \left \{
\begin{array}{ll} 1 \mbox{ if }x \in A, \\ 0 \mbox{ if }x \not \in
A.\\
\end{array}\right .$$ 

We want to express the indicator function of a simple cone as an
integer linear combination of the indicator functions of unimodular simple
cones.  There is a nice valuation from the algebra of
indicator functions of polyhedra to the field of rational functions \citep{BarviPom}, and many of its properties can be
used in the calculation. For example, the valuation is zero when 
the polyhedron contains a line.

\begin{theorem} \citep[Theorem 3.1]{BarviPom} 
There is a valuation $f$ from the algebra of indicator functions of
rational polyhedra into the field of multivariate rational functions
such that for any polyhedron $P$, $f([P])=\sum_{\alpha \in P \cap
\Z^d} x^\alpha$.
\end{theorem}

Therefore once we have a unimodular cone decomposition, the rational
generating function of the original cone is a signed sum of
``short'' rational functions.  Next we focus on how to decompose
a simple cone into unimodular cones.

 Let $u_1, u_2, \ldots , u_d$ be linearly independent integral vectors
which generate a simple cone $K$.  We denote the \textsl{index} of
$K$ by ind$(K)$ which tells how far $K$ is from being unimodular.
That is, ind$(K) =|\det (u_1| u_2 | \ldots | u_d)|$ which is the
volume of the parallelepiped spanned by $u_1, u_2, \ldots , u_d$. It
is also equal to the number of lattice points inside the half-open
parallelepiped. $K$ is unimodular if and only if the index of $K$ is
$1$. Now we discuss how we implemented the following key result of
Barvinok:

\begin{theorem} \citep[Theorem 4.2]{BarviPom} \label{break}
Fix the dimension $d$.  Then, there exists a polynomial time algorithm
with a given rational polyhedral cone $K \subset \real^d$, which
computes unimodular cones $K_i$, $i \in I = \{1,
2, \ldots , l\}$, and numbers $\epsilon_i \in \{ -1, 1 \}$ such that
$$[K] = \sum_{i \in I} \epsilon_i [K_i].$$
\end{theorem}

Let $K$ be a rational pointed simple cone.  Consider the closed
parallelepiped $$\Gamma = \{ \alpha_1u_1+\alpha_2u_2 + \ldots +
\alpha_du_d : \mbox{ } |\alpha_j| \leq (\ind (K))^{- \frac{1}{d}}
\mbox{, } j = 1, 2, \ldots , d\}.$$ Note that this parallelepiped
$\Gamma$ is centrally symmetric and one can show that the volume of
$\Gamma$ is $2^d$. Minkowski's First Theorem \citep{schrijver}
guarantees that because $\Gamma \subset \real^d$ is a centrally
symmetric convex body with volume $\geq 2^d$, there exists a non-zero
lattice point $w$ inside of $\Gamma$. We will use $w$ to build the
decomposition.

We need to find $w$ explicitly. We take essentially the approach
suggested by \citet{Dyer}. We require a
subroutine that computes the shortest vector in a lattice. For fixed dimension this can be done in polynomial time using lattice basis reduction (this follows trivially from \citet[Corollary 6.4b][page 72]{schrijver}).  It is worth observing that when the dimension is not fixed the problem becomes NP-hard \citep{ajtai}.  
We use the basis reduction algorithm
of Lenstra, Lenstra, and Lov\'asz \citep{grolosch,schrijver} to find a short vector. Given $A$, an
integral $d \times d$ matrix whose columns generate a lattice, LLL's
algorithm outputs $A'$, a new $d \times d$ matrix, spanning the same
lattice generated by $A$. The column vectors of $A'$, $u'_1, u'_2,
..., u'_d$, are short and nearly orthogonal to each other, and each
$u'_i$ is an approximation of the shortest vector in the lattice, in
terms of Euclidean length. It is well-known \citep{schrijver} that
there exists a unique unimodular matrix $U$ such that $A U = A'$.

The method proposed in \citet{Dyer} to find $w$ is
the following: Let $A = (u_1|u_2| \ldots | u_d)$, where the $u_i$
are the rays of the simple cone we wish to decompose. Compute the
reduced basis of $A^{-1}$ using the LLL algorithm. Let $A'= (u_1'|u_2'| \ldots | u_d')$ be the
reduced basis of $A^{-1}$.  Dyer and Kannan observed that we can find
the smallest vector with respect to the $l^{\infty}$ norm by searching
over all linear integral combinations of the column vectors of $A'$
with small coefficients. We call this search the {\em enumeration step}.  
This enumeration step can be performed in polynomial time in fixed dimension.
We will briefly describe the enumeration step. 
First we introduce some notation.  
Let $u_1', \dots , u_d'$ be linearly independent integral vectors in $\Z^d$.
Let $\| \cdot \|_2$ be the $l^2$ norm and let $\| \cdot \|_\infty$ be the 
infinity norm.

We will need to recall the Gram-Schmidt process that computes a set of orthogonal 
vectors $u^*_j$, $1\leq j \leq d$,
from independent vectors $u_j'$, $1\leq j \leq d$.  
In particular we need some values from this process.  The vectors 
$u^*_j$, and real numbers 
$\mu_{j, k}\mbox{, } 1 \leq k < j \leq d$ are computed from $u_j'$
by the recursive process: $$u^*_1=u_1'$$
$$u^*_j=u_j'-\sum_{k=1}^{j-1} \mu_{j, k} u^*_k\mbox{, } \, \, 2 \leq j \leq d$$
$$\mu_{j, k}=\frac{u_j' \cdot u^*_k}{ \| u^*_k\|_2 ^2}.$$
Letting $w_i := u_i^* /\|u_i^*\|_2$, there exists real numbers $u_i(j)$ 
such that
\begin{equation}\label{eqkannan1}
u_i' = \sum_{j=1}^d u_i(j) w_j.
\end{equation}
Note that $u_i(j) = \mu_{i, j} \|u_j^*\|_2$ 
for $1 \leq k < j \leq d$ and $u_i(i) = \|u^*_i\|_2$. 
These $u_i(j)$ will be used below.
Let $L(u_1', \dots , u_d')$ be the lattice generated by 
$u_1', \dots , u_d'$.  Then 
we denote $L_j(u_1', \dots , u_d')$ be the projection of $L(u_1', \dots , u_d')$ 
orthogonal to the vector space $V_j$ spanned by $u_1', \dots , u_j'$. 

Now we are ready to describe the process of the enumeration step.
Let $\lambda$ be a shortest vector in the
lattice spanned by $A'$ with respect to the $l^{\infty}$ norm. 
Then we can write $\lambda$ as an integral linear 
combination of columns of $A'$.  Let $\lambda = \sum_{i=1}^d \alpha_i u'_i$,
where $\alpha = (\alpha_1, \dots , \alpha_d) \in \Z^d$. 
The goal is to find some finite set $T \subset \Z^d$ such that $\alpha \in T$ and
the cardinality of $T$ is polynomial size in fixed dimension. $T$ will be
contained inside a certain parallelepiped.  Then we can
search $\lambda$ by enumerating all lattice points inside $T$.

Suppose $A' = (u_1'| \dots | u_d')$ form the reduced basis obtained by LLL algorithm.
Let $m:= \min \{j: u_j(j) \geq u_1(1) \} - 1$. 
Now we will apply the inequalities 
\begin{equation}\label{eqkannan3}
\| x \|_\infty \leq \| x \|_2,
\end{equation}

\begin{equation}\label{eqkannan4}
\| x \|_2 \leq \sqrt{d} \| x \|_\infty.
\end{equation}
We are going to prove that a shortest vector of $L(u_1', \dots , u_m')$ is a shortest
vector of $L(u_1', \dots , u_d')$ with respect to the $l^\infty$ norm.  Any 
vector in $L(u_1', \dots , u_d') \backslash L(u_1', \dots , u_m')$ must have
$l^2$ norm at least $u_1(1)$.  Since $u_1(1) = \|u_1'\|_2$ it must have $l^\infty$ norm 
at least $u_1(1)$ which is at least the $l^\infty$ norm of $u_1'$ by 
(\ref{eqkannan3}).  

We will show how to construct $T$.
Let $y = \sum_{i=1}^d \alpha_i u_i'$ be a candidate for a shortest vector
with respect to the $l^\infty$ norm.
Applying the fact that $|\alpha_i u_i(i)| \leq \|y\|_2$ 
(using the same trick as on page 423, \cite{kannan}),
we have $|\alpha_i u_i(i)| / \sqrt{d} \leq  \|y\|_2 / \sqrt{d} 
\leq  \|u_1'\|_\infty $
for any candidate vector for a shortest vector with respect to the 
$l^\infty$ norm. Therefore, we have
$$|\alpha_i|u_i(i)  / \sqrt{d} \leq
 \|u_1' \|_\infty \leq u_1(1).$$  From this
$$|\alpha_i| \leq  \sqrt{d} u_1(1)/u_i(i) \mbox{ for } i = 1, \dots, m, $$
which defines a parallelepiped in the variables $\alpha_i$ such that, 
$$Q := \{ \alpha \in \R^d: -\sqrt{d} u_1(1)/u_i(i) \leq \alpha_i \leq \sqrt{d} u_1(1)/u_i(i) \mbox{ for } i = 1, \dots , d\}.$$
Finally we set $T := Q \cap \Z^d$.

Now we are going to show that $Q$ contains polynomially many
lattice points. For each $\alpha_i$, there exist at most 
$1+2 \sqrt{d}u_1(1)/u_i(i)$ candidates.  So the total number of candidates is 
$$\prod_{i=1}^m (1+2 \sqrt{d} u_1(1)/u_i(i)).$$  
With the fact that $u_1(1) \geq u_i(i)$ (by the definition of $m$), we have 
$$\prod_{i=1}^m (1+2 \sqrt{d} u(1)_1/u_i(i)) \leq 3^m d^{m/2} \prod_{i=1}^m(u_1(1)/u_i(i)).$$
We derive the following from Minkowski's theorem 
$$u_1(1)^m \leq (2m)^{m/2} \det(L(u_1', \dots, u_m')),$$
$$\det(L(u_1', \dots, u_m')) = \prod_{i=1}^m u_i(i).$$ 
Therefore, we have $\prod_{i=1}^m(u_1(1)/u_i(i)) 
\leq (2m)^{m/2}$.  This implies that 
$$\prod_{i=1}^m (1+2 \sqrt{d} u(1)_1/u_i(i))
\leq (3d)^{d},$$ which is a constant if we fix $d$. With this method, we 
can compute a shortest vector $\lambda$ with respect to the $l^\infty$ norm
in polynomial time in fixed dimension by the enumeration step.

After we compute $\lambda$ in polynomial time in fixed dimension,
we know that
there exists a unique unimodular matrix $U$ such that $A' = A^{-1}U$.
Minkowski's theorem for the  $l^{\infty}$ norm implies 
that for the non-singular matrix $A'$, there exists a non-zero integral vector
$z$ such that $\|\lambda \|_\infty =\|A'z\|_{\infty} \leq |\det(A')|^{1/d}$.
See statement 23 in page 81 
in \citet{schrijver}. We can set

$$\|\lambda \|_{\infty} \leq |\det(A')|^{1/d} = |\det(A^{-1} U)|^{1/d}
=
|\det(A^{-1}) \det(U)|^{1/d}
$$

$$ = |\det(A^{-1})|^{1/d} = |\det(A)|^{-1/d} = |\ind(K)|^{-1/d}.$$

Since $A^{-1}$ and $A'$ span the same lattice, there exists an
integral vector $w \in \real^d$ such that $\lambda = A^{-1} w$. Then,
we have $$w = A \lambda.$$ Note that $w$ is
a non-zero integral vector which is a linear integer combination of
the generators $u_i$ of the cone $K$ {\em with possibly negative
coefficients}, and with coefficients at most $|\ind(K)|^{-1/d}$.
Therefore, we have found a non-zero integral vector $w \in
\Gamma$. In {\tt LattE}, we try to avoid the enumeration step
because it is very costly. Instead, we choose $\lambda$ to be the
shortest of the columns in $A'$. This may not be the smallest vector,
but for practical purposes, it often decreases the index $| \ind(K)|$
just like for a shortest vector. Experimentally
we have observed that we rarely use the enumeration step.

In the next step of the algorithm, for $i = 1, 2, \ldots , d$, we set
$$K_i = \co\{ u_1, u_2, \ldots , u_{i-1}, w , u_{i+1}, \ldots ,
u_d \}.$$ Now, we have to show that for each $i$, $\ind(K_i)$ is
smaller than $\ind(K)$.  Let $w = \sum_{i = 1}^d \alpha_i u_i$.
Then, we have

\begin{eqnarray*}
\ind(K_i)  &=& |\det((u_1|u_2|\ldots | u_{i-1}| w |u_{i+1}|
\ldots|u_d))| \,  \\
 &=& |\alpha_i||\det((u_1|u_2|\ldots | u_{i-1}| u_i |u_{i+1}| \ldots
|u_d))|  \,  \\
&=&|\alpha_i|\ind(K) \leq (\ind(K))^{\frac{d-1}{d}}  \, . \\ 
\end{eqnarray*}





There is one more technical condition that $w$ needs to satisfy. This
is that $w$ and $u_1,\dots,u_d$ belong to an open half-space \citep[Lemma 5.2]{bar}. This is easy to achieve as either the $w$ we
found or $-w$ satisfy this condition.  We can now decompose the
original cone $K$ into cones $K_i$ for $i= 1, 2, \ldots , d$, of
smaller index, $[K]=\sum \pm [K_i]$. This sum of indicator functions
carries signs which depend on the position of $w$ with respect to the
interior or exterior of $K$. We iterate this process until $K_i$
becomes a unimodular cone for $i = 1, 2, \ldots , d$. For implementing
Barvinok's decomposition of cones, we use the package {\tt NTL} by \citet{shoup} to compute the reduced basis of a cone and
to compute with matrices and determinants. All our calculations were
done in exact long integer arithmetic using the routines integrated in {\tt NTL}.  Here is the
pseudo-code of the algorithm and an example.

\begin{Algorithm}(Barvinok's Decomposition of a Simple Cone)

\noindent
\textbf{Input:}  A simple cone $K = \co \{u_1, u_2, 
\ldots , u_d \}$ given by its generators.

\noindent
\textbf{Output:} A list of unimodular cones and numbers $\epsilon_i$ 
as in Theorem \ref{break}.
\begin{tabbing}
\quad \= \quad \= \quad \= \quad \kill
Set two queues Uni and NonUni.\\
\keyw{if} $K$ is unimodular\\
\> \keyw{then} Uni $= \mbox{Uni } \cup \{K\}$.\\  
\keyw{else} NonUni $= \mbox{NonUni } \cup \{K\}$.\\
\keyw{while} NonUni is not empty \keyw{do}\\
\> Take a cone $K \in $ NonUni and set $A = (u_1, u_2, \ldots 
, u_d)$\\
\> \> to be a matrix  whose columns are the rays of $K$.\\
\> Compute the smallest vector $\lambda$ in the lattice,\\
\> \> with respect to  $l^{\infty}$, which is spanned by  the column vectors of $A^{-1}$.\\
\> Find a non-zero integral vector $z$ such that $\lambda = A^{-1}z$.\\
\> \keyw{if} vectors $z, u_1, u_2, \ldots , u_d$ are in an open half plane\\
\> \> \keyw{then} set $z := z$.\\
\> \keyw{else} set $z := -z$.\\
\> \keyw{for} $i = 1, 2, \ldots , d$ \keyw{do}\\
\> \> set $K_i = \co \{ u_1, \ldots , u_{i-1}, z, u_{i+1}, \ldots,
u_d \}$\\
\> \>  and set $A_i = (u_1, \ldots , u_{i-1}, z, u_{i+1}, \ldots,
u_d )$.\\
\>\keyw{for} $i = 1, 2, \ldots , d$ \keyw{do}\\
\> \> \keyw{if} $\det(A_i)$ and $\det(A)$ have the same sign\\
\> \> \> \keyw{then} assign $\epsilon_{K_i} = \epsilon_K$.\\
\> \> \keyw{else} $\epsilon_{K_i} = -\epsilon_K$.  \\
\> \keyw{for} $i = 1, 2, \ldots , d$ \keyw{do}\\
\> \>  \keyw{if} $K_i$ is unimodular\\
\> \> \> \keyw{then}  Uni = Uni $\cup \{K_i\}$. \\
\> \> \keyw{else} NonUni = NonUni $\cup \{K_i\}$.\\
\keyw{return} all elements in Uni.
\end{tabbing}
\end{Algorithm}

\begin{figure}[htpb]
\begin{center}
\includegraphics[width = 10cm]{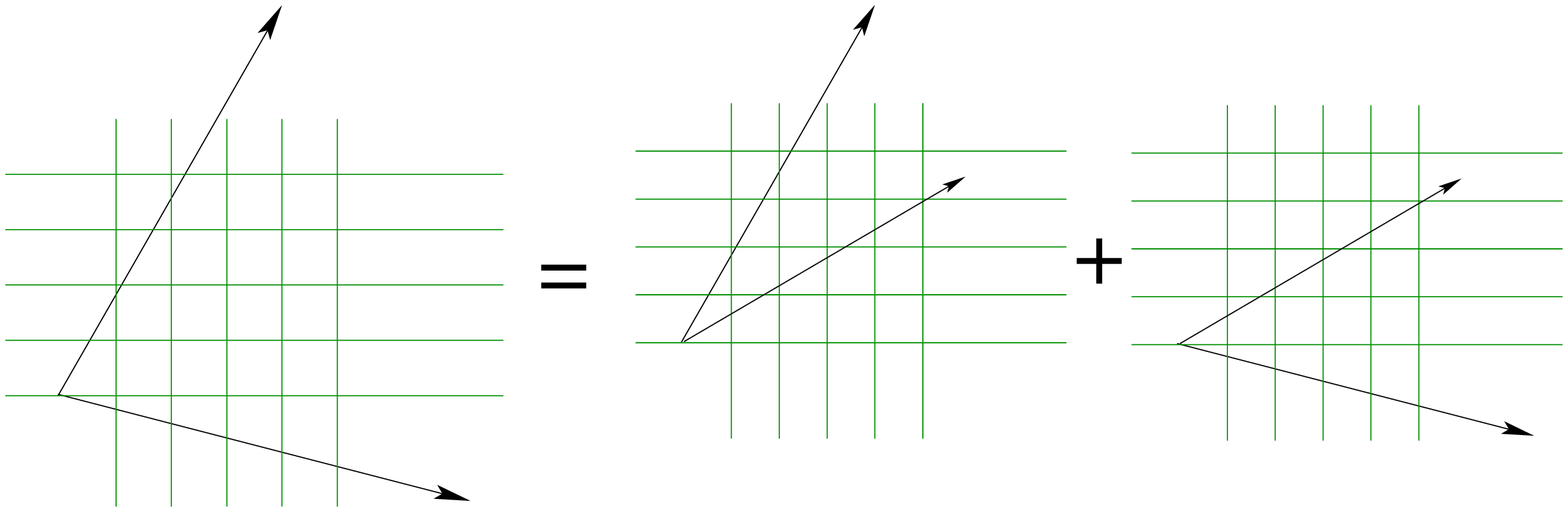}
\caption{Contribution of lower dimensional cones} \label{contribution}
\end{center}
\end{figure}

It is very important to remark that, in principle, 
one also needs to keep track of lower dimensional cones present in the
decomposition for the purpose of writing the inclusion-exclusion
formula of the generating function $f(K, z)$. For example in Figure
\ref{contribution} we have counted a ray twice, and thus it needs
to be removed.

But this is actually not necessary thanks to a {\em Brion's
polarization trick} \citep[Remark 4.3]{BarviPom}: Let $K^*$ be
the dual cone to $K$. Apply the iterative procedure above to $K^*$
instead of $K$, ignoring the lower dimensional cones. This can be done 
because once we polarize the result back, the
contribution of the lower dimensional cones is zero with respect to
the valuation that assigns to an indicator function its generating
function counting the lattice points \citep[Corollary 2.8]{BarviPom}. In the current implementation of {\tt LattE} we do
the following:

\begin{enumerate}
\item Find the vertices of the polytope and their defining supporting cones. 

\item Compute the polar cone to each of the cones.

\item Apply the Barvinok decomposition to each of the polars.

\item Polarize back the cones to obtain a decomposition, into
full-dimensional unimodular cones, of the original supporting cones.

\item Recover the generating function of each cone and, by Brion's
theorem, of the whole polytope.
\end{enumerate}

Here is an example of how we carry out the decomposition. 

\noindent
\begin{example}
Let $K$ be a cone generated by $(2, 7)^T$ and $(1, 0)^T$.  Let
$$A = \left( \begin{array}{cc}
2& 1\\
7& 0\\
\end{array}\right).$$  Then, we have $\det(A) = -7$ and $$A^{-1} = \left( \begin{array}{cc}
0 & \frac{1}{7}\\
1 & \frac{-2}{7}\\
\end{array}\right).$$
The reduced basis $A'$ of $A^{-1}$ and the unimodular matrix $U$ for 
the
transformation from $A^{-1} $ to $A'$ are:
$A' = \left( \begin{array}{cc}
\frac{1}{7} &\frac{3}{7}\\
\frac{-2}{7} & \frac{1}{7}\\
\end{array}\right),$
and 
$U = \left( \begin{array}{cc} 
0 & 1\\
1 & 3\\
\end{array}\right).$
By enumerating the column vectors, we can verify that $(\frac{-2}{7},
\frac{1}{7})^T$ is the smallest vector with respect to $l^{\infty}$ in
the lattice generated by the column vectors of $A^{-1}$.  So, we have $z =
(1, 0)^T$.  Then, we have two cones:
$$\left( \begin{array}{cc}
2 & 0\\
7 & 1\\
\end{array}\right) \mbox{ and } \left( \begin{array}{cc} 
0 & 1\\
1 & 0\\
\end{array}\right).$$
The second cone is unimodular of index $-1$ which is the 
same sign as the determinant of $A$.  Thus, Uni $= \mbox{ Uni } \cup \{\left( 
\begin{array}{cc} 
0 & 1\\
1 & 0\\
\end{array} \right)\},$ and assign to it $\epsilon = 1$.  The first cone 
has 
determinant $2$.  So, we assign $\epsilon = -1$.  Since the first cone is 
not unimodular, we have NonUni = NonUni $\cup \{\left( 
\begin{array}{cc} 
2 & 0\\
7 & 1\\
\end{array} \right)\}.$  Set $$A =\left( \begin{array}{cc}
2 & 0\\
7 & 1\\
\end{array}\right).$$
Then, we have $\det(A) = 2$ and $$A^{-1} = \left( \begin{array}{cc}
\frac{1}{2} & 0\\
\frac{-7}{2} & 1\\
\end{array}\right) \mbox{, } A' = \left( \begin{array}{cc}  
\frac{1}{2} & \frac{1}{2}\\
\frac{-1}{2} & \frac{1}{2}\\
\end{array}\right) \mbox{ and } U = \left( \begin{array}{cc}
1 & 1\\
3 & 4\\
\end{array}\right).$$ Since $\lambda = (\frac{1}{2}, \frac{-1}{2})^T$ is
the smallest vector with respect to $l^{\infty}$, we have $z = (1,
3)^T$.  So, we get two cones: $$\left( \begin{array}{cc}
2 & 1\\
7 & 3\\
\end{array}\right) \mbox{ and } \left( \begin{array}{cc} 
1 & 0\\
3 & 1\\
\end{array}\right).$$ 
The first matrix has negative determinant which is not the same sign
as the determinant of its parent matrix $A$.  Since $\epsilon_{A} =
-1$, we assign to the first cone $\epsilon = 1$ and the second one has
positive determinant, so we assign to it $\epsilon = 1$.  Since both
of them are unimodular, we take them into Uni and since NonUni is
empty, we end while loop and print all elements in Uni.

\noindent
This gives a full decomposition:
$$ \co \{ \left( \begin{array}{c}
2\\
7\\
\end{array}\right) ,
\left( \begin{array}{c}
1\\
0\\
\end{array}\right) \}$$
$$ = \ominus
\co\{ \left( \begin{array}{c}
1\\
3\\
\end{array}\right),  
\left( \begin{array}{c}
0\\
1\\
\end{array}\right) \} \oplus \co\{ \left( \begin{array}{c}
0\\
1\\
\end{array}\right), 
\left( \begin{array}{c}
1\\
0\\
\end{array}\right) \} \oplus \co\{ \left( \begin{array}{c}
2\\
7\\
\end{array}\right), 
\left( \begin{array}{c}
1\\
3\\
\end{array}\right) \}.$$ \eopf

\begin{figure}[htpb]
\begin{center}
\includegraphics[width = 10cm]{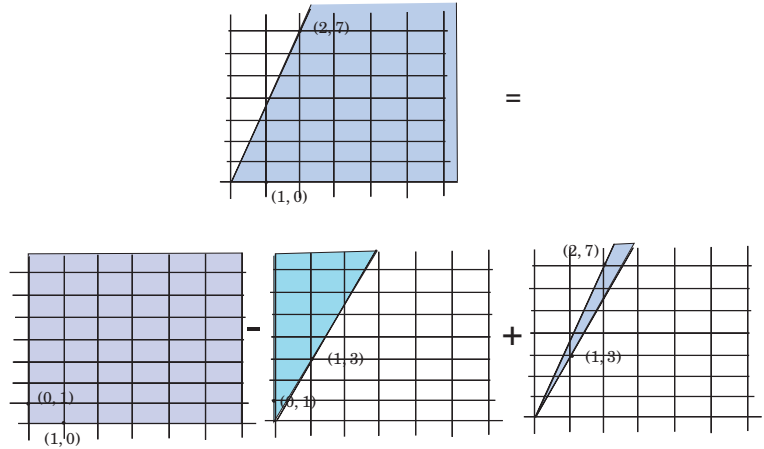}
\caption{Example of Barvinok's decomposition}
\end{center}
\end{figure}

\end{example}

From the previous example, we notice that the determinant of each cone
gets much smaller in each step. This is not an accident as Theorem
\ref{break} guarantees that the cardinality of the index set $I$ of
cones in the decomposition is bounded polynomially in terms of the
determinant of the input matrix.  We have looked experimentally at how
many levels of iteration are necessary to carry out the
decomposition. We observed experimentally that it often grows linearly with
the dimension. We tested two kinds of instances. We used random square
matrices whose entries are between 0 and 9, thinking of their columns
as the generators of a cone centered at the origin. We tested from $2
\times 2$ matrices all the way to $8 \times 8$ matrices, and we tested
fifteen random square matrices for each dimension. We show the results
in Table \ref{randomcones}. For computation, we used a $1$ GHz Pentium
PC machine running Red Hat Linux.

\begin{table} 
\begin{center}
\begin{tabular}{|c|c|c|c|c|} \hline
Dimension & Height of tree  & \# of cones & $|$ determinant$| $& Time 
(seconds)\\ \hline 
2 &  1.33 & 2.53 & 11.53 & 0 \\ \hline
3 & 2.87 & 12.47 & 55.73 & 0.005\\ \hline
4 & 3.87 & 65.67  & 274.667& 0.153  \\ \hline
5 & 5.87 & 859.4 & 3875.87 & 0.25 \\ \hline
6 & 7.47 & 10308 & 19310.4 & 3.67 \\ \hline
7 & 8.53 & 91029.4 & 72986.3 & 41.61 \\ \hline 
8 & 10.67  & 2482647.533  & 1133094.733  & 2554.478  \\ \hline

\end{tabular}
\caption{Averages of 15 random matrices for computational experiences} 
\label{randomcones}
\end{center}
\end{table}

The second set of examples comes from the Birkhoff polytope $B_n$ of
doubly stochastic matrices \citep{schrijver}. Each vertex of
the polytope is a permutation matrix which is a $0/1$ matrix whose 
column sums and row sums are all $1$ \citep{schrijver}. We decompose
the cone with vertex at the origin and whose rays are the $n!$ 
permutation matrices. The results are reported in Table \ref{birkhoff}.


\begin{table}
\begin{center}
\begin{tabular}{|c|c|c|c|} \hline
Dimension & \# of vertices & \# of unimodular cones at a vertex cone &
Time (seconds)\\ \hline
$B_3=4$ & 6 & 3 & 0.05\\ \hline
$B_4=9$ & 24 & 16 & 0.15  \\ \hline
$B_5=16$ & 120 & 125 & 0.5 \\ \hline
$B_6=25$ & 720 & 1296 & 7.8 \\ \hline

\end{tabular}
\caption{The numbers of unimodular cones for the Birkhoff polytopes} 
\label{birkhoff}
\end{center}
\end{table} 

\subsection{From cones to rational functions and counting}
\label{limit}

Once we decompose all cones into simple unimodular cones, it is
easy to find the generating function attached to the $i$th cone
$K_i$. In the denominator there is a product of binomials of the form
$(1-z^{B_{ij}})$ where $B_{ij}$ is the $j$th ray of the cone $K_i$.
Thus the denominator is the polynomial $\prod (1-z^{B_{ij}})$. How
about the numerator? The cone $K_i$ is unimodular, thus it must have a
single monomial $z^{A_i}$, corresponding to the unique lattice point
inside the fundamental parallelepiped of $K_i$. Remember that the
vertex of $K_i$ is one of the vertices of our input polytope. If that
vertex $v$ has all integer coordinates then $A_i=v$, or else $v$ can
be written as a linear combination $\sum \lambda_j B_{ij}$ where all
the $\lambda_i$ are rational numbers and can be found by solving a
system of equations (remember the $B_{ij}$ form a vector space basis
for $\real^d$). The unique lattice point inside the parallelepiped of
the cone $K_i$ is simply $\sum \lceil \lambda_j \rceil B_{ij}$ \citep[Lemma 4.1]{BarviPom}.

Brion's theorem says the sum of the rational functions coming from the
unimodular cones at the vertices is a polynomial with one monomial per
lattice point inside the input polytope. One might think that to
compute the number of lattice points inside of a given convex
polyhedron, one could directly substitute the value of $1$ at each of
the variables. Unfortunately, $(1,1,\dots,1)$ is a singularity of all
the rational functions. Instead we discuss the method used in {\tt
LattE} to compute this value, which is different from that presented
by Barvinok \citep{BarviPom}. The typical generating function of
lattice points inside a unimodular cone forms:

$$ E[i] \frac{z^{A_i}}{\prod (1-z^{B_{ij}})},$$

where $z^a$ is monomial in $d$ variables, each $A_i$ (cone vertex)
and $B_{ij}$ (a generator of cone $i$) are integer vectors of length
$d, i$ ranges over all cones given, $j$ ranges over the generators of
cone $i$, and $E[i]$ is 1 or -1. Adding these rational functions and
simplifying would yield the polynomial function of the lattice point
of the polytope. Now this is practically impossible as the number of
monomials is too large. But calculating the number of monomials in
this polynomial is equivalent to evaluating the limit as $z_i$ goes to
$1$ for all $i$. We begin by finding an integer vector $\lambda$ and making the
substitution $z_i \rightarrow t^{\lambda_i}$. This is with the
intention of obtaining a univariate polynomial. To do this, $\lambda$
must be picked such that there is no zero denominator in any cone
expression, i.e. no dot product of $\lambda$ with a $B_{ij}$ can be
zero. Barvinok showed that such a $\lambda$ can be picked in polynomial time
by choosing points on the moment curve. Unfortunately, this method
yields large values in the entries of $\lambda$. Instead we try random
vectors with small integer entries, allowing small increments if
necessary, until we find  $\lambda$. Since we are essentially trying
to avoid a measure zero set, this process terminates very quickly in
practice.

After substitution, we have expressions of the form $\pm
t^{N_i}/\prod (1-t^{D_{ij}})$, where $N_i$ and $D_{ij}$ are
integers. Notice that this substitution followed by summing these
expressions yields the same polynomial as would result from first
summing and then substituting. This follows from the fact that we can
take Laurent series expansions, and the sum of Laurent series is equal
to the Laurent series of the sum of the original expressions.

Also, note that we have the following identity:
$$\sum_{\alpha \in P \cap \Z^d} z^{\alpha} = \sum_{i=1}^{\# \ of \
cones} E[i] \frac{z^{A_i}}{\prod(1-z^{B_{ij}})}.$$ After substitution
we have the following univariate (Laurent) polynomial such that:
$$\sum_{\alpha \in P \cap \Z^d} t^{\sum_{i=1}^d \lambda_i \alpha_i}=\sum^{\# \ of \ cones}_{i=1} E[i] \frac{t^{N_i}}{\prod (1-t^{D_{ij}})}.$$

With the purpose of avoiding large exponents in the numerators, we
factor out a power of $t$, say $t^c$. Now we need to evaluate the sum
of these expressions at $t=1$, but we cannot evaluate these
expressions directly at $t=1$ because each has a pole there. Consider
the Laurent expansion of the sum of these expressions about $t=1$.
The expansion must evaluate at $t=1$ to the finite number
$\sum_{\alpha \in P \cap \Z^d} 1$. It is a Taylor expansion and its
value at $t=1$ is simply the constant coefficient.  If we expand each
expression about $t=1$ individually and add them up, it will yield the
same result as adding the expressions and then expanding (again the
sum of Laurent expansions is the Laurent expansion of the sum of the
expressions). Thus, to obtain the constant coefficient of the sum, we
add up the constant coefficients of the expansions about $t=1$ of each
summand.  Computationally, this is accomplished by substituting
$t=s+1$ and expanding about $s=0$ via a polynomial division.  Summing
up the constant coefficients with proper accounting for $E[i]$ and
proper decimal accuracy yields the desired result: the number of
lattice points in the polytope. Before the substitution $t=s+1$ we
rewrite each rational function in the sum (recall $t^c$ was factored
to keep exponents small);

$$\sum E[i] \frac{t^{N_i-c}}{\prod (1-t^{D_{ij}})}=\sum E'[i]
\frac{t^{N'_i}}{\prod (t^{D'_{ij}}-1)},$$ involves in such a way that
$D'_{ij}>0$ for all $i,j$. This requires that the powers of $t$ at each
numerator to be modified, and the sign $E[i]$ is also adjusted to $E'[i]$. Then
the substitution $t=s+1$ yields
$$ \sum E'[i] \frac{(1+s)^{N'_i}}{\prod ((1+s)^{D'_{ij}}-1)},$$ 
where it is evident that, in each summand, the pole $s=0$ has an order
equal to the number of factors in the denominator. This is the same as
the number of rays in the corresponding cone and we denote this number
by $d$.

Thus the summand for cone $i$ can be rewritten as
$E'[i]s^{-d}P_i(s)/Q_i(s)$ where $P_i(s)=(1+s)^{N_i}$ and
$Q_i(s)=\prod^d ((1+s)^{D'_{ij}}-1)/s)$.  $P_i(s)/Q_i(s)$ is a
Taylor polynomial whose $s^d$ coefficient is the contribution we are
looking for (after accounting for the sign $E'[i]$ of course). The
coefficients of the quotient $P_i(s)/Q_i(s)$ can be obtained
recursively as follows: Let $Q_i(s) =b_0 + b_1 s + b_2 s^2 + \ldots $
and $P_i(s) = a_0 + a_1 s+ a_2 s^2 + \ldots $ and let
$\frac{P_i(s)}{Q_i(s)} = c_0 + c_1 s + c_2 s^2 +\ldots $. Therefore,
we want to obtain $c_d$ which is the coefficient of the constant term
of $P_i/Q_i$.  So, how do we obtain $c_d$ from $Q_i(s)$ and $P_i(s)$?
We obtain this by the following recurrence relation:

$$ c_0 = \frac{a_0}{b_0}, $$
$$ c_k = \frac{1}{b_0}(a_k - b_1 c_{k-1} - b_2c_{k-2} - \ldots - b_k
c_0) \mbox{ for } k = 1, 2, \ldots.$$

In order to obtain $c_d$, only the coefficients $a_0$, $a_1, \ldots ,
a_d$ and $b_0$, $b_1, \ldots , b_d$ are required.

\begin{example}
(A triangle). Let us consider three points in $2$
dimensions such that $V_1 = (0, 1)$, $V_2 = (1, 0)$, and $V_3 = (0,
0)$.  Then, the convex hull of $V_1$, $V_2$, and $V_3$ is a triangle
in $2$ dimensions.  We want to compute the number of lattice points by
using the residue theorem.  Let $K_i$ be the vertex cone at $V_i$ for
$i = 1, 2, 3$.  Then, we have the rational functions:

$$f(K_1, (x, y)) = \frac{y}{(1-y^{-1})(1 - x y^{-1})}, f(K_2, (x, y)) = \frac{x}{(1-x^{-1})(1 - x^{-1} y)},$$ $$ f(K_3, (x, y)) = \frac{1}{(1-x)(1-y)}.$$

We choose a vector $\lambda$ such that the inner products of $\lambda$
and the generators of $K_i$ are not equal to zero.  We choose $\lambda
= (1, -1)$ in this example. Then, reduce multivariate to univariate
with $\lambda$, so that we have:

$$f(K_1, t) = \frac{t^{-1}}{(1-t)(1 - t^{2})}, f(K_2, t) = \frac{t}{(1-t^{-1})(1 - t^{-2})}, f(K_3, t) = \frac{1}{(1-t)(1-t^{-1})}.$$ 

We want to have all the denominators to have
positive exponents.  We simplify them in order to eliminate negative
exponents in the denominators with simple algebra.  Then, we have:

$$f(K_1, t) = \frac{t^{-1}}{(1-t)(1 - t^{2})}, f(K_2, t) =
\frac{t^4}{(1-t)(1 - t^{2})}, f(K_3, t) = \frac{-t}{(1-t)(1-t)}.$$ 

We factor out $t^{-1}$ from each rational function, so that we obtain:

$$f(K_1, t) = \frac{1}{(1-t)(1 - t^{2})}, f(K_2, t) =
\frac{t^5}{(1-t)(1 - t^{2})}, f(K_3, t) = \frac{-t^2}{(1-t)(1-t)}.$$  

We substitute $t = s+1$ and simplify them to the form 
$\frac{P(s)}{s^d Q(s)}$:

$$f(K_1, s) = \frac{1}{s^2(2 + s)}, f(K_2, s) = \frac{1+5s+10s^2+10s^3+5s^4+s^5}{s^2(2 + s)},$$ $$ f(K_3, s) = \frac{-(1+2s+s^2)}{s^2}.$$

Now we use the recurrence relation to obtain the coefficient of
the constant terms.  Then, for $f(K_1, s)$, we have $c_2 = \frac{1}{8}$.
For $f(K_2, s)$, we have $c_2 = \frac{31}{8}$.  For $f(K_3, s)$, we have
$c_2 = - 1$.  Thus, if we sum up all these coefficients, we have $3$,
which is the number of lattice points in this triangle. \eopf

\end{example}

{\tt LattE} produces the sum of rational functions which converges to the generating function of the lattice points of
an input polytope. This generating function is a multivariate
polynomial of finite degree. As we saw in Subsection \ref{limit}
it is possible to count the number of lattice points without expanding
the rational functions into the sum of monomials. Suppose that instead
of wanting to know the number of lattice points we simply wish to {\em
decide} whether there is one lattice point inside the polytope or
not. The integer feasibility problem is an important and difficult
problem \citep{aardaletal1,schrijver}. Obviously, one can simply
compute the residues and then if the number of lattice points is
non-zero, clearly, the polytope has lattice points. 

Before we end our description of {\tt LattE}, we must comment on how
we deal with polytopes that are not full-dimensional (e.g.
transportation polytopes).  Given the lower-dimensional polytope
$P=\{x\in\R^n:Ax=a, Bx\leq b\}$ with the $d\times n$ matrix $A$ of
full row-rank, we will use the equations to transform $P$ into a
polytope $Q=\{x\in\R^{n-d}:Cx\leq c\}$ in fewer variables, whose
integer points are in one-to-one correspondence to the integer points
of $P$. This second polytope will be the input to the main part of
{\tt LattE}. The main idea of this transformation is to find the
general integer solution $x=x_0+\sum_{i=1}^{n-d} \lambda_i g_i$ to
$Ax=a$ and to substitute it into the inequalities $Bx\leq b$, giving a
new system $Cx\leq c$ in $n-d$ variables
$\lambda_1,\ldots,\lambda_{n-d}$.

It is known that the general integer solution $Ax=a$ can be found
via the Hermite normal form $H=(R|0)$ of $A$ \citep{schrijver}.
Here, $R$ is a lower-triangular matrix and $H=AU$ for some
unimodular matrix $U$. Moreover, as $A$ is supposed to have full
row-rank, $R$ is a non-singular $d\times d$ matrix. Let $U_1$ be
the matrix consist of the first $d$ columns of $U$ and $U_2$
consisting of the remaining $n-d$ columns of $U$. Now we have
$AU_1=R$ and $AU_2=0$ and the columns of $U_2$ give the generators
$\{g_1,\ldots,g_{n-d}\}$ of the integer null-space of $A$. Thus,
it remains to determine a special integer solution $x_0$ to
$Ax=a$.

To do this, first find an integer solution $y_0$ to $Hy=(R|0)y=a$,
which is easy due to the triangular structure of $R$. With
$x_0=Uy_0$, we get $Ax_0=AUy_0=Hy_0=a$ and have found all pieces
of the general integer solution $x=x_0+\sum_{i=1}^{n-d} \lambda_i
g_i$ to $\{x\in\Z^n:Ax=a\}$.

\section{Computational experience and performance for counting} \label{experiments}

One can download {\tt LattE} via a web page
\url{www.math.ucdavis.edu/~latte}. You can also find there the files of all the
experiments presented in this section.  At the moment we
have been able to handle polytopes of dimension 30 and several
thousands vertices.  It is known that the theoretical upper bound of
the number of unimodular cones is $2^{dh}$, where $h = \lfloor
\frac{\log \log 1.9 - \log \log D }{\log(d-1/d)} \rfloor$ and where
$D$ is the volume of the fundamental parallelepiped of the input cone
\citep{bar}. If we fix the dimension this upper bound becomes
polynomial time. Unfortunately, if we do not fix the dimension, this
upper bound becomes exponential.  In practice this might be costly and
some families of polytopes have large numbers of unimodular cones.  The
cross polytope family, for instance, has many unimodular cones and behaves badly.  For example, for the cross polytope in $6$
dimensions, with cross6.ine input file \citep{fukuda}, {\tt LattE} took
$147.63$ seconds to finish computing.  The number of lattice points of
this polytope is obviously $13$. Also, for the cross polytope in $8$
dimensions, with cross8.ine input file \citep{fukuda}, {\tt LattE} took
85311.3 seconds to finish computing, even though this polytope has
only $16$ vertices and the number of lattice points of this polytope
is $17$.  For all computations, we used a $1$ GHz Pentium PC machine
running Red Hat Linux.

Here is a short description of how to use {\tt
LattE}. 

For computations involving a polytope $P$ described by a
system of inequalities $Ax\leq b$, where $A\in\Z^{m\times d}$, 
$A=(a_{ij})$, and $b\in\Z^m$, the {\tt LattE} readable input file
would be as follows: 
\begin{verbatim}
m d+1
b  -A
\end{verbatim}

\noindent
\textbf{EXAMPLE.}  
Let $P=\{(x,y): x\leq 1, y\leq 1, x+y\leq 1, x\geq 0, y\geq 0\}$.
Thus
\[
\begin{array}{ccc}
A=\left(
\begin{array}{rr} 
 1 &  0 \\ 
 0 &  1 \\ 
 1 &  1 \\
-1 &  0 \\ 
 0 & -1 \\ 
\end{array} 
\right) 
& , &
b = \left( 
\begin{array}{r} 
1 \\ 
1 \\ 
1 \\ 
0 \\
0 \\ 
\end{array} 
\right)
\end{array}
\]
and the {\tt LattE} input file would be as such:
\begin{verbatim}
5 3
1 -1  0
1  0 -1
1 -1 -1
0  1  0
0  0  1
\end{verbatim}

In {\tt LattE}, polytopes are represented by {\bf linear constraints},
i.e. equalities or inequalities. By default a constraint is an
inequality of type $ax\leq b$ unless we specify, by using a single
additional line, the line numbers of constraints that are linear
equalities. 

\noindent
\textbf{EXAMPLE.}  
Let $P$ be as in the previous example, but require $x+y=1$ instead of
$x+y\leq 1$, thus, 
$P=\{(x,y): x\leq 1, y\leq 1, x+y=1, x\geq 0, y\geq 0\}$.
Then the {\tt LattE} input file that describes $P$ would be as such:
\begin{verbatim}
5  3
1 -1  0
1  0 -1
1 -1 -1
0  1  0
0  0  1
linearity 1 3
\end{verbatim}
The last line states that among the $5$ inequalities one is to be
considered an equality, the third one.

For bigger examples it quickly becomes cumbersome to state all
nonnegativity constraints for the variables one by one. Instead, you
may use another short-hand.

\noindent
\textbf{EXAMPLE.}  
Let $P$ be as in the previous example, then the {\tt LattE} input file
that describes $P$ could also be described as such: 
\begin{verbatim}
3  3
1 -1  0
1  0 -1
1 -1 -1
linearity 1 3
nonnegative 2 1 2
\end{verbatim}
The last line states that there are two nonnegativity constraints and
that the first and second variables are required to be nonnegative. 
{\bf NOTE} that the first line reads ``3 3'' and not ``5 3'' as above!

We now report on computations with convex rational polytopes. We used
a $1$ GHz Pentium PC machine running Red Hat Linux. We begin with the
class of {\em multiway contingency tables}. A {\em $d$-table of size
$(n_1,\dots,n_d)$} is an array of non-negative integers
$v=(v_{i_1,\dots,i_d})$, $1\leq i_j\leq n_j$. For $0\leq m <d$, an
{\em $m$-marginal} of $v$ is any of the $d\choose m$ possible
$m$-tables obtained by summing the entries over all but $m$
indices. For instance, if $(v_{i,j,k})$ is a $3$-table then its
$0$-marginal is
$v_{+,+,+}=\sum_{i=1}^{n_1}\sum_{j=1}^{n_2}\sum_{k=1}^{n_3}
v_{i,j,k}$, its $1$-marginals are
$(v_{i,+,+})=(\sum_{j=1}^{n_2}\sum_{k=1}^{n_3}v_{i,j,k})$ and likewise
$(v_{+,j,+})$, $(v_{+,+,k})$, and its $2$-marginals are
$(v_{i,j,+})=(\sum_{k=1}^{n_3} v_{i,j,k})$ and likewise $(v_{i,+,k})$,
$(v_{+,j,k})$.

Such tables appear naturally in statistics and operations research
under various names such as {\em multi-way contingency tables}, or
{\em tabular data}.  We consider the {\em table counting problem}:
{\em given a prescribed collection of marginals, how many $d$-tables
are there that share these marginals?}  Table counting has several
applications in statistical analysis, in particular for independence
testing, and has been the focus of much research (see \citep{DG} and
the extensive list of references therein).  Given a specified
collection of marginals for $d$-tables of size $(n_1,\dots,n_d)$
(possibly together with specified lower and upper bounds on some of
the table entries) the associated {\em multi-index transportation
polytope} is the set of all non-negative {\em real valued} arrays
satisfying these marginals and entry bounds. The counting problem can
be formulated as that of counting the number of integer points in the
associated multi-index transportation polytope.  We begin with a small
example of a three-dimensional table of format $ 2 \times 3 \times 3$
given below. The data displayed in Table \ref{census} have been
extracted from the 1990 decennial census and is used in
\citep{FMMS}. For the $2$-marginals implied by these data we get the
answer of 441 in less than a second.

\begin{table}[ht]
\begin{center}
\begin{tabular}{|l|c|c|c||c|}
\multicolumn{5}{c}{Gender $=$ Male} \\ 
\multicolumn{1}{c}{} & \multicolumn{4}{c}{Income Level} \\ \cline{2-5}
\multicolumn{1}{c|}{Race} & $\le \$10,000$ & $> \$10000 \mbox{ and} \le
\$25000$ & $> \$25000 $ & Total \\ \hline
White & 96 & 72 & 161 & 329 \\ \hline
Black & 10 & 7 & 6 & 23 \\ \hline
Chinese & 1 & 1 & 2 & 4 \\ \hline\hline
Total & 107 & 80 & 169 & 356 \\ \hline
\multicolumn{5}{c}{} \\[11pt] 
\multicolumn{5}{c}{Gender $=$ Female} \\ 
\multicolumn{1}{c}{} & \multicolumn{4}{c}{Income Level} \\ \cline{2-5}
\multicolumn{1}{c|}{Race} & $\le \$10,000$ & $> \$10000 \mbox{ and} \le
\$25000$ & $> \$25000 $ & Total \\ \hline
White & 186 & 127 & 51 & 364 \\ \hline
Black & 11 & 7 & 3 & 21 \\ \hline
Chinese & 0 & 1 & 0 & 1 \\ \hline\hline
Total & 197 & 135 & 54 & 386 \\ \hline
\end{tabular}
\end{center}
\caption{Three-way cross-classification of gender, race, and income for a
selected U.S. census tract. {\em Source}: 1990 Census Public Use Microdata
Files.}\label{census}
\end{table}

 We present now an example of a $3 \times 3 \times 3$ table with
fairly large $2$-marginals. They are displayed in Table \ref{3by3by3}.
{\tt LattE} took only $19.67$ seconds of CPU time.  The number of
lattice points inside of this polytope is
$ 2249847900174017152559270967589010977293 $.

\begin{table}
\begin{center}
\begin{tabular}{|c|c|c|} \hline
      164424  &    324745  &    127239  \\ \hline
      262784  &    601074  &   9369116  \\ \hline
      149654  &   7618489  &   1736281  \\ \hline
\end{tabular}
\end{center}

\begin{center}
\begin{tabular}{|c|c|c|} \hline
      163445   &    49395 &     403568  \\ \hline
     1151824   &   767866 &    8313284  \\ \hline
     1609500   &  6331023 &    1563901  \\ \hline

\end{tabular}
\end{center}

\begin{center}
\begin{tabular}{|c|c|c|} \hline
      184032  &    123585 &     269245  \\ \hline
      886393  &   6722333 &     935582  \\ \hline
     1854344  &    302366 &    9075926  \\ \hline

\end{tabular}
\end{center}
\caption{2-Marginals for the $3 \times 3 \times 3$ example.} \label{3by3by3}
\end{table}

 Next we present an example of a $3 \times 3 \times 4$ table with
large 2-marginals. The $2$-marginals are displayed in Table
\ref{3by3by4}.  The CPU time for this example was $44$ minutes $42.22$
seconds. The number of lattice points inside of this polytope is

$ 4091700129572445106288079361219676736812805058988286839062994.$

\begin{table}
\begin{center}
\begin{tabular}{|c|c|c|c|} \hline
      273510  &    273510  &    273510  &    191457  \\ \hline
      273510  &    273510  &    547020  &    191457  \\ \hline
      273510  &    547020  &    273510  &    191457  \\ \hline
\end{tabular}
\end{center}

\begin{center}
\begin{tabular}{|c|c|c|} \hline

      464967  &    273510  &    273510   \\ \hline 
      547020  &    273510  &    464967  \\ \hline
      410265  &    601722  &    273510  \\ \hline

\end{tabular}
\end{center}

\begin{center}
\begin{tabular}{|c|c|c|} \hline

      273510  &    273510  &    273510   \\ \hline
      410265  &    547020  &    136755   \\ \hline
      547020  &    136755  &    410265   \\ \hline
      191457  &    191457  &    191457   \\ \hline

\end{tabular}
\end{center}
\caption{2-Marginals for the $3 \times 3 \times 4$ example.} \label{3by3by4}
\end{table}

The next family of examples are some hard knapsack-type
problems. Suppose we have a set of positive relatively prime integers
$\{a_1, a_2, \ldots , a_d \}$. Denote by $a$ the vector $(a_1, a_2,
\ldots , a_d)$.  Consider the following problem: does there exist a
non-negative integral vector $x$ satisfying $ax = a_0$ for some
positive integer $a_0$?  We take several examples from
\citep{aardaletal3} which have been found to be extremely hard to solve
by commercial quality branch-and-bound software. This is very
surprising since the number of variables is at most 10. It is not very
difficult to see that if the right-hand-side value $a_0$ is large
enough, the equation will surely have a non-negative integer solution.
The {\em Frobenius number} for a knapsack problem is the largest value
$a_0$ such that the knapsack problem is infeasible. Aardal and Lenstra
\citep{aardaletal3} solved them using the reformulation in
\citep{aardaletal1}.  Their method works significantly better than
branch-and-bound using {\tt CPLEX} 6.5. Here we demonstrate that our
implementation of Barvinok's algorithm is fairly fast and, on the
order of seconds, we resolved the first 15 problems in Table 1 of
\citep{aardaletal3} and verified all are infeasible except {\bf prob9},
where there is a mistake. The vector $(3480,1,4,4,1,0,0,0,0,0)$ is a
solution to the right-hand-side $13385099.$ In fact, using {\tt LattE}
we know that the exact number of solutions is $838908602000$.  For
comparison we named the problems exactly as in Table 1 of
\citep{aardaletal3}. We present our results in Table \ref{knapsack}.
It is very interesting to know the
number of lattice points if we add $1$ to the Frobenius number for
each problem.  In Table \ref{knap2}, we find the number of solutions
if we add $1$ to the Frobenius number on each of the (infeasible)
problems. The speed is practically the same as in the previous
case. In fact the speed is the same regardless of the right-hand-side
value $a_0$.

\begin{sidewaystable}
\centering
\begin{tabular}{|c|cccccccccc|c|c|} \hline

&&&&&&&&& &&Frobenius \# & Time (m, s)\\ \hline
cuww1&12223 &12224 &36674 &61119 &85569& & & & & & 89643481 & 0.55s\\ 
\hline
cuww2&12228 &36679 &36682 &48908 &61139 &73365 & & & & &89716838& 1.78s  
\\ \hline
cuww3&12137 &24269 &36405 &36407 &48545 &60683 & & & & &58925134 & 1.27s\\ 
\hline
cuww4&13211 &13212 &39638 &52844 &66060 &79268 &92482 & & & &104723595& 
2.042s \\ \hline
cuww5&13429 &26850 &26855 &40280 &40281 &53711 &53714 &67141 & & & 
45094583& 16.05s\\ \hline
pro1& 25067 &49300 &49717 &62124 &87608 &88025 &113673 &119169 & & & 
33367335& 47.07s  \\ \hline 
prob2& 11948 &23330 &30635 &44197 &92754 &123389 &136951 &140745 & & 
&14215206 & 1m0.58s \\ \hline
prob3  &39559 &61679 &79625 &99658 &133404 &137071 &159757 &173977 & & 
&58424799 & 1m28.3s \\ \hline
prob4& 48709 &55893 &62177 &65919 &86271 &87692 &102881 &109765 & & 
&60575665 & 59.04s \\ \hline
prob5& 28637 &48198 &80330 &91980 &102221 &135518 &165564 &176049 & & & 
62442884& 1m41.78s \\ \hline
prob6 &20601 &40429 &40429 &45415 &53725 &61919 &64470 &69340 &78539 
&95043 &
22382774 & 3m45.86s  \\ \hline
prob7 &18902 &26720 &34538 &34868 &49201 &49531 &65167 &66800 &84069 
&137179 &
27267751 &  2m57.64s\\  \hline
prob8& 17035 &45529 &48317 &48506 &86120 &100178 &112464 &115819 &125128 
&129688
& 21733990 & 8m29.78s\\ \hline
prob10 & 45276 &70778 &86911 &92634 &97839 &125941 &134269 &141033 &147279 
&153525 &106925261 & 4m24.67s \\ \hline
\end{tabular}
\centering
\caption{Infeasible knapsack problems.} \label{knapsack}
\end{sidewaystable}

\begin{table}
\begin{center}
\begin{tabular}{|c|c|c|}\hline
problem & RHS & \# of lattice points. \\ \hline
cuww1& 89643482 & 1 \\ \hline
cuww2 & 89716839 & 1 \\ \hline
cuww3 & 58925135 & 2 \\ \hline
cuww4 & 104723596 & 1 \\ \hline
cuww5 & 45094584  & 1 \\ \hline
prob1 & 33367336 & 859202692 \\ \hline
prob2 & 14215207  & 2047107 \\ \hline
prob3 & 58424800   &  35534465752\\ \hline
pro4 & 60575666   & 63192351 \\ \hline
pro5 & 62442885   & 21789552314 \\ \hline
pro6 & 22382775   &  218842 \\ \hline
pro7 & 27267752   & 4198350819898 \\ \hline
pro8 & 21733991   & 6743959 \\ \hline
pro10 & 106925262  & 102401413506276371 \\ \hline
\end{tabular}
\caption{The number of lattice points if we add 1 to the Frobenius
number.}\label{knap2}
\end{center}
\end{table}

Already counting the lattice points of large width convex polygons is
a non-trivial task if one uses brute-force enumeration (e.g. list one
by one the points in a bounding box of the polygon and see whether it
is inside the polygon). Here we experiment with very large convex {\em
almost} regular $n$-gons. Regular $n$-gons cannot have rational
coordinates, but we can approximate them to any desired accuracy by
rational polygons. In the following experiment we take regular
$n$-gons, from $n=5$ to $n=12$ centered at the origin (these have only
a handful of lattice points). We take a truncation of the coordinates
up to 3, 9, and 15 decimal digits, then we multiply by a large enough
power of 10 to make those vertex coordinates integral and we count the
number of lattice points in the dilation. All experiments take less
than a second.

{\tiny
\begin{table} 
\begin{center}
{\small
\begin{tabular}{|c|c|c|c|} \hline
 & $10^3$ (seconds) & $10^9$ (seconds)& $10^{15}$ (seconds) \\ \hline
5gon & 2371673(0.136) & 2377641287748905186(0.191) & 2377641290737895844565559026875(0.289) \\ \hline
6gon & 2596011(0.153) & 2598076216000000011(0.193) & 2598076211353321000000000000081(0.267) \\ \hline
7gon & 2737110(0.175) & 2736410188781217941(0.318) & 2736410188638105174143840143912(0.584s) \\ \hline
8gon & 2820021(0.202) & 2828427120000000081(0.331) & 2828427124746200000000000000201(0.761) \\ \hline
9gon & 2892811(0.212) &  2892544245156317460(0.461) & 2892544243589428566861745742966(0.813) \\ \hline
10gon & 2931453(0.221) &  2938926257659276211(0.380) & 2938926261462380264188126524437(0.702) \\ \hline
11gon & 2974213(0.236) & 2973524496796366173(0.745) &  2973524496005786351949189500315(1.858) \\ \hline
12gon & 2997201(0.255) & 3000000004942878881(0.466) & 3000000000000005419798779796241(0.696) \\ \hline
\end{tabular}
\caption{The numbers of the approximated regular polygons. We show
the number of lattice points in different dilation factors (powers
of ten) and time of computation.}
}
\end{center}
\end{table}
}

The next two sets of examples are families that have been studied quite
extensively in the literature and provide us with a test for speed.
In the first case we deal with {\em two-way contingency tables.} The
polytope defined by a two-way contingency table is called the {\em
transportation polytope}.  We present the results in Table
\ref{transport}.  The second family consists of flow polytopes for the
complete $4$-vertex and the complete $5$-vertex tournaments (directed complete  graphs).  Consider the
directed complete graph $K_l$ for $l \in \N$ and $l \geq 3$.  We
assign a number to each node of the graph.  Then, we orient the arcs
 from the node of smaller index to the node of bigger index. 
Let $N$ be the node set of the
complete graph $K_l$, let $w_i$ be a weight assigned to node $i$ for $i =
1, 2, \ldots , l$, and let $A$ be the arc set of $K_l$.  Then, we have
the following constraints, with as many variables as arcs:

$$\sum_{(j,i) \mbox{arc enters} i} x_{ji} - \sum_{(i,j) \mbox{arc has tail} i} x_{ij} = w_i,$$
$$ x_{ij} \geq 0 \mbox{ } \forall  (i,j) \in
A.$$ These equalities and inequalities define a polytope and this
polytope is the special case of a {\em flow polytope}. The results for the complete
graphs $K_4$ and $K_5$, with different weight vectors, are shown in Tables 
\ref{k4} and \ref{k5} respectively.

These two families of polytopes have been studied by several authors\\
\citep{baldonideloeravergne, Beck2, deloerasturmfels, mount} and thus
are good for testing the performance of {\tt LattE}. We used several
examples of transportation polytopes, as presented in the table
below. In general, {\tt LattE} runs at comparable performance to the
software of \citep{baldonideloeravergne, Beck2} for generic vectors
$(a,b)$ but is slower for degenerate inputs (those that do not give a
simple polytope). The reason seems to be that at each non-simplex vertex {\tt
LattE} needs to triangulate each cone which takes considerable time in
problems of high dimension.

\begin{sidewaystable}
\centering
{\tiny
\begin{tabular}{|p{2.7in}|r|c|} \hline 
Margins  & \# of lattice points & Time (seconds)\\ \hline 
[220, 215, 93, 64], \newline 
[108, 286, 71, 127] &  1225914276768514 & 1.048 \\ \hline
[109, 127, 69, 109], \newline [119, 86, 108, 101] & 993810896945891 & 1.785 \\ \hline
[72, 67, 47, 96], \newline [70, 70, 51, 91] &  25387360604030 &  1.648\\ \hline
[179909, 258827, 224919, 61909], \newline [190019, 90636, 276208, 168701] &  13571026063401838164668296635065899923152079 &  1.954  \\ \hline
[229623, 259723, 132135, 310952],\newline  [279858, 170568, 297181, 184826] &646911395459296645200004000804003243371154862 &  1.765 \\ \hline
[249961, 232006, 150459, 200438],\newline  [222515, 130701, 278288, 201360] & 319720249690111437887229255487847845310463475 &  1.854  \\ \hline 
[140648, 296472, 130724, 309173], \newline [240223, 223149, 218763, 194882]&322773560821008856417270275950599107061263625 &  1.903 \\ \hline
[65205, 189726, 233525, 170004],\newline  [137007, 87762, 274082, 159609]  & 6977523720740024241056075121611021139576919 &1.541 \\ \hline
[251746, 282451, 184389, 194442], \newline [146933, 239421, 267665, 259009] & 861316343280649049593236132155039190682027614 & 1.880 \\ \hline
[138498, 166344, 187928, 186942], \newline [228834, 138788, 189477, 122613] & 63313191414342827754566531364533378588986467 & 1.973 \\ \hline
[20812723, 17301709, 21133745, 27679151], \newline [28343568, 18410455,
19751834, 20421471] & 
665711555567792389878908993624629379187969880179721169068827951
& 2.917\\ 
\hline

[15663004, 19519372, 14722354, 22325971], \newline [17617837, 25267522,
20146447, 9198895] &
63292704423941655080293971395348848807454253204720526472462015 &
3.161
\\
\hline

[13070380, 18156451, 13365203, 20567424], \newline [12268303, 20733257,
17743591, 14414307] &
43075357146173570492117291685601604830544643769252831337342557 &
2.990
\\
\hline

\end{tabular}
\caption{Testing for $4 \times 4$ transportation polytopes.}
}
\label{transport}
\end{sidewaystable} 

\begin{table}
\begin{center}
\begin{tabular}{|l|r|c|} \hline 
Weights on nodes  & \# of lattice points & Time (seconds)\\ \hline 
[-6, -8, 5, 9] &  223 & 0.288 \\ \hline
[-9, -11, 12, 8] & 330 & 0.286 \\ \hline
[-1000, -1, 1000, 1] &  3002 & 0.287 \\ \hline
[-4383, 886, 2777, 720] &  785528058 &  0.287  \\ \hline
[-4907, -2218, 3812, 3313] & 20673947895 &  0.288 \\ \hline
[-2569, -3820, 1108, 5281] & 14100406254 &  0.282  \\ \hline 
[-3842, -3945, 6744, 1043]& 1906669380  &  0.281 \\ \hline
[-47896, -30744, 46242, 32398]  &  19470466783680  & 0.282 \\ \hline
[-54915, -97874, 64165, 88624] & 106036300535520 & 0.281 \\ \hline
[-69295, -62008, 28678, 102625] & 179777378508547 & 0.282 \\ \hline
[-3125352, -6257694, 926385, 8456661] & 34441480172695101274 & 0.509 \\
\hline

[-2738090, -6701290, 190120, 9249260] & 28493245103068590026 & 0.463 \\
\hline

[-6860556, -1727289, 934435, 7653410] & 91608082255943644656 & 0.503 \\
\hline

\end{tabular}
\caption{Testing for the complete graph $K_4$.} \label{k4}
\end{center}
\end{table}

{\tiny
\begin{table}
{\small
\begin{tabular}{|l|r|c|} \hline 
Weights on nodes  & \# of lattice points & secs\\ \hline 
[-12, -8, 9, 7, 4]&  14805 & 0.319 \\ \hline
[-125, -50, 75, 33, 67] & 6950747024 & 0.325 \\ \hline
[-763, -41, 227, 89, 488] &  222850218035543 & 0.325 \\ \hline
[-11675, -88765, 25610, 64072, 10758] &  563408416219655157542748  &  0.319 \\ \hline
[-78301, -24083, 22274, 19326, 60784] & 1108629405144880240444547243 &  0.336 \\ \hline
[-52541, -88985, 1112, 55665, 84749] & 3997121684242603301444265332 &  0.331  \\ \hline 
[-71799, -80011, 86060, 39543, 26207]& 160949617742851302259767600 & 0.316  \\ \hline
[-45617, -46855, 24133, 54922, 13417]  & 15711217216898158096466094 & 0.285 \\ \hline
[-54915, -97874, 64165, 86807, 1817] & 102815492358112722152328 & 0.277 \\ \hline
[-69295, -62008, 28678, 88725, 13900] & 65348330279808617817420057 & 0.288 \\ \hline
[-8959393, -2901013, 85873, 533630, 11240903] &
6817997013081449330251623043931489475270 & 0.555 \\ \hline
[-2738090, -6701290, 190120, 347397, 8901863] &
277145720781272784955528774814729345461 & 0.599 \\ \hline
[-6860556, -1727289, 934435, 818368, 6835042] &
710305971948234346520365668331191134724 & 0.478 \\ \hline

\end{tabular}
\centering
\caption{Testing for the complete graph $K_5$. Time is given in seconds} \label{k5}
}
\end{table}
}

The following experiment is from \cite{mcmc2}.  
This is from an actual data in German survey.  For 2,262 German citizens they asked the following questions: if you order the following items, how you order them? 
\begin{enumerate}
\item Maintain order.
\item Give people more say.
\item Fight rising prices.
\item Protect freedom of speech.
\end{enumerate}

Then we have the data in Table \ref{S4:data}.

\begin{table}
\begin{center}
\begin{tabular}{||c|c||c|c||c|c||c|c||}\hline
1234 & 137 & 2134 & 48 & 3124 & 330 & 4123 & 21 \\ \hline
1243 & 29  & 2143 & 23 & 3142 & 294 & 4132 & 30 \\ \hline
1324 & 309 & 2314 & 61 & 3214 & 117 & 4213 & 29 \\ \hline
1342 & 255 & 2341 & 55 & 3241 & 69  & 4231 & 52 \\ \hline
1423 & 52  & 2413 & 33 & 3412 & 70  & 4312 & 35 \\ \hline
1432 & 93  & 2431 & 39 & 3421 & 34  & 4321 & 27 \\ \hline
Total & 875 & & 279 && 914 && 194 \\ \hline
\end{tabular}
\caption{Permutation $S_4$ problem from \cite{mcmc2}.}\label{S4:data}
\end{center}
\end{table}

We take sum for the number of people who picked $j \in \{1, 2, 3, 4\}$ for the $i$th order for all $j = 1, 2, 3, 4$ and for all $i = 1, 2, 3, 4$.  Then we have the following condition given in Table \ref{S4}.  

\begin{table}
\begin{center}
\begin{tabular}{|c|c|c|c|}\hline
875 & 279 & 914 & 194 \\ \hline
746 & 433 & 742 & 341 \\ \hline
345 & 773 & 419 & 725 \\ \hline
296 & 777 & 187 & 1002 \\ \hline
\end{tabular}
\caption{Marginal conditions for the permutation problem from \cite{mcmc2}.}\label{S4}
\end{center}
\end{table}

Then we want to compute how many functions $f$ such that $f: S_4 \to \N$ and $\sum_{\sigma \in S_4} f(\sigma) = 2262$.

The solution to this problem is: Total number of functions is\\ $11606690287805167142987310121$ and 
CPU Time is 523.12 sec.

The last experiment in this section
is from Ian Dinwoodie.
Ian Dinwoodie 
communicated to us the problem of counting all $7 \times 7$
contingency tables whose entries are nonnegative integers $x_i$, with
diagonal entries multiplied by a constant as presented in Table
\ref{dinwoodie}. The row sums and column sums of the entries are given
there too. Using {\tt LattE} we obtained the exact answer {\em
{8813835312287964978894}}.

\begin{table}
\begin{center}
\begin{tabular}{|c|c|c|c|c|c|c||c|} \hline
  2{\color{blue}{$x_1$}}  &   $x_2$   &   $x_3$ &  $x_4$ &  $x_5$ & $x_6$ & $x_7$ & 205\\ \hline
     $x_2$   &  2{\color{blue}{$x_8$}}   &   $x_{9}$ &  $x_{10}$ &  $x_{11}$ & $x_{12}$ & $x_{13}$ & 600 \\ \hline
     $x_3$  & $x_9$   &  2{\color{blue}{$x_{14}$}}  &  $x_{15}$ &  $x_{16}$ & $x_{17}$ & $x_{18}$ & 61\\ \hline
     $x_4$   & $x_{10}$ & $x_{15}$      &  2{\color{blue}{$x_{19}$}} &  $x_{20}$ & $x_{21}$ & $x_{22}$ & 17 \\ \hline
     $x_5$  & $x_{11}$ & $x_{16}$     &$x_{20}$  & 2{\color{blue}{$x_{23}$}}  & $x_{24}$ & $x_{25}$ & 11 \\ \hline
     $x_6$  &$x_{12}$  & $x_{17}$    & $x_{21}$  & $x_{24}$  & 2{\color{blue}{$x_{26}$}} & $x_{27}$ & 152 \\ \hline
     $x_7$ &  $x_{13}$ & $x_{18}$    & $x_{22}$  & $x_{25}$  & $x_{27}$  & 2{\color{blue}{$x_{28}$}} & 36 \\ \hline
    205 & 600 & 61 & 17 & 11 & 152 &36 & 1082 \\ \hline
\end{tabular}
\caption{The conditions for retinoblastoma RB1-VNTR genotype data from
 the Ceph database.} \label{dinwoodie} 
\end{center}
\end{table}

\section{New Ehrhart (quasi-)polynomials}\label{formulas}
 
Given a rational polytope $P \subset \R^d$, the function \[ i_P(t) :=
  \# \left( t P \cap \Z^d \right),\] for a positive integer $t$, was
  first studied by E. Ehrhart \citep{Ehrhart} and has received a lot of
  attention in combinatorics.  It is known to be a polynomial when all
  vertices of $P$ are integral and it is a quasi-polynomial for
  arbitrary rational polytopes. It is called the {\em Ehrhart
  quasi-polynomial} in honor of its discoverer 
  \citep[Chapter 4]{stanley}.  A function $f: \N \rightarrow \C$ is a
  quasi-polynomial if there is an integer $N > 0$ and polynomials
  $f_0, \ldots , f_{N-1}$ such that $f(s) = f_i(s) \mbox{ if } s
  \equiv i \mod N$.  The integer $N$ is called a {\em quasi-period} of
  $f$.  Therefore, by counting the number of lattice points for
  sufficiently many dilations of a rational polytope, we can
  interpolate its Ehrhart quasi-polynomial.

Using {\tt LattE}, {\tt Maple}, and interpolation, we have calculated the Ehrhart polynomials and
  quasi-polynomials for polytopes that are slices or nice truncations
  of the unit $d$-cube.  To the best of our knowledge these values
  were not known before. For example, the $24$-cell polytope centered
  at the origin with smallest integer coordinates has Ehrhart
  polynomial $i_{24\_cell}(s) = 8s^4 + \frac{32s^3}{3} + 8s^2 +
  \frac{16s}{3} + 1.$ In Table \ref{hypersi}, we see the Ehrhart
  polynomials for the hypersimplices $\Delta(n,k)$. They are defined
  as the slice of the $n$-cube by the hyperplane of equation $\sum x_i=k$ 
  with $k \leq n$. Note that $\Delta(n,k)=\Delta(n,n-k)$ because of the
symmetries of the regular cube. The hypersimplices form one of the most famous families of $0/1$-polytopes.
It is known that hypersimplices are {\em compressed polytopes} \citep{hibi}.
This means that their Ehrhart polynomials can be recovered from the $f$-vectors
of any of their reverse lexicographic triangulations. Instead, we recovered 
them explicitly for the first time  using {\tt LattE} and interpolation.


\begin{sidewaystable}
\centering
{\tiny
\begin{tabular}{|c|c|c|} \hline
n & k  & the Ehrhart polynomial $P(s)$ \\ \hline 
4 & 1 & $\frac{s^3}{6} + s^2 + \frac{11s}{6} + 1$  \\ \hline
4 & 2 & $\frac{2x^3}{3} + 2s^2 + \frac{7s}{3} + 1$ \\ \hline
5 & 1 & $\frac{s^4}{24} + \frac{5s^3}{12} + \frac{35s^2}{24} +\frac{25s}{12} + 1 $\\ \hline
5 & 2 & $\frac{11s^4}{24} + \frac{25s^3}{12} + \frac{85s^2}{24} + 
\frac{35s}{12} + 1$ \\ \hline
6 & 1 & $\frac{s^5}{120} + \frac{s^4}{8} + \frac{17s^3}{24} + \frac{15 
s^2}{8} + \frac{137s}{60} + 1$  \\ \hline
6 & 2 & $\frac{13s^5}{60} + \frac{3s^4}{2} + \frac{47s^3}{12} + 5s^2 + \frac{101 s}{30} + 1$ \\ \hline
6 & 3 & $\frac{11s^5}{20}+\frac{11s^4}{4} + \frac{23s^3}{4} + \frac{25s^2}{4} +\frac{37s}{10} + 1$   \\ \hline
7 & 1 & $\frac{s^6}{720} + \frac{7s^5}{240} + \frac{35s^4}{144} + \frac{49s^3}{48} + \frac{203s^2}{90} + \frac{49s}{20} + 1$ \\ \hline 
7 & 2  & $\frac{19s^6}{240}+ \frac{63s^5}{80} + \frac{49s^4}{16} + \frac{287s^3}{48} + \frac{763s^2}{120} + \frac{56s}{15} + 1$ \\ \hline
7 & 3 & $\frac{151s^6}{360} + \frac{161s^5}{60} + \frac{256s^4}{36} + \frac{21s^3}{2} + \frac{3199s^2}{360} + \frac{259s}{60} + 1$ \\ \hline
8 & 1 &
$\frac{s^7}{5040}+\frac{s^6}{180}+\frac{23s^5}{360}+\frac{7s^4}{18} +
\frac{967s^3}{720} + \frac{469s^2}{180}+\frac{363s}{140} + 1$ \\ \hline
8 & 2 &
$\frac{s^7}{42}+\frac{29s^6}{90}+\frac{53s^5}{30}+\frac{91s^4}{18}+
\frac{49s^3}{6}+\frac{343s^2}{45}+\frac{283s}{70}+1$ \\ \hline
8 & 3 &
$\frac{397s^7}{1680}+\frac{359s^6}{180}+\frac{281s^5}{40}+\frac{245s^4}{18}
+ \frac{1273s^3}{80}+\frac{2051s^2}{180}+\frac{2027s}{420}+1$ \\ \hline
8 & 4 & $\frac{151s^7}{315}+\frac{151s^6}{45}+\frac{463s^5}{45} +
\frac{161s^4}{9}+\frac{862s^3}{45}+\frac{574s^2}{45}+\frac{533s}{105}+1$ 
\\ \hline
9 & 1 &$\frac{s^8}{40320}+\frac{s^7}{1120}+\frac{13s^6}{960}+\frac{9s^5}{80}+\frac{1069s^4}{1920}+\frac{267s^3}{160}+ \frac{29531s^2}{10080}+\frac{761s}{280}+1$ \\ \hline
9&2& $\frac{247s^8}{40320}+\frac{121s^7}{1120} + \frac{763s^6}{960}+\frac{253s^5}{80}+\frac{14203s^4}{1920}+\frac{1667s^3}{160}+\frac{88721s^2}{10080}+\frac{1207s}{280}+1$ \\ \hline
9 & 3 & $\frac{477s^8}{4480}+\frac{1311s^7}{1120}+\frac{1731s^6}{320}+\frac{1107s^5}{80}+\frac{13899s^4}{640}+\frac{3477s^3}{160}+\frac{15419s^2}{1120}+\frac{1473s}{280}+1$ \\ \hline
9 & 4 & $\frac{15619s^8}{40320}+\frac{3607s^7}{1120}+\frac{11311s^6}{960}+\frac{1991s^5}{80}+\frac{63991s^4}{1920}+\frac{4669s^3}{160}+\frac{166337s^2}{10080}+\frac{1599s}{280}+1$ \\ \hline

10 & 1 & $\frac{s^9}{362880}+\frac{s^8}{8064}+\frac{29s^7}{12096}+\frac{5s^6}{192}+\frac{3013s^5}{17280}+\frac{95s^4}{128}+\frac{4523s^3}{2268}+\frac{6515s^2}{2016}+\frac{7129s}{2520}+1$ \\ \hline

10 & 2 & $\frac{251s^9}{181440}+\frac{31s^8}{1008}+\frac{1765s^7}{6048}+\frac{37s^6}{24}+\frac{42863s^5}{8640}+\frac{481s^4}{48}+\frac{115205s^3}{9072}+\frac{4993s^2}{504}+\frac{5729s}{1260}+1$ \\ \hline

10 & 3 & $\frac{913s^9}{22680}+\frac{1135s^8}{2016}+\frac{5071s^7}{1512}+\frac{179s^6}{16}+\frac{3128s^5}{135} +\frac{2999s^4}{96}+\frac{63041s^3}{2268}+\frac{ 8069s^2}{504} +\frac{3553s}{630}+1$ \\ \hline

10 & 4 & $\frac{44117s^9}{181440}+\frac{2489s^8}{1008}+\frac{66547s^7}{6048}+\frac{683s^6}{24}+\frac{409361s^5}{ 8640}+\frac{2543s^4}{48}+\frac{363947s^3}{9072}+\frac{10127s^2}{504}+\frac{7883s}{1260}+1$ \\ \hline

10 & 5 & $ \frac{15619s^9}{36288}+\frac{15619s^8}{4032}+\frac{94939s^7}{6048}+\frac{3607s^6}{96}+\frac{101311s^5}{1728}+\frac{11911s^4}{192}+\frac{25394s^3}{567}+\frac{21689s^2}{1008}+\frac{1627s}{252}+1$ \\ \hline

11 & 1 & $\frac{s^{10}}{3628800}+\frac{11s^9}{725760}+\frac{11s^8}{30240}+\frac{121s^7}{24192}+\frac{7513s^6}{172800}+\frac{8591s^5}{34560}+\frac{341693s^4}{362880}+\frac{84095s^3}{36288}+\frac{177133s^2}{50400}+\frac{7381s}{2520}+1$ \\ \hline

11 & 2 & $\frac{1013s^{10}}{3628800}+\frac{ 5533s^9}{725760}+\frac{2189s^8}{24192}+\frac{14795s^7}{24192}+\frac{447689s^6}{172800}+\frac{246697s^5}{34560}+\frac{14597s^4}{1134}+\frac{543763s^3}{36288}+\frac{91949s^2}{8400}+\frac{1199s}{252}+1$ \\ \hline

11 & 3 & $\frac{299s^{10}}{22680}+\frac{16621s^9}{72576}+\frac{41591s^8}{24192}+\frac{88693s^7}{12096}+\frac{170137s^6}{8640}+\frac{604109s^5}{17280}+\frac{3043997s^4}{72576}+\frac{308473s^3}{9072}+\frac{60929s^2}{3360}+\frac{15059s}{2520}+1$ \\ \hline

11 & 4 & $\frac{56899s^{10}}{453600}+\frac{565631s^9}{362880}+\frac{205733s^8}{24192}+\frac{326491s^7}{12096}+\frac{2400629s^6}{43200}+\frac{1348787s^5}{17280}+\frac{5535695s^4}{72576}+\frac{468655s^3}{9072}+\frac{1185701s^2}{ 50400}+\frac{16973s}{2520}+1$ \\ \hline

11 & 5 & $\frac{655177s^{10}}{1814400}+\frac{336083s^9}{90720}+\frac{2078791s^8}{120960}+\frac{287639s^7}{6048}+\frac{7525771s^6}{86400}+\frac{95557s^5}{864}+\frac{35914087s^4}{362880}+\frac{1125575s^3}{18144}+\frac{443179s^2}{16800}+\frac{17897s}{2520}+1$ \\ \hline

12 & 1 & $\frac{s^{11}}{39916800}+\frac{s^{10}}{604800}+\frac{s^9}{20736}+\frac{11s^8}{13440}+\frac{10831s^7}{1209600}+\frac{1903s^6}{28800}+\frac{242537s^5}{725760}+\frac{139381s^4}{120960}+\frac{341747s^3}{129600}+\frac{190553s^2}{50400 }+\frac{83711s}{ 27720}+1 $ \\ \hline

12 & 2 & $\frac{509s^{11}}{9979200}+\frac{169s^{10}}{100800}+\frac{551s^9}{22680}+\frac{2057s^8}{10080}+\frac{332249s^7}{302400}+\frac{18997s^6}{4800}+\frac{876959s^5}{90720}+\frac{80179s^4}{5040 }+\frac{244681s^3}{14175}+\frac{150293s^2}{12600}+\frac{68591s}{13860}+1$ \\ \hline

12 & 3 & $\frac{50879s^{11}}{13305600}+\frac{6979s^{10}}{86400}+\frac{60271s^9}{80640}+\frac{32153s^8}{8064}+\frac{5483809s^7}{403200}+\frac{897259s^6}{28800}+\frac{11875111s^5}{241920}+\frac{185339s^4}{3456}+\frac{451173s^3}{11200}+\frac{338503s^2}{16800}+\frac{58007s}{9240}+1$ \\ \hline

12 & 4 & $\frac{1093s^{11}}{19800}+\frac{62879s^{10}}{75600}+\frac{20893s^9}{3780}+\frac{10813s^8}{504}+\frac{684323s^7}{12600}+\frac{340967s^6}{3600}+\frac{5258s^5}{45}+\frac{38819s^4}{378}+\frac{1202029s^3}{18900}+\frac{42218s^2}{1575}+\frac{1103s}{154}+1$ \\ \hline

12 & 5 & $\frac{1623019s^{11}}{6652800}+\frac{882773s^{10}}{302400}+\frac{1908073s^9}{120960}+\frac{1028401s^8}{20160}+\frac{7395023s^7}{67200}+\frac{2401619s^6}{14400}+\frac{4398559s^5}{24192}+\frac{8661917s^4}{60480}+\frac{12163441s^3}{151200}+\frac{782969s^2}{25200}+\frac{8861s}{1155}+1 $ \\ \hline

12 & 6 & $\frac{655177s^{11}}{1663200}+\frac{655177s^{10}}{151200}+\frac{5507s^9}{252}+\frac{336083s^8}{ 5040}+\frac{6898277s^7}{50400}+\frac{1430341s^6}{7200}+\frac{3152491s^5}{15120}+\frac{1200463s^4}{7560}+\frac{30291s^3}{350}+\frac{68321s^2}{2100}+\frac{18107s}{2310}+1$ \\ \hline

\end{tabular}
\caption{The Ehrhart polynomials for the hypersimplices $\Delta(n, k)$ } 
\label{hypersi}
}
\end{sidewaystable}

We also have the Ehrhart quasi-polynomials of some truncated unit cubes.

\begin{proposition}
The Ehrhart quasi-polynomial for the truncated unit cube in Figure 
\ref{tru1}, where its vertices are at $1/3$ and $2/3$ of the way along
edges of the cube, is given by:

$$i_{tru\_cube1}(s) = \left \{ \begin{array}{ll} \frac{77s^3}{81} + 
\frac{23s^2}{9} + \frac{19s}{9} + 1 \mbox{ if } s \equiv 0 \mod 3, \\ 
 \frac{77s^3}{81} + \frac{61s^2}{27} - \frac{7s}{27}
- \frac{239}{81} \mbox{ if } s \equiv 1 \mod 3, \\
 \frac{77s^3}{81} + \frac{65s^2}{27} +
\frac{29s}{27} - \frac{31}{81} \mbox{ if } s \equiv 2 \mod 3.\\
\end{array}\right .$$


\begin{figure}[ht]
\begin{center}
\includegraphics[width=5 cm]{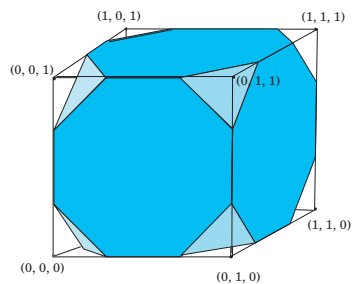}  
\caption{The truncated cube.} \label{tru1}
\end{center} 
\end{figure}
\end{proposition}

\begin{proposition}
The Ehrhart quasi-polynomial for the cuboctahedron (Figure \ref{tru2}) is:

$$i_{tru\_cube2}(s) = \left \{ \begin{array}{ll} \frac{5s^3}{6} + 2s^2 + 
\frac{5s}{3} + 1 \mbox{ if } s\equiv 0 \mod 2, \\
 \frac{5s^3}{6} + \frac{3s^2}{2} - \frac{5s}{6} - \frac{3}{2} \mbox{ if } 
s \equiv 1 \mod 2.\\ 
\end{array}\right .$$


\begin{figure}[ht]
\begin{center}
\includegraphics[width=5 cm]{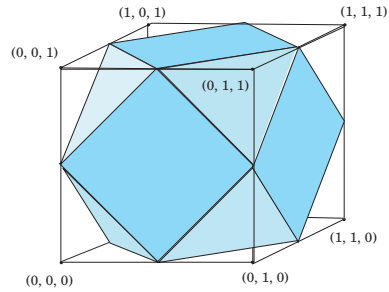} 
\caption{The cuboctahedron.} \label{tru2}
\end{center} 
\end{figure}

\end{proposition}

\begin{proposition}
The Ehrhart quasi-polynomial for the truncated regular simplex, where
the vertices are at $1/3$ and $2/3$ of the way along the simplex edges
(see Figure \ref{tru3}), is given by:

$$i_{tru\_simplex}(s) = \left \{ \begin{array}{ll}
\frac{23s^3}{81} + \frac{7s^2}{9} + \frac{13s}{9}
+ 1 \mbox{ if } s \equiv 0 \mod 3, \\
 \frac{23s^3}{81} + \frac{19s^2}{27} +
\frac{5s}{27} - \frac{95}{81} \mbox{ if } s \equiv 1 \mod 3, \\
 \frac{23s^3}{81} + \frac{17s^2}{27} +
\frac{23s}{27} + \frac{41}{81} \mbox{ if } s \equiv 2 \mod 3.\\
\end{array}\right .$$


\begin{figure}[ht]
\begin{center}
\includegraphics[width=5 cm]{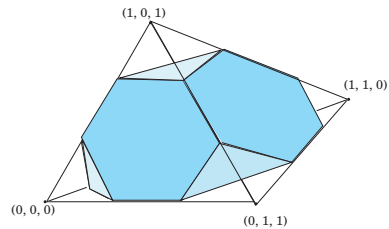} 
\caption{The truncated simplex.} \label{tru3}
\end{center} 
\end{figure}
\end{proposition}

\section{Computational results via Homogenized Barvinok's algorithm} \label{theory}

We have demonstrated the practical relevance of Barvinok's cone
decomposition approach for counting lattice points and deriving
formulas. Several other algorithms are available to carry out the same
kind of enumeration. It is important to implement them all in the same
computer system for comparison of performance and to corroborate that
the answers are correct. Some problems are solvable by some
methods but not by others.

We would like to remind a reader that one of the bottle necks of the 
original Barvinok algorithm in \citep{bar} is the fact that a polytope
may have too many vertices.  Originally, we visit all vertices of the input
polytope to compute Barvinok's short rational function and this can be
costly in terms of computational time.  For example, the well-known polytope of semi-magic cubes in
the $ 4\times 4 \times 4$ case has over two million vertices, but only
64 linear inequalities describe the polytope. 
Algorithm \ref{homogbarvinok} shows that Homogenized Barvinok's algorithm
works with only a single cone.  In this section, we will show some
practical results with Homogenized Barvinok's algorithm.

A \emph{normal semigroup} $S$ is the intersection of the lattice
$\jdlZ^d$ with a rational convex polyhedral cone in $\R^d$. 
Each pointed affine semigroup $S \subset \Z^d$ can be {\em graded}. This means that there is a linear map $\,deg: S \rightarrow \N$ with
$deg(x)=0$ if and only if $x=0$.  Given a pointed graded affine
semigroup, we define $S_r$ to be the set of elements with degree $r$, i.e.
$S_r=\{x \in S : deg(x)=r\}$.  The {\em Hilbert series} of $S$ is the
formal power series $H_S(t)=\sum_{k=0}^\infty \# (S_r)
t^r$, where $\# (S_r)$ is the cardinality of $S_r$. Algebraically, this is just the Hilbert series of the semigroup
ring $\C[S]$. It is a well-known property that $H_S$ is
represented by a rational function of the form

$$\frac{Q(t)}{(1-t^{s_1})(1-t^{s_2})\dots (1-t^{s_d})}$$

where $Q(t)$ is a polynomial of degree less than $s_1+\dots+s_d$ (see
Chapter 4 \citep{stanley}).  
The first challenge is to compute the Hilbert series of magic cubes.
Several other methods had been tried to
compute the Hilbert series explicitly (see \citep{ahmed} for
references).  One of the most well-known challenges was that of
counting the $5 \times 5$ magic squares of magic sum $n$. Similarly
several authors had tried before to compute the Hilbert series of the
$3 \times 3 \times 3 \times 3$ magic cubes. It is not difficult
to see this is equivalent to determining an Ehrhart series. Using
Algorithm \ref{homogbarvinok} we finally present the solution, which had been
inaccessible using Gr\"obner bases methods. For comparison, the reader
familiar with Gr\"obner bases computations should be aware that the
$5\times 5$ magic squares problem required a computation of a
Gr\"obner bases of a toric ideal of a matrix $A$ with 25 rows and over
4828 columns. Our attempts to handle this problem with {\tt CoCoA} and
{\tt Macaulay2} were unsuccessful.  We now give the numerator and then
the denominator of the rational functions computed with the software
{\tt LattE}:

\begin{theorem} 
\label{magic}

\item The generating function $\,\sum_{n \geq 0} f(n) t^n \,$ for the number $f(n)$ of $5\times
5$ magic squares of magic sum $n$ is given by the rational function
$p(t)/q(t)$ with numerator

{\footnotesize

$p(t)=
{t}^{76}+28\,{t}^{75}+639\,{t}^{74}+11050\,{t}^{73}+136266\,{t}^{72}+
1255833\,{t}^{71}+9120009\,{t}^{70}+54389347\,{t}^{69}+\newline
274778754\,{t}^{68}+1204206107\,{t}^{67}+4663304831\,{t}^{66}+
16193751710\,{t}^{65}+51030919095\,{t}^{64}+ \newline
147368813970\,{t}^{63}+393197605792\,{t}^{62}+
975980866856\,{t}^{61}+
2266977091533\,{t}^{60}+
4952467350549\,{t}^{59}+\newline
10220353765317\,{t}^{58}+
20000425620982\,{t}^{57}+37238997469701\,{t}^{56}+
66164771134709\,{t}^{55}+112476891429452\,{t}^{54}+\newline
183365550921732\,{t}^{53}+
287269293973236\,{t}^{52}+433289919534912\,{t}^{51}+
630230390692834\,{t}^{50}\newline +885291593024017\,{t}^{49}+
1202550133880678\,{t}^{48}+1581424159799051\,{t}^{47}+
2015395674628040\,{t}^{46}+\newline 2491275358809867\,{t}^{45}+
2989255690350053\,{t}^{44}+3483898479782320\,{t}^{43}+
3946056312532923\,{t}^{42}+\newline
4345559454316341\,{t}^{41}+
4654344257066635\,{t}^{40}+4849590327731195\,{t}^{39}+
4916398325176454\,{t}^{38}+\newline 4849590327731195\,{t}^{37}+
4654344257066635\,{t}^{36}+
4345559454316341\,{t}^{35}+
3946056312532923\,{t}^{34}+\newline 3483898479782320\,{t}^{33}+
2989255690350053\,{t}^{32}+2491275358809867\,{t}^{31}+
2015395674628040\,{t}^{30}+\newline
1581424159799051\,{t}^{29}+
1202550133880678\,{t}^{28}+885291593024017\,{t}^{27}+
630230390692834\,{t}^{26}+433289919534912\,{t}^{25}+
287269293973236\,{t}^{24}+183365550921732\,{t}^{23}+
112476891429452\,{t}^{22}+66164771134709\,{t}^{21}+37238997469701
\,{t}^{20}+\newline  20000425620982\,{t}^{19}+
10220353765317\,{t}^{18}+
4952467350549\,{t}^{17}+2266977091533\,{t}^{16}+975980866856\,{t}^
{15}+\newline
393197605792\,{t}^{14}+147368813970\,{t}^{13}+51030919095\,{t
}^{12}+16193751710\,{t}^{11}+4663304831\,{t}^{10}+1204206107\,{t}^
{9}+\newline 274778754\,{t}^{8}+54389347\,{t}^{7}+9120009\,{t}^{6}+1255833
\,{t}^{5}+136266\,{t}^{4}+11050\,{t}^{3}+639\,{t}^{2}+28\,t+1
$ \ \ \ \ \ and denominator

$ q(t)=
\left ({t}^{2}-1\right )^{10}\left ({t}^{2}+t+1\right )
^{7}\left ({t}^{7}-1\right )^{2}
\left ({t}^{6}+{t}^{3}+1\right )\left ({t}^{5}+{t}^{3}+{t}^{2}+t+1
\right )^{4}\left (1-t\right )^{3}\left ({t
}^{2}+1\right )^{4}
$.

}

The generating function $\,\sum_{n \geq 0} f(n) t^n \,$ for the number
$f(n)$ of $3\times 3 \times 3 \times 3$ magic cubes with magic sum $n$
is given the rational function $r(t)/s(t)$ where

{\footnotesize
$
r(t)={t}^{54}+150\,{t}^{51}+5837\,{t}^{48}+63127\,{t}^{45}+331124\,{t}^{42}
+1056374\,{t}^{39}+2326380\,{t}^{36}+3842273\,{t}^{33}+5055138\,{t}^{
30}+5512456\,{t}^{27}+5055138\,{t}^{24}+3842273\,{t}^{21}+2326380\,{t}
^{18}+1056374\,{t}^{15}+331124\,{t}^{12}+63127\,{t}^{9}+5837\,{t}^{6}+
150\,{t}^{3}+1
$ \ \ \ \ \ and

$
s(t)=\left ({t}^{3}+1\right )^{4}\left ({t}^{
12}+{t}^{9}+{t}^{6}+{t}^{3}+1\right )\left (1-{t}^{3}\right )^{9}
\left ({t}^{6}+{t}^{3}+1\right )
$.
}
\end{theorem}




\section{Computational results via the BBS algorithm and the digging algorithm}
\label{Computational Experiments}

In this section we report our experience solving hard
knapsack problems from \cite{aardaletal3,cuww}. See Table
\ref{knapsacks} for the data used here. Their form is $\maximize c\cdot x \mbox{ subject
to } ax=b, x\geq 0, x\in\Z^d,$ where $b\in\Z$ and where $a\in\Z^d$
with $\gcd(a_1,\ldots,a_d)=1$.  For the cost vector $c$, we choose the
first $d$ components of the vector
$(213,-1928,-11111,-2345,9123,-12834,-123,122331,0,0)$. 
We compared
the original digging algorithm, the single cone digging algorithm, and the BBS algorithm, which are implemented in {\tt LattE}
(available at \url{www.math.ucdavis.edu/~latte}), with {\tt
CPLEX} version 6.6. The computations were done on a $1$ GHz Pentium PC
running Red Hat Linux.
Table \ref{Running times} provides the optimal values and an optimal
solution for each problem. As it turns out, there is exactly one optimal
solution for each problem. 

\begin{table}[tbh]

\centering
{\tiny
\begin{tabular}{|l|rrrrrrrrrr|r|} \hline
\hspace{-4pt}Problem && &&& $a$ &&& && &$b$\\ \hline
\hspace{-4pt}cuww1 & \hspace{-6pt} 12223 & \hspace{-6pt} 12224 & \hspace{-6pt} 36674 & \hspace{-6pt} 61119 & \hspace{-6pt} 85569 &      &      &      & & & \hspace{-6pt} 89643482 \\ \hline
\hspace{-4pt}cuww2 & \hspace{-6pt} 12228 & \hspace{-6pt} 36679 & \hspace{-6pt} 36682 & \hspace{-6pt} 48908 & \hspace{-6pt} 61139 & \hspace{-6pt} 73365 &      &      & & & \hspace{-6pt} 89716839 \\ \hline
\hspace{-4pt}cuww3 & \hspace{-6pt} 12137 & \hspace{-6pt} 24269 & \hspace{-6pt} 36405 & \hspace{-6pt} 36407 & \hspace{-6pt} 48545 & \hspace{-6pt} 60683 &      &      & & & \hspace{-6pt} 58925135 \\ \hline
\hspace{-4pt}cuww4 & \hspace{-6pt} 13211 & \hspace{-6pt} 13212 & \hspace{-6pt} 39638 & \hspace{-6pt} 52844 & \hspace{-6pt} 66060 & \hspace{-6pt} 79268 & \hspace{-6pt} 92482 &      & & &  \hspace{-6pt} 104723596 \\ \hline
\hspace{-4pt}cuww5 & \hspace{-6pt} 13429 & \hspace{-6pt} 26850 & \hspace{-6pt} 26855 & \hspace{-6pt} 40280 & \hspace{-6pt} 40281 & \hspace{-6pt} 53711 & \hspace{-6pt} 53714 & \hspace{-6pt} 67141 & & & \hspace{-6pt} 45094584 \\ \hline
\hspace{-4pt}prob1 & \hspace{-6pt} 25067 & \hspace{-6pt} 49300 & \hspace{-6pt} 49717 & \hspace{-6pt} 62124 & \hspace{-6pt} 87608 & \hspace{-6pt} 88025 & \hspace{-6pt} 113673 & \hspace{-6pt} 119169 & \hspace{-6pt}  & \hspace{-6pt}  & \hspace{-6pt} 33367336 \\ \hline
\hspace{-4pt}prob2 & \hspace{-6pt}11948 & \hspace{-6pt} 23330 & \hspace{-6pt} 30635 & \hspace{-6pt} 44197 & \hspace{-6pt} 92754 & \hspace{-6pt} 123389 & \hspace{-6pt} 136951 & \hspace{-6pt} 140745 & \hspace{-6pt}  & \hspace{-6pt}  & \hspace{-6pt} 14215207 \\ \hline
\hspace{-4pt}prob3 & \hspace{-6pt}39559 & \hspace{-6pt} 61679 & \hspace{-6pt} 79625 & \hspace{-6pt} 99658 & \hspace{-6pt} 133404 & \hspace{-6pt} 137071 & \hspace{-6pt} 159757 & \hspace{-6pt} 173977 & \hspace{-6pt}  & \hspace{-6pt}  & \hspace{-6pt} 58424800 \\ \hline
\hspace{-4pt}prob4 & \hspace{-6pt}48709 & \hspace{-6pt} 55893 & \hspace{-6pt} 62177 & \hspace{-6pt} 65919 & \hspace{-6pt} 86271 & \hspace{-6pt} 87692 & \hspace{-6pt} 102881 & \hspace{-6pt} 109765 & \hspace{-6pt}  & \hspace{-6pt}  & \hspace{-6pt} 60575666\\ \hline
\hspace{-4pt}prob5 & \hspace{-6pt}28637 & \hspace{-6pt} 48198 & \hspace{-6pt} 80330 & \hspace{-6pt} 91980 & \hspace{-6pt} 102221 & \hspace{-6pt} 135518 & \hspace{-6pt} 165564 & \hspace{-6pt} 176049 & \hspace{-6pt}  & \hspace{-6pt}  & \hspace{-6pt}  62442885\\ \hline
\hspace{-4pt}prob6 & \hspace{-6pt}20601 & \hspace{-6pt} 40429 & \hspace{-6pt} 40429 & \hspace{-6pt} 45415 & \hspace{-6pt} 53725 & \hspace{-6pt} 61919 & \hspace{-6pt} 64470 & \hspace{-6pt} 69340 & \hspace{-6pt} 78539 & \hspace{-6pt} 95043 & \hspace{-6pt}  22382775 \\ \hline
\hspace{-4pt}prob7 & \hspace{-6pt}18902 & \hspace{-6pt} 26720 & \hspace{-6pt} 34538 & \hspace{-6pt} 34868 & \hspace{-6pt} 49201 & \hspace{-6pt} 49531 & \hspace{-6pt} 65167 & \hspace{-6pt} 66800 & \hspace{-6pt} 84069 & \hspace{-6pt} 137179 & \hspace{-6pt}  27267752 \\  \hline
\hspace{-4pt}prob8 & \hspace{-6pt}17035 & \hspace{-6pt} 45529 & \hspace{-6pt} 48317 & \hspace{-6pt} 48506 & \hspace{-6pt} 86120 & \hspace{-6pt} 100178 & \hspace{-6pt} 112464 & \hspace{-6pt} 115819 & \hspace{-6pt} 125128 & \hspace{-6pt} 129688 & \hspace{-6pt}  21733991 \\ \hline
\hspace{-4pt}prob9 & \hspace{-6pt}3719 & \hspace{-6pt} 20289 & \hspace{-6pt} 29067 & \hspace{-6pt} 60517 & \hspace{-6pt} 64354 & \hspace{-6pt} 65633 & \hspace{-6pt} 76969 & \hspace{-6pt} 102024 & \hspace{-6pt} 106036 & \hspace{-6pt} 119930 & \hspace{-6pt}  13385100 \\ \hline
\hspace{-4pt}prob10 &\hspace{-6pt} 45276 & \hspace{-6pt} 70778 & \hspace{-6pt} 86911 & \hspace{-6pt} 92634 & \hspace{-6pt} 97839 & \hspace{-6pt} 125941 & \hspace{-6pt} 134269 & \hspace{-6pt} 141033 & \hspace{-6pt} 147279 & \hspace{-6pt} 153525 & \hspace{-7pt} 106925262 \\ \hline
\end{tabular}
}
\centering
\caption{knapsack problems.} \label{knapsacks} 
\end{table}


\begin{sidewaystable}
{\small
\centering

\begin{tabular}{|l|r|l|r|r|r|r|}\hline
& & & Runtime for & Runtime for & Runtime for & Runtime for \\
\hspace{-6pt} Problem & \hspace{-6pt}    Value & \hspace{-6pt}                 Solution  &  Digging (Original)  &  Digging (S. Cone)   & \hspace{-8pt} BBS & \hspace{-6pt} CPLEX 6.6\\ \hline
\hspace{-6pt} cuww1 & \hspace{-6pt}   1562142 & \hspace{-6pt} [7334 0 0 0 0]            & \hspace{-6pt}     0.4 sec. &  0.17 sec. & \hspace{-8pt}     414 sec. & \hspace{-6pt} $>$ 1.5h (OM) \\ \hline
\hspace{-6pt} cuww2 & \hspace{-6pt}  -4713321 & \hspace{-6pt} [3 2445 0 0 0 0]          & \hspace{-6pt} $>$ 3.5h  &  $>$ 3.5h   & \hspace{-8pt}   6,600 sec. & \hspace{-6pt} $>$ 0.75h (OM)\\ \hline
\hspace{-6pt} cuww3 & \hspace{-6pt}   1034115 & \hspace{-6pt} [4855 0 0 0 0 0]          & \hspace{-6pt}     1.4 sec. & 0.24 sec. &\hspace{-8pt}   6,126 sec. & \hspace{-6pt} $>$ 0.75h (OM)\\ \hline
\hspace{-6pt} cuww4 & \hspace{-6pt} -29355262 & \hspace{-6pt} [0 0 2642 0 0 0 0]        & \hspace{-6pt} $>$ 1.5h  &   $>$ 1.5h   & \hspace{-8pt}  38,511 sec. & \hspace{-6pt} $>$ 0.75h (OM)\\ \hline
\hspace{-6pt} cuww5 & \hspace{-6pt}  -3246082 & \hspace{-6pt} [1 1678 1 0 0 0 0 0]      & \hspace{-6pt} $>$ 1.5h  & 147.63 sec.    & \hspace{-8pt}  $>$ 80h        & \hspace{-6pt} $>$ 0.75h (OM)\\ \hline
\hspace{-6pt} prob1 & \hspace{-6pt}   9257735 & \hspace{-6pt} [966 5 0 0 1 0 0 74]      & \hspace{-6pt}    51.4 sec. & 18.55 sec. & \hspace{-8pt}  $>$  3h        & \hspace{-6pt} $>$ 1h (OM) \\ \hline
\hspace{-6pt} prob2 & \hspace{-6pt}   3471390 & \hspace{-6pt} [853 2 0 4 0 0 0 27]      & \hspace{-6pt}    24.8 sec. & 6.07 sec & \hspace{-8pt}  $>$ 10h        & \hspace{-6pt} $>$ 0.75h (OM)\\ \hline
\hspace{-6pt} prob3 & \hspace{-6pt}  21291722 & \hspace{-6pt} [708 0 2 0 0 0 1 173]     & \hspace{-6pt}    48.2 sec. & 9.03  sec. & \hspace{-8pt}  $>$ 12h        & \hspace{-6pt} $>$ 1.5h (OM) \\ \hline
\hspace{-6pt} prob4 & \hspace{-6pt}   6765166 & \hspace{-6pt} [1113 0 7 0 0 0 0 54]     & \hspace{-6pt}    34.2 sec. & 9.61 sec. & \hspace{-8pt}  $>$  5h        & \hspace{-6pt} $>$ 1.5h (OM) \\ \hline
\hspace{-6pt} prob5 & \hspace{-6pt}  12903963 & \hspace{-6pt} [1540 1 2 0 0 0 0 103]    & \hspace{-6pt}    34.5 sec. & 9.94 sec. & \hspace{-8pt}  $>$  5h        & \hspace{-6pt} $>$ 1.5h (OM) \\ \hline
\hspace{-6pt} prob6 & \hspace{-6pt}   2645069 & \hspace{-6pt} [1012 1 0 1 0 1 0 20 0 0] & \hspace{-6pt}   143.2 sec. & 19.21 sec. & \hspace{-8pt}  $>$  4h        & \hspace{-6pt} $>$ 2h (OM)\\ \hline
\hspace{-6pt} prob7 & \hspace{-6pt}  22915859 & \hspace{-6pt} [782 1 0 1 0 0 0 186 0 0] & \hspace{-6pt}   142.3 sec.  & 12.84 sec. & \hspace{-8pt}  $>$  4h        & \hspace{-6pt} $>$ 1h (OM)\\ \hline
\hspace{-6pt} prob8 & \hspace{-6pt}   3546296 & \hspace{-6pt} [1 385 0 1 1 0 0 35 0 0]  & \hspace{-6pt}   469.9 sec. & 49.21 sec. & \hspace{-8pt}  $>$  3.5h      & \hspace{-6pt} $>$ 2.5h (OM)\\ \hline
\hspace{-6pt} prob9 & \hspace{-6pt}  15507976 & \hspace{-6pt} [31 11 1 1 0 0 0 127 0 0] & \hspace{-6pt} 1,408.2 sec. & 283.34 sec. & \hspace{-8pt}  $>$ 11h        & \hspace{-6pt} 4.7 sec. \\ \hline
\hspace{-6pt} prob10& \hspace{-6pt}  47946931 & \hspace{-6pt} [0 705 0 1 1 0 0 403 0 0] & \hspace{-6pt}   250.6 sec. & 29.28 sec. & \hspace{-8pt}  $>$ 11h        & \hspace{-6pt} $>$ 1h  (OM)\\ \hline
\end{tabular}
\caption{{\small Optimal values, optimal solutions, and running times for each
  problem. OM:= Out of memory. }} \label{Running times}\label{IPtable4}

\centering
}
\end{sidewaystable}

\begin{table}
\begin{center}
\begin{tabular}{|c|c|c|c|c|}\hline
problem & Original & Original & Single Cone & Single Cone \\ 
 & Digging (A) &  Digging (B) & Digging (A) & Digging (B) \\ \hline
cuww 1& 110 & 0 & 25 & 0 \\ \hline
cuww 2 & 386 & $>$ 2,500,000 & 79 & $>$ 2,500,000 \\ \hline
cuww 3 & 346 & 0 & 49 & 0 \\ \hline
cuww 4 & 364 & $>$ 400,000 & 51 &$>$ 400,000 \\ \hline
cuww 5 & 2,514 & $>$ 100,000 & 453 & 578,535 \\ \hline
prob 1 & 10,618 &  74,150 & 1,665 & 74,150  \\ \hline
prob 2 &  6,244 & 0 & 806 & 0 \\ \hline
prob 3 &  12,972 & 0 & 2,151 & 0 \\ \hline
prob 4 &  9,732 & 0 & 1,367 & 0  \\ \hline
prob 5 &  8,414 & 1 & 2,336 & 1  \\ \hline
prob 6 &  26,448  & 5 & 3,418 & 5  \\ \hline
prob 7 &  20,192 & 0 & 2,015 & 0  \\ \hline
prob 8 & 62,044 & 0 & 6,523 & 0    \\ \hline
prob 9 &   162,035  & 3,558 &   45,017 & 3,510 \\ \hline
prob 10 &  38,638  & 256 & 5,128 & 256 \\ \hline
\end{tabular}
\caption{Data for the digging algorithm. $A:=$ number of unimodular cones and $B:=$ number of digging levels}\label{IPtable5}
\end{center}
\end{table}

With one exception, CPLEX 6.6. could not solve the given problems. Note that
whenever the digging algorithm found the optimal value, it did so much
faster than the BBS algorithm. This is interesting, as the worst-case
complexity for the digging algorithm is exponential even for fixed
dimension, while the BBS has polynomial complexity in fixed dimension.
The digging algorithm fails to find a solution for problems cuww2,
cuww3, and cuww5. What happens is that the expansion step becomes
costly when more coefficients have to be computed. In these three
examples, we computed coefficients for more than 2,500,000, 400,000,
and 100,000 powers of $t$; all turning out to be $0$. The Digging
algorithm is slower than CPLEX in problem prob9 because during the
execution of Barvinok's unimodular cone decomposition (see pages 15
and 16 of \cite{BarviPom}) more than 160,000 cones are generated,
leading to an enormous rational function for $f(P;t)$. Moreover, for
prob9 more than 3,500 coefficients turned out to be $0$, before a
non-zero leading coefficient was detected.  Finally, in problems
cuww1, cuww3, prob2, prob3, prob4, prob6, and prob8, no digging was
necessary at all, that is, Lasserre's condition did not fail here. For
all other problems, Lasserre's condition did fail and digging steps
were necessary to find the first non-vanishing coefficient in the
expansion of $f(P;t)$.

\begin{table}[tbh]
\begin{center}

\begin{tabular}{|c|c|c|c|}\hline
problem & (C) & (D)& (E)\\ \hline
cuww 1& 125,562 & 4,829.3 & 26 \\ \hline
cuww 2 & 2,216,554 & 88,662.16 & 25  \\ \hline
cuww 3 & 2,007,512 & 80,300.48 & 25  \\ \hline
cuww 4 & 10,055,730 & 402,229.2 &  25 \\ \hline
cuww 5 & NA & NA & $\leq$ 26\\ \hline
\end{tabular}
\caption{Data for the BBS algorithm. $C:=$ Total num. of unimodular cones, $D:=$ Average num. of unimodular cones per an iteration, $E:=$ Total number of iterations, and $NA:=$ not available. }\label{IPtable2}

\end{center}
\end{table}

\begin{table}
\begin{center}
{\small
\begin{tabular}{|c|c|c|c|}\hline
problem & The optimal value & An optimal solution & (N) \\ \hline
cuww 1& 1562142 & [7334 0 0 0 0] & 1 \\ \hline
cuww 2 & -4713321 & [3 2445 0 0 0 0] & 1 \\ \hline
cuww 3 & 1034115 & [4855 0 0 0 0 0] & 1 \\ \hline
cuww 4 &  -29355262 & [0 0 2642 0 0 0 0] & 1\\ \hline
cuww 5 &  -3246082 & [1 1678 1 0 0 0 0 0] & 1 \\ \hline
prob 1 & 9257735 & [966 5 0 0 1 0 0 74] & 1 \\ \hline
prob 2 &  3471390 & [853 2 0 4 0 0 0 27] & 1 \\ \hline
prob 3 &  21291722& [708 0 2 0 0 0 1 173] & 1\\ \hline
prob 4 &  6765166 &  [1113 0 7 0 0 0 0 54] & 1 \\ \hline
prob 5 & 12903963 & [1540 1 2 0 0 0 0 103]  & 1 \\ \hline
prob 6 & 2645069 & [1012 1 0 1 0 1 0 20 0 0]  & 1 \\ \hline
prob 7 &  22915859 & [782 1 0 1 0 0 0 186 0 0]  & 1 \\ \hline
prob 8 &  3546296 & [1 385 0 1 1 0 0 35 0 0]  & 1 \\ \hline
prob 9 &  15507976 & [31 11 1 1 0 0 0 127 0 0] & 1 \\ \hline
prob 10 & 47946931 & [0 705 0 1 1 0 0 403 0 0] & 1 \\ \hline
\end{tabular}
\caption{The optimal value and an optimal solution for each problem. $N:=$ num. of solutions}\label{IPtable3}
}
\end{center}
\end{table}

\newpage
\appendix
\pagestyle{myheadings} 
\markright{  \rm \normalsize APPENDIX. \hspace{0.5cm}
 User manual of {\tt LattE}}

\chapter{User manual of {\tt LattE}}\label{App}

\thispagestyle{myheadings} 

\section{Introduction}

\subsection{What is {\tt LattE}?} \label{intro}

The name ``{\tt LattE}'' is an abbreviation for ``{\bf Latt}ice point 
{\bf E}numeration.'' So what exactly does {\tt LattE} do? The 
software's main function is to count the lattice points contained in 
convex polyhedra defined by linear equations and inequalities with 
integer coefficients. The polyhedra can be of any (reasonably small) 
dimension, and {\tt LattE} uses an algorithm that runs in polynomial
time for fixed dimension: Barvinok's algorithm \cite{BarviPom}. To
learn more about the exact details of our implementation and
algorithmic techniques involved, the interested reader can consult 
\cite{latte,newlatte,latte3} and the references listed therein. Here we
give a rather short description of the mathematical objects used by
{\tt LattE}, {\bf Barvinok's Rational Functions}:

\noindent
Given a convex polyhedron $P = \{u\in\R^d:Au\leq b\}$, where $A$ and
$b$ are integral, the fundamental object that we compute is a short 
representation of the infinite power series:
\[
f(P;x) \quad = \sum_{\alpha\in P\cap\Z^d} x_1^{\alpha_1}
x_2^{\alpha_2} \ldots x_d^{\alpha_d}.
\]
Here each lattice point is given by one monomial. Note that this can be 
a rather long sum, in fact for a polyhedral cone it can be infinite, but 
the good news is that it admits short representations.

\noindent {\bf Example:} Let $P$ be the quadrangle with vertices 
$V_1=(0,0)$, $V_2=(5,0)$, $V_3=(4,2)$, and $V_4=(0,2)$. 

\begin{figure}[thb]
\begin{center}
     \includegraphics[scale=.79]{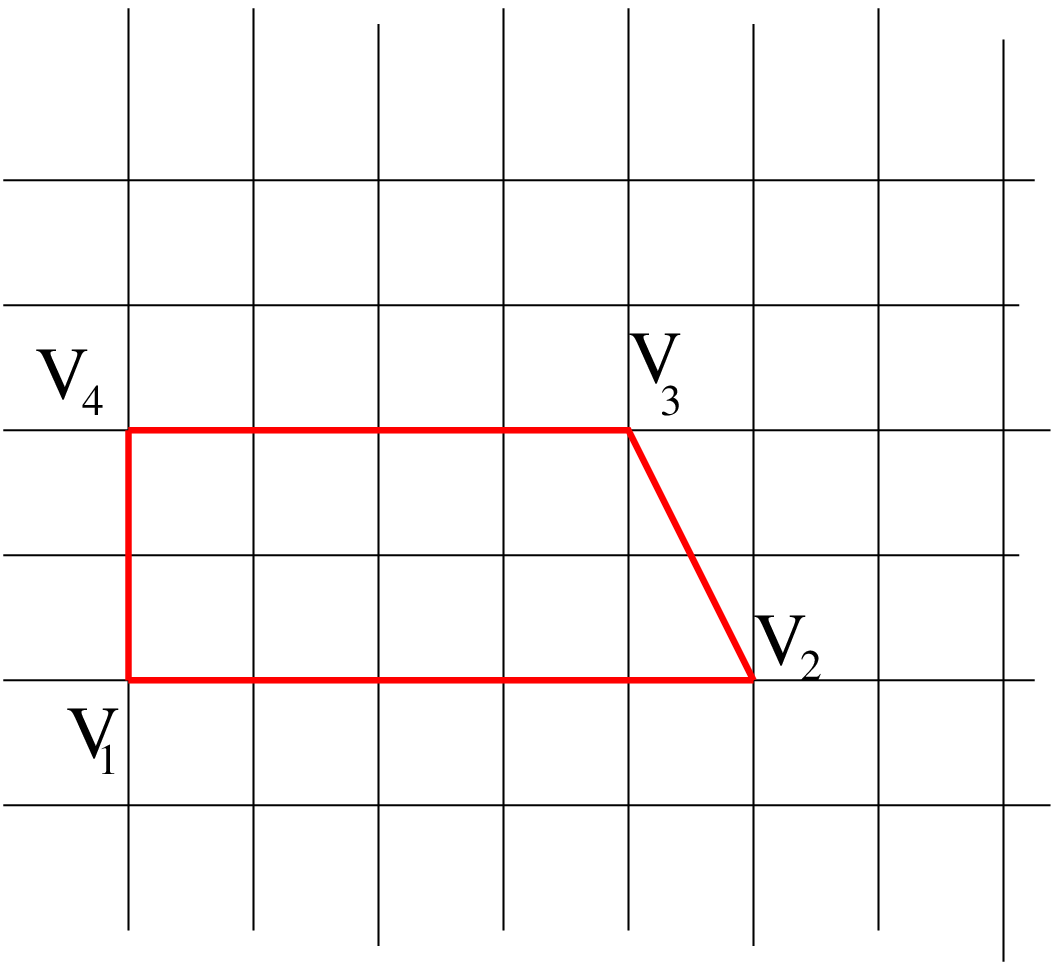}
\end{center}
\end{figure}

{\small
\noindent
$f(P;x,y)={x}^{5}+{x}^{4}y+{x}^{4}+{x}^{4}{y}^{2}+y{x}^{3}+{x}^{3}+
{x}^{3}{y}^{2}+y{x}^{2}+{x}^{2}+{x}^{2}{y}^{2}+xy+x+x{y}^{2}+y+1+
{y}^{2}$
}

The fundamental theorem of Barvinok (circa 1993, see \cite{BarviPom})
says that you can write $f(P;x)$ as a sum of short rational functions,
in polynomial time when the dimension of the polyhedron is fixed.
In our running example we easily see that the 16 monomial polynomial
can be written as shorter rational function sum:

\noindent $f(P;x,y)=f(K_{V_1};x,y)+f(K_{V_2};x,y)+f(K_{V_3};x,y)+f(K_{V_4};x,y)$ 

\noindent where 

$f(K_{V_1};x,y)={\frac {1}{\left (1-x\right )\left (1-y\right )}} \quad f(K_{V_2};x,y)=\frac{({x}^{5}+{x}^{4}y)}{ (1-{x}^{-1}) (1-y^2x^{-1})}$

$f(K_{V_3};x,y)=\frac{({x}^{4}{y}^{2}+{x}^{4})}{ (1-{x}^{-1})
(1-xy^{-2})}
 \quad  f(K_{V_4};x,y)=\frac{y^{2}}{(1-{y}^{-1} )(1-x )}$
\vskip .5cm

$ f(P; 1,1)=16$

Counting the lattice points in convex polyhedra is a powerful tool which 
allows many applications in areas such as Combinatorics, Statistics, 
Optimization, and Number Theory. 

\subsection{What can {\tt LattE} compute?}

In the following we list the operations that {\tt LattE} v1.1 can
perform on bounded convex polyhedra (more commonly referred to as 
\textit{polytopes}). For the reader's convenience, we already include
the basic commands to actually do the tasks. Let us assume that a
description of a polytope $P$ is given in the file ``fileName'' (see
Section \ref{Input Files} for format) and that a cost vector is
specified in the file ``fileName.cost'' (needed for the optimization
part, see Section \ref{Input Files} for format).

Tasks performed by {\tt LattE} v1.1:
\begin{enumerate}
\item Count the number of lattice points in $P$.
\begin{verbatim}
     ./count fileName
\end{verbatim} 
\item Count the number of lattice points in $nP$, the dilation of $P$
  by the integer factor $n$.
\begin{verbatim}
     ./count dil n fileName
\end{verbatim} 
\item Calculate a rational function that encodes the \textit{Ehrhart
  series} associated with the polytope. By definition, the $n$-th
  coefficient in the Ehrhart series equals the number of lattice
  points in $nP$. For more details on Ehrhart counting functions see,
  for example, Chapter 4 of \cite{stanley}. 
\begin{verbatim}
     ./ehrhart fileName
\end{verbatim} 
\item Calculate the first $n+1$ terms of the Ehrhart series associated
  with the polytope.  
\begin{verbatim}
     ./ehrhart n fileName
\end{verbatim} 
\item Maximize or minimize a given linear function of the lattice
  points in $P$.
\begin{verbatim}
     ./maximize fileName
     ./minimize fileName
\end{verbatim} 
\end{enumerate}	

In addition to these basic functions, there are more specific calls to
{\tt LattE}. For example to use the homogenized Barvinok algorithm
instead of the original one in order to count the lattice points. These
details will be explained in Section \ref{Running LattE}.

\newpage

\section{Downloading and Installing {\tt LattE}}

{\tt LattE} is downloadable from the following website:

\makebox[12 cm]{http://www.math.ucdavis.edu/$\sim$latte/downloads/}

\textbf{Step 1: Create directory for {\tt LattE}}

\begin{verbatim}
     mkdir latte
\end{verbatim}

\textbf{Step 2: Download ``latte\_v1.1.tar.gz'' to directory ``latte''}

\begin{verbatim}
     Download ``latte_v1.1.tar.gz'' from
\end{verbatim}

\makebox[12 cm]{http://www.math.ucdavis.edu/$\sim$latte/downloads/}
 
(If you have never downloaded a file from the internet: A click with
your right mouse button onto the file name on the webpage should do
the trick. In any case, if you do not succeed, ask your system
administrator, a friend, or send us an email.) 

\textbf{Step 3: Change to directory for ``latte''}

\begin{verbatim}
     cd latte
\end{verbatim}

\textbf{Step 4: Unzip and untar the archive} 

\begin{verbatim}
     gunzip latte_v1.1.tar.gz
     tar xvf latte_v1.1.tar
\end{verbatim}

\textbf{Step 5: Make ``install'' executable} 

\begin{verbatim}
     chmod 700 install
\end{verbatim}

\textbf{Step 6: Install {\tt LattE}} 

\begin{verbatim}
     ./install
\end{verbatim}

\newpage

\section{Input Files}\label{Input Files}

\subsection{{\tt LattE} Input Files}

\subsubsection{Inequality Description}
For computations involving a polytope $P$ described by a
system of inequalities $Ax\leq b$, where $A\in\Z^{m\times d}$, 
$A=(a_{ij})$, and $b\in\Z^m$, the {\tt LattE} readable input file
would be as follows: 
\begin{verbatim}
m d+1
b  -A
\end{verbatim}

\textbf{EXAMPLE.}
Let $P=\{(x,y): x\leq 1, y\leq 1, x+y\leq 1, x\geq 0, y\geq 0\}$.
Thus
\[
\begin{array}{ccc}
A=\left(
\begin{array}{rr} 
 1 &  0 \\ 
 0 &  1 \\ 
 1 &  1 \\
-1 &  0 \\ 
 0 & -1 \\ 
\end{array} 
\right) 
& , &
b = \left( 
\begin{array}{r} 
1 \\ 
1 \\ 
1 \\ 
0 \\
0 \\ 
\end{array} 
\right)
\end{array}
\]
and the {\tt LattE} input file would be as such:
\begin{verbatim}
5 3
1 -1  0
1  0 -1
1 -1 -1
0  1  0
0  0  1
\end{verbatim}

\subsubsection{Equations}
In {\tt LattE}, polytopes are represented by {\bf linear constraints},
i.e. equalities or inequalities. By default a constraint is an
inequality of type $ax\leq b$ unless we specify, by using a single
additional line, the line numbers of constraints that are linear
equalities. 

\textbf{EXAMPLE.}
Let $P$ be as in the previous example, but require $x+y=1$ instead of
$x+y\leq 1$, thus, 
$P=\{(x,y): x\leq 1, y\leq 1, x+y=1, x\geq 0, y\geq 0\}$.
Then the {\tt LattE} input file that describes $P$ would be as such:
\begin{verbatim}
5  3
1 -1  0
1  0 -1
1 -1 -1
0  1  0
0  0  1
linearity 1 3
\end{verbatim}
The last line states that among the $5$ inequalities one is to be
considered an equality, the third one.

\subsubsection{Nonnegativity Constraints}
For bigger examples it quickly becomes cumbersome to state all
nonnegativity constraints for the variables one by one. Instead, you
may use another short-hand.

\textbf{EXAMPLE.}
Let $P$ be as in the previous example, then the {\tt LattE} input file
that describes $P$ could also be described as such: 
\begin{verbatim}
3  3
1 -1  0
1  0 -1
1 -1 -1
linearity 1 3
nonnegative 2 1 2
\end{verbatim}
The last line states that there are two nonnegativity constraints and
that the first and second variables are required to be nonnegative. 
{\bf NOTE} that the first line reads ``3 3'' and not ``5 3'' as above! 

\subsubsection{Cost Vector}
The functions maximize and minimize solve the integer linear programs
\[
\max\{c^\intercal x: x\in P\cap\Z^d\}
\]
and
\[
\min\{c^\intercal x: x\in P\cap\Z^d\}.
\]
Besides a description of the polyhedron $P$, these functions need a
linear objective function given by a certain cost vector $c$. If the
polyhedron is given in the file ``fileName''
\begin{verbatim}
4  4
1 -1  0  0
1  0 -1  0
1  0  0 -1
1 -1 -1 -1
linearity 1 4
nonnegative 3 1 2 3
\end{verbatim}
the cost vector must be given in the file ``fileName.cost'', as for
example in the following three-dimensional problem: 
\begin{verbatim}
1 3
2 4 7
\end{verbatim}
The first two entries state the size of a $1\times n$ matrix (encoding
the cost vector), followed by the $1\times n$ matrix itself. Assuming
that we call maximize, this whole data encodes the integer program
\[
\max\{2x_1+4x_2+7x_3: x_1+x_2+x_3=1, x_1,x_2,x_3\in\{0,1\}\}.
\]

\subsection{{\tt cdd} Input Files}
In addition to the formats described above, {\tt LattE} can also
accept input files in standard {\tt cdd} format. (See Subsection 
\ref{Command Syntax} for
details on how to run {\tt LattE} on a {\tt cdd} input file.) Below is
an example of {\tt cdd} input that is readable into {\tt LattE}.
\begin{verbatim}
H-representation
begin
4  4  integer
2 -2  4 -1
3 -2 -2  3
6  2 -4 -3
1  2  2  1
end
\end{verbatim}
It is important to note that {\tt LattE} can only read 
{\em integer} input.  Clearly, {\tt cdd}'s rational data files can be
converted into integer files by multiplying by the right constants. In
the packaged release of {\tt LattE} we include a binary version of
{\tt cdd}.

\newpage

\section{Running {\tt LattE}}\label{Running LattE}
 	
\subsection{Command Syntax}\label{Command Syntax}
The basic syntax to invoke the various functions of {\tt LattE} is:
\begin{verbatim}
./count fileName
./ehrhart fileName
./maximize fileName
./minimize fileName
\end{verbatim}
Note that the last two functions require a cost vector specified in
the file ``fileName.cost''!

Additionally, a variety of options can be used. All options should be
space-delimited in the command. 

One option that can be set in addition to the options given below is
``cdd'' which tells {\tt LattE} to read its input from a {\tt cdd}
input file. Thus, the above invocations for {\tt cdd} input files would be
\begin{verbatim}
./count cdd fileName
./ehrhart cdd fileName
./maximize cdd fileName
./minimize cdd fileName
\end{verbatim}

\subsection{Counting}
\begin{itemize}
\item Count the number of lattice points in polytope $P$, where $P$
  is given in ``fileName''.
\begin{verbatim}
     ./count fileName
\end{verbatim} 
\item Count the number of lattice points in $nP$, the dilation of $P$
  by the integer factor $n$.
\begin{verbatim}
     ./count dil n fileName
\end{verbatim} 
\item Count the number of lattice points in the interior of the
  polytope $P$, where $P$ is given in ``fileName''.
\begin{verbatim}
     ./count int fileName
\end{verbatim} 
\item Use the homogenized Barvinok algorithm \cite{latte3} to count
  the number of lattice points in the polytope $P$, where $P$ is given
  in ``fileName''. Use if number of vertices of $P$ is big compared to
  the number of constraints. 
\begin{verbatim}
     ./count homog fileName
\end{verbatim} 
\end{itemize}

\subsection{Ehrhart Series}
\begin{itemize}
\item Compute the Ehrhart series encoded as a rational function for the
  polytope given in ``fileName''. Writes the unsimplified rational function
  to file ``fileName.rat''.
\begin{verbatim}
     ./ehrhart fileName
\end{verbatim} 
\item Compute the Ehrhart series encoded as a rational function for
  the polytope given in ``fileName''. {\bf NEEDS} Maple for
  simplification of terms. Writes the simplified rational function to
  file ``fileName.rat''. 
\begin{verbatim}
     ./ehrhart simplify fileName
\end{verbatim} 
\item Compute the Taylor series expansion of Ehrhart generating
  function up to degree $n$ for the polytope given in ``fileName''. 
\begin{verbatim}
     ./ehrhart n fileName
\end{verbatim} 
\end{itemize}

\subsection{Optimizing}
This functions {\bf NEEDS} a cost vector specified in ``fileName.cost''!!!
\begin{itemize}
\item Maximizes/Minimizes given linear cost function over the lattice
  points in the polytope given in ``fileName''. Digging algorithm
  \cite{latte3} is used. Optimal point and optimal value is returned. 
\begin{verbatim}
     ./maximize fileName
     ./minimize fileName
\end{verbatim} 
\item Maximizes/Minimizes given linear cost function over the lattice
  points in the polytope given in ``fileName''. Binary search
  algorithm is used. Only optimal value is returned. 
\begin{verbatim}
     ./maximize bbs fileName
     ./minimize bbs fileName
\end{verbatim} 
\end{itemize}

\newpage

\section{A Brief Tutorial}
In this section we invite the reader to follow along a few examples
that show how to use {\tt LattE} and also how to counter-check
results.

\subsection{Counting Magic Squares}
Our first example deals with counting magic $4\times 4$ squares. We 
call a $4\times 4$ array of nonnegative numbers a magic square if the
sums of the $4$ entries along each row, along each column and along
the two main diagonals equals the same number $s$, the magic
constant. Let us start with counting magic $4\times 4$ squares that
have the magic constant $1$. Associating variables $x_1,\ldots,x_{16}$ with
the $16$ entries, the conditions of a magic $4\times 4$ square of
magic sum $1$ can be encoded into the following input file
``EXAMPLES/magic4x4'' for {\tt LattE}.
\begin{verbatim}
10 17
1 -1 -1 -1 -1  0  0  0  0  0  0  0  0  0  0  0  0
1  0  0  0  0 -1 -1 -1 -1  0  0  0  0  0  0  0  0
1  0  0  0  0  0  0  0  0 -1 -1 -1 -1  0  0  0  0
1  0  0  0  0  0  0  0  0  0  0  0  0 -1 -1 -1 -1
1 -1  0  0  0 -1  0  0  0 -1  0  0  0 -1  0  0  0
1  0 -1  0  0  0 -1  0  0  0 -1  0  0  0 -1  0  0
1  0  0 -1  0  0  0 -1  0  0  0 -1  0  0  0 -1  0
1  0  0  0 -1  0  0  0 -1  0  0  0 -1  0  0  0 -1
1 -1  0  0  0  0 -1  0  0  0  0 -1  0  0  0  0 -1
1  0  0  0 -1  0  0 -1  0  0 -1  0  0 -1  0  0  0
linearity 10 1 2 3 4 5 6 7 8 9 10
nonnegative 16 1 2 3 4 5 6 7 8 9 10 11 12 13 14 15 16
\end{verbatim}
Now we simply invoke the counting function of {\tt LattE} by typing:
\begin{verbatim}
    ./count EXAMPLES/magic4x4
\end{verbatim}

The last couple of lines that {\tt LattE} prints to the screen
look as follows:
\begin{verbatim}
Total Unimodular Cones: 418
Maximum number of simplicial cones in memory at once: 27

*****  Total number of lattice points: 8  ****

Computation done.
Time: 1.24219 sec
\end{verbatim}
Therefore, there are exactly $8$ magic $4\times 4$ squares that
have the magic constant $1$. This is not yet impressive, as we could
have done that by hand. Therefore, let us try and find the
corresponding number for the magic constant $12$. Since this problem
is a dilation (by factor $12$) of the original problem, we do not have
to create a new file. Instead, we use the option ``dil'' to indicate
that we want to count the number of lattice points of a dilation of
the given polytope:
\begin{verbatim}
    ./count dil 12 EXAMPLES/magic4x4
\end{verbatim}
The last couple of lines that {\tt LattE} prints to the screen
look as follows:
\begin{verbatim}
Total Unimodular Cones: 418
Maximum number of simplicial cones in memory at once: 27

*****  Total number of lattice points: 225351  ****

Computation done.
Time: 1.22656 sec
\end{verbatim}
Therefore, there are exactly $225351$ magic $4\times 4$ squares that
have the magic constant $12$. (We would NOT want to do THAT one by
hand, would we?!) 

Here is some amazing observation: the running time of {\tt LattE}
is roughly the same for counting magic squares of sum $1$ and of sum
$12$. This phenomenon is due to the fact that the main part of the
computation, the creation of the generating function that encodes all
lattice points in the polytope, is nearly identical in both cases.

Although we may be already happy with these simple counting results,
let us be a bit more ambitious and and let us find a counting formula
that, for given magic sum $s$, returns the number of magic 
$4\times 4$ squares that have the magic constant $s$.

For this, simply type (note that {\tt LattE} invokes {\tt Maple} to
simplify intermediate expressions):
\begin{verbatim}
    ./ehrhart simplify EXAMPLES/magic4x4
\end{verbatim}
The last couple of lines that {\tt LattE} prints to the screen looks
as follows:
\begin{verbatim}
Rational function written to EXAMPLES/magic4x4.rat

Computation done. 
Time: 0.724609 sec
\end{verbatim}
We are informed that this call created a file ``EXAMPLES/magic4x4.rat''
containing the Ehrhart series as a rational function:
{\small
\begin{verbatim}
(t^8+4*t^7+18*t^6+36*t^5+50*t^4+36*t^3+18*t^2+4*t+1)/(-1+t)^4/(-1+t^2)^4
\end{verbatim}
}
Now we could use {\tt Maple} (or your favorite computer algebra
software) to find a series expansion of this expression. 
\begin{eqnarray*}
& & 
\frac{t^8+4*t^7+18*t^6+36*t^5+50*t^4+36*t^3+18*t^2+4*t+1}{(-1+t)^4(-1+t^2)^4}\\
& = & 1+8t^1+48t^2+200t^3+675t^4+1904t^5+4736t^6+10608t^7+21925t^8+\\
& & 42328t^9+77328t^{10}+134680t^{11}+225351t^{12}+364000t^{13}+570368t^{14}+\\
& & 869856t^{15}+{O}(t^{16})\\
\end{eqnarray*}
The summands $8t$ and $225351t^{12}$ reconfirm our previous
counts.

Although this rational function encodes the full Ehrhart series, it is
not always as easy to compute as for magic $4\times 4$ squares. As it
turns out, adding and simplifying rational functions, although in just
one variable $t$, can be extremely costly due to the high powers in
$t$ and due to long integer coefficients that appear.

However, even if we cannot compute the full Ehrhart series, we can at
least try and find the first couple of terms of it. 
\begin{verbatim}
    ./ehrhart 15 EXAMPLES/magic4x4
\end{verbatim}
The last couple of lines that {\tt LattE} prints to the screen
look as follows:
\begin{verbatim}
Memory Save Mode: Taylor Expansion:
1
8t^1
48t^2
200t^3
675t^4
1904t^5
4736t^6
10608t^7
21925t^8
42328t^9
77328t^10
134680t^11
225351t^12
364000t^13
570368t^14
869856t^15
Computation done.
Time: 1.83789 sec
\end{verbatim}
Again, our previous counts are reconfirmed.

Nice, but the more terms we want to compute the more time-consuming
this task becomes. Clearly, if we could find sufficiently many
terms, we could compute the full Ehrhart series expansion in terms of
a rational function by interpolation.

\subsection{Counting Lattice Points in the $24$-Cell}
Our next example deals with a well-known combinatorial object, the
$24$-cell. Its description is given in the file ``EXAMPLES/24\_cell'':
\begin{verbatim}
24 5 
2 -1  1 -1 -1
1  0  0 -1  0
2 -1  1 -1  1
2 -1  1  1  1
1  0  0  0  1
1  0  1  0  0
2  1 -1  1 -1
2  1  1 -1  1
2  1  1  1  1
1  1  0  0  0
2  1  1  1 -1
2  1  1 -1 -1
2  1 -1  1  1
2  1 -1 -1  1
2  1 -1 -1 -1
1  0  0  1  0
2 -1  1  1 -1
1  0  0  0 -1
2 -1 -1  1 -1
1  0 -1  0  0
2 -1 -1  1  1
2 -1 -1 -1  1
2 -1 -1 -1 -1
1 -1  0  0  0
\end{verbatim}
Now we invoke the counting function of {\tt LattE} by typing:
\begin{verbatim}
    ./count EXAMPLES/24_cell
\end{verbatim}
The last couple of lines that {\tt LattE} prints to the screen
look as follows:
\begin{verbatim}
Total Unimodular Cones: 240
Maximum number of simplicial cones in memory at once: 30

*****  Total number of lattice points: 33  ****

Computation done. 
Time: 0.429686 sec
\end{verbatim}
Therefore, there are exactly $33$ lattice points in the $24$-cell. We
get the same result by using the homogenized Barvinok algorithm:
\begin{verbatim}
    ./count homog EXAMPLES/24_cell
\end{verbatim}
The last couple of lines that {\tt LattE} prints to the screen
look as follows:
\begin{verbatim}
Memory Save Mode: Taylor Expansion:

****  Total number of lattice points is: 33  ****

Computation done. 
Time: 0.957031 sec
\end{verbatim}
But how many of these $33$ points lie in the interior of the
$24$-cell?
\begin{verbatim}
    ./count int EXAMPLES/24_cell
\end{verbatim}
The last couple of lines that {\tt LattE} prints to the screen
look as follows:
\begin{verbatim}
Reading .ext file...


*****  Total number of lattice points: 1 ****
\end{verbatim}
Therefore, there only one of the $33$ lattice points in the $24$-cell
lies in the interior.

\subsection{Maximizing Over a Knapsack Polytope}
Finally, let us solve the problem ``cuww1'' \cite{cuww,latte3}. Its
description is given in the file ``EXAMPLES/cuww1'': 
\begin{verbatim}
1 6
89643482 -12223 -12224 -36674 -61119 -85569
linearity 1 1
nonnegative 5 1 2 3 4 5
\end{verbatim}
The cost function can be found in the file ``EXAMPLES/cuww1.cost'':
\begin{verbatim}
1 5
213 -1928 -11111 -2345 9123 
\end{verbatim}
Now let us maximize this cost function over the given knapsack
polytope. Note that by default, the digging algorithm as described in
\cite{latte3} is used.
\begin{verbatim}
    ./maximize EXAMPLES/cuww1
\end{verbatim}
The last couple of lines that {\tt LattE} prints to the screen
look as follows:
{\small
\begin{verbatim}
Finished computing a rational function. 
Time: 0.158203 sec.

There is one optimal solution. 		

No digging.
An optimal solution for [213 -1928 -11111 -2345 9123] is: [7334 0 0 0 0].
The projected down opt value is: 191928257104
The optimal value is: 1562142.
The gap is: 7995261.806
Computation done.
Time: 0.203124 sec.
\end{verbatim}
}
The solution $(7334,0,0,0,0)$ is quickly found. Now let us try to
find the optimal value again by a different algorithm, the binary
search algorithm.
\begin{verbatim}
    ./maximize bbs EXAMPLES/cuww1
\end{verbatim}
The last couple of lines that {\tt LattE} prints to the screen
look as follows:
\begin{verbatim}
Total of Iterations: 26
The total number of unimodular cones: 125562
The optimal value: 1562142

The number of optimal solutions: 1
Time: 0.042968
\end{verbatim}
Note that we get the same optimal value, but no optimal solution is
provided.


\newpage
\pagestyle{myheadings} 
\markright{  \rm \normalsize BIBLIOGRAPHY. \hspace{0.5cm}
}

\end{document}